\documentclass[final,leqno,onefignum,onetabnum]{siamltex}

\usepackage{fullpage}
\usepackage{upgreek}
\usepackage[mathscr]{eucal}
\usepackage{color}
\usepackage[usenames,dvipsnames,svgnames,table]{xcolor}
\usepackage{amsmath,amssymb,amsopn,mathtools}
\usepackage{graphicx}
\usepackage{hyperref}       
\hypersetup{
	colorlinks=true,   
	citecolor=blue,     
	filecolor=blue,     
	linkcolor=red,    
	urlcolor=blue}      

\usepackage{relsize}
\numberwithin{equation}{section}
\usepackage{enumitem}

\newlist{Assumption}{enumerate}{1}
\setlist[Assumption]{label=A\arabic*}

\definecolor{Blue}{rgb}{0,0,1}
\definecolor{Red}{rgb}{1,0,0}
\definecolor{Green}{rgb}{0,1,0}

\newcommand{\fontDiscrete}{\mathcal}

\newcommand{\kmax}{k_\text{max}}
\newcommand{\bds}[1]{{\boldsymbol{#1}}}

\newcommand{\defeq}{:=}

\newcommand{\optDummy}{v}

\newcommand{\DEIM}{\text{DEIM}}
\newcommand{\SNS}{\text{SNS}}
\newcommand{\GNAT}{\text{GNAT}}
\newcommand{\STGNAT}{\text{ST-GNAT}}
\newcommand{\methodAcronym}{\text{SNS}}

\newcommand{\conditionnumber}{\kappa}

\newcommand{\temperature}{u}

\newcommand{\stateContinuousDummyEntryNo}{\boldsymbol{u}}

\newcommand{\stateContinuousDummy}[1]{\bds{\stateContinuousDummyEntryNo}} 
\newcommand{\stateTwoContinuousDummyEntryNo}{\boldsymbol{v}}

\newcommand{\stateTwoContinuousDummy}[1]{\bds{\stateTwoContinuousDummyEntryNo}}

\newcommand{\crossareaSymb}{A}
\newcommand{\externalforceSymb}{q}
\newcommand{\externalforce}{\boldsymbol \externalforceSymb}

\newcommand{\projection}{\mathbb{P}}
\newcommand{\projections}{\projection_\star}

\newcommand{\ndof}{\nspacedof}

\newcommand{\unitvec}{\bds{e}}
\newcommand{\unitvecArg}[1]{\unitvec_{#1}}

\newcommand{\zero}{\bds{0}}
\newcommand{\argmin}[1]{\underset{#1}{\text{argmin}}}
\newcommand{\range}[1]{\text{range}(#1)}

\newcommand{\orthobasis}{\boldsymbol \Upsilon}
\newcommand{\orthogonalmat}{\boldsymbol{Q}}
\newcommand{\orthogonalmatst}{\bar{\orthogonalmat}}
\newcommand{\righttrianglemat}{\boldsymbol R}

\newcommand{\ndofPorts}[1]{n^p}

\newcommand{\card}[1]{|#1|}

\newcounter{remctr}
\newcounter{propctr}
\newcounter{proposctr}
\newcommand{\RR}[1]{\ensuremath{\mathbb{R}^{ #1 }}}
\newcommand{\NN}{\mathbb{N}}
\newcommand{\RRplus}[1]{\ensuremath{\mathbb{R}_+^{ #1 }}}
\newcommand{\natNo}{\NN}
\newcommand{\nat}[1]{\natNo(#1)}

\newcommand{\innat}[1]{\in\nat{#1}}

\newcommand{\Span}[1]{\mathrm{span}\{#1\}}
\newcommand{\FE}{\mathrm{FE}}
\newcommand{\BE}{\mathrm{BE}}
\newcommand{\BDF}{\mathrm{BDF}}
\newcommand{\ADM}{\mathrm{AM}}
\newcommand{\ADB}{\mathrm{AB}}
\renewcommand{\Range}[1]{\mathrm{Ran}\left(#1\right)}

\newcommand{\identity}[1]{\boldsymbol I_{#1}}

\newcommand{\mass}{\Mbold}
\newcommand{\nonsingularmat}{\Kbold}
\newcommand{\pexit}{P_\text{exit}}
\newcommand{\specificheat}{\gamma}
\newcommand{\energypermass}{\epsilon}
\newcommand{\densitySymb}{\rho}
\newcommand{\velocitySymb}{u}
\newcommand{\energydensity}{e}
\newcommand{\pressureSymb}{p}
\newcommand{\specificgasconstant}{R}
\newcommand{\speedofsoundSymb}{c}
\newcommand{\machSymb}{M}
\newcommand{\temperatureSymb}{T}
\newcommand{\resSymb}{r}
\newcommand{\resRedSymb}{{\mathsf r}}
\newcommand{\fluxSymb}{f}
\newcommand{\paramSymb}{\mu}
\newcommand{\solSymb}{x}
\newcommand{\obliqueprojector}{\mathcal P}

\newcommand{\spaceSymb}{s}
\newcommand{\spaceSymbExt}{e}
\newcommand{\spaceSymbExtra}{E}
\newcommand{\timeSymb}{t}

\newcommand{\weightmatSymb}{A}

\newcommand{\basismatspaceSymb}{\Phi}
\newcommand{\basisvecspaceSymb}{\phi}

\newcommand{\samplematSymb}{Z}
\newcommand{\samplevecSymb}{z}
\newcommand{\dummySymb}{v}
\newcommand{\snapshotSymb}{X}

\newcommand{\nsmall}{n}
\newcommand{\nsample}{s}
\newcommand{\nparam}{\nsmall_\mu}
\newcommand{\ntrain}{\nsmall_\text{train}}

\newcommand{\nrestrain}{\nsmall_\text{res}}
\newcommand{\nressample}{\nsmall_{\samplevecSymb}}

\newcommand{\nbasisflux}{\nsmall_\fluxSymb}
\newcommand{\nbasisres}{\nsmall_\resSymb}

\newcommand{\nbasisspace}{{\nsmall_\spaceSymb}}
\newcommand{\nbasisspaceext}{{\nsmall_\spaceSymbExt}}
\newcommand{\nbasisspaceextra}{{\nsmall_\spaceSymbExtra}}

\newcommand{\nbasisst}{{\nsmall_{\spaceSymb\timeSymb}}}

\newcommand{\nbig}{N}
\newcommand{\nspacedof}{\nbig_\spaceSymb}
\newcommand{\ntimedof}{{\nbig_\timeSymb}}

\newcommand{\paramDomain}{\mathcal D}
\newcommand{\paramDomainTrain}{\paramDomain_\text{train}}

\newcommand{\ones}{\onebold}

\newcommand{\totaltime}{T}

\newcommand{\fluxSubspace}{\fontDiscrete F}
\newcommand{\spatialSubspace}{\fontDiscrete S}

\newcommand{\spatiotemporalSubspace}{\fontDiscrete{ST}}

\newcommand{\Abe}{\boldsymbol A_\text{BE}}

\newcommand{\resSymbst}{\bar{\resSymb}}
\newcommand{\resRedSymbst}{\bar{\resRedSymb}}
\newcommand{\res}{\boldsymbol \resSymb}

\newcommand{\resArg}[1]{\boldsymbol \resSymb^{#1}}

\newcommand{\resRedApprox}{\tilde{\boldsymbol \resRedSymb}}
\newcommand{\resRedApproxArg}[1]{\resRedApprox^{#1}}

\newcommand{\resn}{\resArg{n}}
\newcommand{\resst}{\boldsymbol \resSymbst}

\newcommand{\resRedst}{\boldsymbol \resRedSymbst}

\newcommand{\param}{\boldsymbol \paramSymb}

\newcommand{\paramDummy}{\boldsymbol \nu}
\newcommand{\dt}{\Delta \timeSymb}

\newcommand{\sol}{\boldsymbol \solSymb}
\newcommand{\solw}{\boldsymbol w}

\newcommand{\solst}{\bar{\sol}}

\newcommand{\flux}{\boldsymbol \fluxSymb}
\newcommand{\fluxArg}[1]{\flux_{#1}}

\newcommand{\fluxst}{\bar{\flux}}

\newcommand{\solDummy}{\boldsymbol w}

\newcommand{\solRedDummyOpt}{\hat {\boldsymbol \optDummy}}
\newcommand{\timeArg}[1]{t^{#1}}
\newcommand{\timeDummy}{\tau}
\newcommand{\timeDomain}{[0,\totaltime]}

\newcommand{\solArg}[1]{\sol_{#1}}
\newcommand{\solFunc}{\sol}

\newcommand{\solFuncArg}[1]{\solFunc(t^{#1};\param)}

\newcommand{\spacetimesubspace}{\fontDiscrete{ST}}

\newcommand{\solapproxFunc}{\solapprox}
\newcommand{\solapproxFuncArg}[1]{\solapproxFunc(t^{#1};\param)}

\newcommand{\solapproxFuncOnlyParam}{\solapproxFunc(\cdot;\param)}

\newcommand{\redsolapproxFunc}{\redsolapprox}
\newcommand{\redsolapproxFuncArg}[1]{\redsolapproxFunc(t^{#1};\param)}

\newcommand{\testbasisArg}[2]{\boldsymbol \xi_{ij}}

\newcommand{\solapproxArg}[1]{\solapprox_{#1}}

\newcommand{\solSTscalar}{y}
\newcommand{\solST}{\boldsymbol{\solSTscalar}}

\newcommand{\solapproxST}{\tilde\solST}

\newcommand{\redsolapproxST}{\hat\solST}

\newcommand{\redsolapproxSTArg}[1]{\hat\solSTscalar_{#1}(\param)}

\newcommand{\stbasismat}{\bar{\boldsymbol \Phi}}
\newcommand{\stbasismatext}{\bar{\boldsymbol \Phi}_e}
\newcommand{\stbasismatextra}{\bar{\boldsymbol \Phi}_E}
\newcommand{\stbasisvec}[1]{\bar{\boldsymbol \phi}_{#1}}

\newcommand{\mseError}{\text{relative error}}

\newcommand{\solapprox}{\tilde\sol}
\newcommand{\solapproxfom}{\solapprox_\star}
\newcommand{\redfluxapprox}{\hat\flux}
\newcommand{\redsolapprox}{\hat\sol}
\newcommand{\redsolapproxArg}[1]{\redsolapprox_{#1}}

\newcommand{\dummy}{\boldsymbol {\dummySymb}}
\newcommand{\reddummy}{\hat{\dummy}}

\newcommand{\redres}{\hat{\res}}
\newcommand{\redresi}[1]{\hat{\res}_{#1}}
\newcommand{\solinit}{\sol^0}

\newcommand{\solinitst}{\bar{\sol}^0}
\newcommand{\extrast}{\bar{\boldsymbol{q}}^0}
\newcommand{\solref}{\solinit(\param)}

\newcommand{\basismatspace}{\boldsymbol{\basismatspaceSymb}}
\newcommand{\basismatspaceext}{\basismatspace_{\spaceSymbExt}}
\newcommand{\basismatspaceextra}{\basismatspace_{\spaceSymbExtra}}

\newcommand{\basisvecspace}{\boldsymbol{\basisvecspaceSymb}}

\newcommand{\basismatres}{\basismatspace_\resSymb}
\newcommand{\basismatrest}{{\bar \basismatspace}_\resSymb}
\newcommand{\basismatrestvec}[1]{{\bar {\boldsymbol{\basisvecspaceSymb}}_{r,#1}}}

\newcommand{\basisfluxArg}[1]{\boldsymbol {\basisvecspaceSymb}_{\fluxSymb,{#1}}}
\newcommand{\basisfluxmat}{\boldsymbol \basismatspaceSymb_\fluxSymb}

\newcommand{\basisresvecArg}[1]{\boldsymbol \phi_{r,#1}}

\newcommand{\weightmat}{\boldsymbol{\weightmatSymb}}

\newcommand{\snapshots}{\boldsymbol \snapshotSymb}

\newcommand{\tuningparam}{\eta}

\newcommand{\samplematNT}{\boldsymbol \samplematSymb}
\newcommand{\samplemat}{\samplematNT^T}
\newcommand{\samplematstNT}{\bar{\samplematNT}}
\newcommand{\samplematst}{\samplematstNT^T}

\newcommand{\bmat}[1]{\begin{bmatrix}#1\end{bmatrix}} 
\newcommand{\pmat}[1]{\begin{pmatrix}#1\end{pmatrix}} 

\def\onebold{\boldsymbol{1}}
\def\Abold{\boldsymbol{A}}
\def\Bbold{\boldsymbol{B}}

\def\Kbold{\boldsymbol{K}}

\def\Mbold{\boldsymbol{M}}

\def\Sbold{\boldsymbol{S}}

\def\Ubold{\boldsymbol{U}}
\def\Vbold{\boldsymbol{V}}
\def\Wbold{\boldsymbol{W}}

\def\Sigmabold{\boldsymbol{\Sigma}}

\def\zerobold{{\bf 0}}

\newtheorem{remark}{Remark}[section]

\usepackage{subfigure, algorithmic, algorithm}
\usepackage{amsmath,amssymb}
\usepackage{epsfig, algorithmic, algorithm}
\usepackage{multirow, booktabs}
\usepackage{siunitx}
\usepackage{xr}
\mathtoolsset{showonlyrefs=true}
\usepackage{xr}

\title{SNS: a solution-based nonlinear subspace method for
time-dependent model order reduction}

\author{Youngsoo Choi\thanks{Computational
    Engineering Division, Lawrence Livermore National Laboratory, Livermore,
  CA 94550 (choi15@llnl.gov)}\and
        Deshawn Coombs\thanks{Department of Mechanical and Aerospace
        Engineering, Syracuse University, Syracuse, NY 13244
      (dmcoombs@syr.edu)}
        \and
        Robert Anderson\thanks{Center for Applied Scientific Computing,
        Lawrence Livermore National Laboratory,
      Livermore, CA 94550 (anderson110@llnl.gov)}
}

\begin{document}
\setlength{\abovedisplayskip}{3pt}
\setlength{\belowdisplayskip}{3pt} 
\setlength{\abovedisplayshortskip}{3pt} 
\setlength{\belowdisplayshortskip}{3pt}

\maketitle

\begin{abstract}
  Several reduced order models have been successfully developed for nonlinear
  dynamical systems. To achieve a considerable speed-up, a hyper-reduction step
  is needed to reduce the computational complexity due to nonlinear terms.  Many
  hyper-reduction techniques require the construction of nonlinear term basis,
  which introduces a computationally expensive offline phase.  A novel way of
  constructing nonlinear term basis within the hyper-reduction process is
  introduced.  In contrast to the traditional hyper-reduction techniques where
  the collection of nonlinear term snapshots is required, the SNS method
  avoids collecting the nonlinear term snapshots. Instead, it uses
  the solution snapshots that are used for building a solution basis,
  which enables avoiding an extra data compression of nonlinear term snapshots.
  As a result, the SNS method provides a more efficient offline strategy than
  the traditional model order reduction techniques, such as the DEIM, GNAT, and
  ST-GNAT methods.  
  The SNS method is theoretically justified by the conforming subspace
  condition and the subspace inclusion relation.  
  It is useful for model order reduction of large-scale nonlinear
  dynamical problems to reduce the offline cost. 
  It is especially useful for ST-GNAT that has shown promising results, such as
  a good accuracy with a considerable online speed-up for hyperbolic problems
  in a recent paper by Choi and Carlberg in SISC 2019,
  because ST-GNAT involves an expensive
  offline cost related to collecting nonlinear term snapshots. 
  Error analysis shows that the oblique projection error bound of the SNS method
  depends on the condition number of the matrix $\mass$ (e.g., a volume matrix
  generated from a discretization of a specific numerical scheme).   
  Numerical results support that
  the accuracy of the solution from the SNS method is comparable to the
  traditional methods and a considerable speed-up (i.e., a factor of two to a hundred) 
  is achieved in the offline
  phase.
\end{abstract}




\section{Introduction}\label{sec:intro}
    Time-dependent nonlinear problems arise in many important disciplines such
    as engineering, science, and technologies.   They are numerically solved if
    it is not possible to solve them analytically. Depending on the complexity
    and size of the governing equations and problem domains,  the problems can
    be computationally expensive to solve. It may take a long time to run one
    forward simulation even with high performance computing.   For example, a
    simulation of the powder bed fusion additive manufacturing procedure shown
    in \cite{khairallah2014mesoscopic} takes a week to finish with $108$ cores.
    Other computationally expensive simulations include the 3D shocked spherical
    Helium bubble simulation appeared in \cite{anderson2018high} and the
    inertial confinement fusion implosion dynamics simulations appeared in
    \cite{blast-library}.
    The computationally expensive simulations
    are not desirable in the context of parameter study, design optimization,
    uncertainty quantification, and inverse problems where several forward
    simulations are needed. A Reduced Order Model (ROM) can be useful in this
    context to accelerate the computationally expensive forward simulations
    with good enough approximate solutions.

    We consider projection-based ROMs for nonlinear dynamical systems.  Such
    ROMs include the Empirical Interpolation method (EIM)
    \cite{barrault2004empirical,grepl2007efficient}, the Discrete Empirical
    Interpolation Method (DEIM) \cite{chaturantabut2010nonlinear,drmac2016new}
    and the Gauss--Newton with Approximated Tensors
    \cite{carlberg2011efficient,carlberg2013gnat}, the Best Point Interpolation Method (BPIM)
    \cite{nguyen2008efficient}, the Missing Point Estimation (MPE)
    \cite{astrid2008missing}, and Cubature Methods
    \cite{an2008optimizing,farhat2014dimensional,farhat2015structure,hernandez2017dimensional}.
    There a hyper-reduction (The term first used in the context of ROMs by
    Ryckelynck in \cite{ryckelynck2005priori}) is necessary to efficiently reduce
    the complexity of nonlinear terms for a considerable speed-up compared to a
    corresponding full order model. The DEIM and GNAT approaches take discrete
    nonlinear term snapshots to build a nonlinear term basis. Then, they select
    a subset of each nonlinear term basis vector to either interpolate or
    data-fit in a least-squares sense. In this way, they reduce the
    computational complexity of updating nonlinear terms in an iterative solver
    for nonlinear problems.  The EIM, BPIM, and MPE approaches take the similar
    hyper-reduction to the ones in DEIM and GNAT except for the fact that they
    build nonlinear term basis functions in a continuous space.  Whether they
    work on a discrete or continuous space, they follow the framework of
    reconstructing ``gappy" data, first introduced in the context of reduced
    order models by \cite{everson1995karhunen}. The cubature methods, in
    contrast, take a different approach. They approximate nonlinear integral as
    a finite sum of positive scalar weights multiplied by the integrand
    evaluated at sampled elements. The cubature methods developed in
    \cite{an2008optimizing,farhat2014dimensional,farhat2015structure} do not require
    building nonlinear term basis. They solve the Non-Negative Least-Squares
    (NNLS) problem directly with the nonlinear term snapshots.  Recently,
    Hernandez, et al., in \cite{hernandez2017dimensional} developed the Empirical
    Cubature Method (ECM) approach that builds nonlinear term basis to solve a
    smaller NNLS problem.  
    The requirement of building nonlinear term basis in hyper-reduction results
    in computational cost and storage in addition to solution basis
    construction. The cost is significant if the corresponding Full Order Models
    (FOMs) are large-scale. The large-scale problem requires large additional
    storage for nonlinear term snapshots and large-scale compression techniques
    to build a basis.  In particular, the cost of the hyper-reduction in the
    recently-developed space--time ROM (i.e., ST-GNAT) \cite{choi2019space} is
    even more significant than aforementioned spatial ROMs (e.g., EIM, DEIM,
    GNAT, BPIM, MPE, and ECM).  The best nonlinear term snapshots of the ST-GNAT
    method are obtained from the corresponding space--time ROMs without
    hyper-reduction (i.e., ST-LSPG) \cite{choi2019space}, which is not practical
    for a large-scale problem due to the cumbersome size.\footnote{The full
    size, implicitly handled by the ST-LSPG approach due to the nonlinear term
    and corresponding Jacobian updates, is the FOM spatial degrees of freedom
    multiplied by the number of time steps.} Therefore, a novel and efficient
    way of constructing nonlinear term basis needs to be developed.

    This paper shows a practical way of avoiding nonlinear term snapshots for
    the construction of nonlinear term basis.  The idea comes from a simple fact
    that the nonlinear terms are related with solution snapshots through
    underlying time integrators.  In fact, many time integrators approximate
    time derivative terms as a linear combination of the solution snapshots. It
    implies that the nonlinear term snapshots belong to the subspace spanned by
    the solution snapshots.  Furthermore, a subspace needed for nonlinear terms
    in a hyper-reduction is determined by the range space of the solution basis
    matrix possibly multiplied by a nonsingular matrix (e.g., the volume matrix).  
    Therefore, the solution snapshots can be used to construct
    nonlinear term basis. This leads to our proposed method, the Solution-based
    Nonlinear Subspace (SNS) method, that provides  two savings for constructing
    nonlinear term basis because of 
      \begin{enumerate}
        \item no additional collection of nonlinear term snapshots (i.e., storage saving).
        \item no additional compression of snapshots (e.g., no additional singular
          value decomposition of nonlinear term snapshots, implying
          computational cost saving).
      \end{enumerate}
    The first saving is especially big for GNAT and ST-GNAT becuase they involve
    expensive collection procedure for their best performace. 
\subsection{Organization of the paper}\label{sec:organization} 
    We start our discussion by describing the time-continuous representation of
    the FOM in Section~\ref{sec:FOM}.  Section~\ref{sec:FOM} also describes the
    time-discrete representation of the FOM with one-step Euler time
    integrators.  The subspace inclusion relation between the subspaces spanned
    by the solution and nonlinear term snapshots is described for the Euler time
    integrators.  Several projection-based ROMs (i.e., the DEIM, GNAT, and
    ST-GNAT models) are considered in Section~\ref{sec:ROMs} for the SNS method
    to be applied.  If a reader is familiar with DEIM, GNAT, and ST-GNAT,
    then the one can skip Section~\ref{sec:FOM} and \ref{sec:ROMs} and start
    with Section~\ref{sec:SNS}.
    The SNS method is described in Section~\ref{sec:SNS} and
    applied to those ROMs. Section~\ref{sec:SNS} also introduces the conforming
    subspace condition to justify the SNS method.
    Section~\ref{sec:erroranalysis} shows an error analysis for the SNS method
    regarding the oblique projection error bound.  It analyzes the effect of the
    nonsingular matrix that is used to form a nonlinear term basis.
    Section~\ref{sec:numericresults} reports numerical results that support
    benefits of the SNS method and Section~\ref{sec:conclusion} concludes the
    paper with summary and future work. Although the SNS method is mainly
    illustrated with the Euler time integrators throughout the paper, it is
    applicable to other time integrators.  Appendix~\ref{sec:appendix} considers
    several other time integrators and the subspace inclusion relation for each
    time integrator.  The following time integrators are included: the
    Adams--Bashforth, Adams--Moulton, BDF, and midpoint Runge--Kutta methods.  

\section{Full order model}\label{sec:FOM}
    A parameterized nonlinear dynamical system is considered, characterized by a system of nonlinear ordinary differential equations (ODEs), which can be considered as a resultant system from semi-discretization of Partial Differential Equations (PDEs) in space domains
 \begin{equation} \label{eq:fom}
  \mass\frac{d\sol}{dt} = \flux(\sol,t; \param),\quad\quad
  \sol(0;\param) = \solinit(\param),
 \end{equation} 
 where $\mass \in \RR{\nspacedof \times \nspacedof}$ denotes a nonsingular matrix, $t\in[0,\totaltime]$ denotes time with the final time $\totaltime\in\RRplus{}$, and $\sol(t;\param)$ denotes the time-dependent, parameterized state implicitly defined as the solution to problem~\eqref{eq:fom} with $\sol:\timeDomain\times \paramDomain\rightarrow \RR{\nspacedof}$.  Further, $\flux: \RR{\nspacedof} \times [0,\totaltime] \times \paramDomain \rightarrow \RR{\nspacedof}$ with $(\solDummy,\timeDummy;\paramDummy)\mapsto\flux(\solDummy,\timeDummy;\paramDummy) $ denotes the scaled velocity of $\mass\sol$, which we assume to be nonlinear in at least its first argument.  The initial state is denoted by $\solinit:\paramDomain\rightarrow \RR{\nspacedof}$, and $\param \in \paramDomain$ denotes the parameters with parameter domain $\paramDomain\subseteq\RR{\nparam}$. 

    A uniform time discretization is assumed throughout the paper, characterized
    by time step $\dt\in\RRplus{}$ and time instances $\timeArg{n} =
    \timeArg{n-1} + \dt$ for $n\innat{\ntimedof}$ with $\timeArg{0} = 0$,
    $\ntimedof\in\natNo$, and $\nat{N}\defeq\{1,\ldots,N\}$.  To avoid
    notational clutter, we introduce the following time discretization-related
    notations: $\solArg{n} \defeq \solFuncArg{n}$, $\solapproxArg{n} \defeq
    \solapproxFuncArg{n}$, $\redsolapproxArg{n} \defeq \redsolapproxFuncArg{n}$,
    and $\fluxArg{n} \defeq \flux(\solFuncArg{n},t^{n}; \param)$. 

  Two different types of time discretization methods are considered: explicit
  and implicit time integrators. As an illustration purpose, we mainly consider
  the forward Euler time integrator for an explicit scheme and the backward
  Euler time integrator for an implicit scheme.  Several other time integrators
  are shown in Appendix~\ref{sec:appendix}. 

    The explicit Forward Euler (FE) method numerically solves Eq.~\eqref{eq:fom}, by time-marching with the following update: 
      \begin{equation} \label{eq:forwardEuler}
        \mass \solArg{n} - \mass \solArg{n-1} =\dt\fluxArg{n-1}.
      \end{equation}
      Eq.~\eqref{eq:forwardEuler} implies the following subspace inclusion:
      \begin{equation}\label{eq:forwardEuler_partial_inclusion}
        \Span{\fluxArg{n-1}} \subseteq \Span{\mass \solArg{n-1},\mass \solArg{n}}.    
      \end{equation}
      By induction, we conclude the following subspace inclusion relation:
      \begin{equation}\label{eq:forwardEuler_total_inclusion}
        \Span{\fluxArg{0},\dots,\fluxArg{\ntimedof-1}} \subseteq \Span{\mass \solArg{0},\ldots,\mass
        \solArg{\ntimedof}},
      \end{equation}
      which shows that the span of nonlinear term snapshots is included in the span of $\mass$-scaled solution snapshots.  The residual function with the forward Euler time integrator is defined as
      \begin{align}\label{eq:residual_FE} 
      \begin{split}
        \resn_{\FE}(\solArg{n};\solArg{n-1},\param) &\defeq 
        \mass(\solArg{n} - \solArg{n-1}) -\dt\fluxArg{n-1} 
      \end{split}
      \end{align} 

    The implicit Backward Euler (BE) method numerically solves Eq.~\eqref{eq:fom}, by solving the following nonlinear system of equations for $\solArg{n}$ at $n$-th time step:
      \begin{equation} \label{eq:backwardEuler}
        \mass \solArg{n} - \mass \solArg{n-1} = \dt\fluxArg{n}.
      \end{equation} 
      Eq.~\eqref{eq:backwardEuler} implies the following subspace inclusion:
      \begin{equation}\label{eq:backwardEuler_partial_inclusion}
        \Span{\fluxArg{n}} \subseteq \Span{\mass \solArg{n-1},\mass \solArg{n}}.    
      \end{equation}
      By induction, we conclude the following subspace inclusion relation:
      \begin{equation}\label{eq:backwardEuler_total_inclusion}
        \Span{\fluxArg{1},\dots,\fluxArg{\ntimedof}} \subseteq \Span{\mass \solArg{0},\ldots,\mass
        \solArg{\ntimedof}},
      \end{equation}
      which shows that the span of nonlinear term snapshots is included in the span of $\mass$-scaled solution snapshots.  The residual function with the backward Euler time integrator is defined as
      \begin{align}\label{eq:residual_BE} 
      \begin{split}
        \resn_{\BE}(\solArg{n};\solArg{n-1},\param) &\defeq 
        \mass(\solArg{n} - \solArg{n-1}) -\dt\fluxArg{n}.
      \end{split}
      \end{align}

  \section{Projection-based reduced order models}\label{sec:ROMs}
    Projection-based reduced order models are considered for nonlinear dynamical systems.  Especially, the ones that require building a nonlinear term basis are of our interest: the DEIM, GNAT, and ST-GNAT approaches.
    \subsection{DEIM}\label{sec:DEIM}
The DEIM approach applies spatial projection using a subspace $\spatialSubspace \defeq \Span{\basisvecspace_i}_{i=1}^\nbasisspace \subseteq \RR{\nspacedof}$ with $\dim(\spatialSubspace)=\nbasisspace\ll\nspacedof$.  Using this subspace, it approximates the solution as $ \sol\approx\solapprox \in\solArg{0}+\spatialSubspace $ (i.e., in a trial subspace) or equivalently 
\begin{equation}\label{eq:spatialROMsolution} 
  \solapprox= \solArg{0} + \basismatspace\redsolapprox 
\end{equation} 
  where $\basismatspace \defeq [\basisvecspace_1 \cdots \basisvecspace_{\nbasisspace}] \in\RR{\nspacedof\times\nbasisspace}$ denotes a basis matrix and $\redsolapprox\in\RR{\nbasisspace}$ with $\redsolapproxArg{0}=\zero$ denotes the generalized coordinates.  Replacing $\sol$ with $\solapprox$ in Eq.~\eqref{eq:fom} leads to the following system of equations with reduced number of unknowns:
 \begin{equation} \label{eq:rom}
  \mass\basismatspace\frac{d\redsolapprox}{dt} =
  \flux(\solArg{0}+\basismatspace\redsolapprox,t; \param).
 \end{equation}
 For constructing $\basismatspace$, Proper Orthogonal Decomposition (POD) is
 commonly used.  POD \cite{berkooz1993proper} obtains
 $\basismatspace$ from a truncated Singular Value Decomposition (SVD)
 approximation to a FOM solution snapshot matrix. It is related to  principal
 component analysis in statistical analysis \cite{hotelling1933analysis} and 
 Karhunen--Lo\`{e}ve expansion \cite{loeve1955} in stochastic analysis.
 POD forms a solution snapshot matrix,
 $\snapshots\defeq\bmat{\solArg{0}^{\param_1} &
 \cdots &
 \solArg{\ntimedof}^{\param_{\nparam}}}\in\RR{\ndof\times\nparam(\ntimedof+1)}$,
 where $\solArg{n}^{\param_k}$ is a solution state at $n$th time step with
 parameter $\param_k$ for $n\innat{\ntimedof}$ and $k\innat{\nparam}$.
 Then, POD computes its thin SVD: 
 \begin{equation}\label{eq:SVD} 
   \snapshots = \Ubold\Sigmabold\Vbold^T,
 \end{equation} 
 where $\Ubold\in\RR{\ndof\times\ntimedof}$ and
 $\Vbold\in\RR{\ntimedof\times\ntimedof}$ are orthogonal matrices and
 $\Sigmabold\in\RR{\ntimedof\times\ntimedof}$ is a diagonal matrix with singular
 values on its diagonals.  Then POD chooses the leading $\nbasisspace$
 columns of $\Ubold$ to set $\basismatspace$ (i.e., $\basismatspace =
 \Ubold(:,1:\nbasisspace)$ in MATLAB notation).  The POD basis minimizes
 $\|\snapshots - \basismatspace\basismatspace^T\snapshots \|_F^2$ over all
 $\basismatspace\in\RR{\ndof\times\nbasisspace}$ with orthonormal columns, where
 $\|\Abold\|_F$ denotes the Frobenius norm of a matrix $\Abold\in\RR{I\times
 J}$, defined as $\|\Abold\|_F = \sqrt{\sum_{i=1}^{I}\sum_{j=1}^{J} a_{ij}}$
 with $a_{ij}$ being an $(i,j)$-th element of $\Abold$.
 Since the objective function does not change if $\basismatspace$ is
 post-multiplied by an arbitrary $\nbasisspace\times\nbasisspace$ orthogonal
 matrix, the POD procedure seeks the optimal $\nbasisspace$–-dimensional
 subspace that captures the snapshots in the least-squares sense.  For more
 details on POD, we refer to \cite{hinze2005proper,kunisch2002galerkin}.

 Note that Eq.~\eqref{eq:rom} has more equations than unknowns (i.e., an
 overdetermined system).  It is likely that there is no solution satisfying
 Eq.~\eqref{eq:rom}.  In order to close the system, the Galerkin projection
 solves the following reduced system with $\redsolapproxArg{0}=\zero$:
 \begin{equation} \label{eq:rom_galerkin}
  \basismatspace^T\mass\basismatspace\frac{d\redsolapprox}{dt} =
  \basismatspace^T\flux(\solArg{0}+\basismatspace\redsolapprox,t; \param).
 \end{equation}

 Applying a time integrator to Eq.~\eqref{eq:rom_galerkin} leads to a fully
 discretized reduced system, denoted as the reduced O$\Delta$E.  
 Note that the reduced O$\Delta$E has $\nbasisspace$ unknowns and
 $\nbasisspace$ equations. If an implicit time integrator is applied, a
 Newton--type method can be applied to solve for unknown generalized coordinates
 each time step. 
 If an explicit time integrator is applied, time marching
 updates will solve the system. 
 However, we cannot expect any speed-up because the size of the
 nonlinear term $\flux$ and its Jacobian, which need to be updated for every
 Newton step, scales with the FOM size.  Thus, to overcome this issue, the DEIM
 approach applies a hyper-reduction technique.  That is, it projects $\flux$
 onto a subspace $\fluxSubspace \defeq
 \Span{\basisfluxArg{i}}_{i=1}^{\nbasisflux}$ and approximates $\flux$ as   
 \begin{equation}\label{eq:DEIM_approx}
   \flux \approx \basisfluxmat\redfluxapprox,
 \end{equation}
 where $\basisfluxmat \defeq [\basisfluxArg{1},\ldots,\basisfluxArg{\nbasisflux}
 ] \in \RR{\ndof\times\nbasisflux}$, $\nbasisflux \ll \ndof$, denotes the nonlinear term basis matrix and $\redfluxapprox \in \RR{\nbasisflux}$ denotes the generalized coordinates of the nonlinear term. The generalized coordinates, $\redfluxapprox$, can be determined by the following interpolation:
 \begin{equation}\label{eq:DEIM_interpolation}
   \samplemat\flux = \samplemat\basisfluxmat\redfluxapprox,
 \end{equation}
 where
 $\samplemat\defeq[\unitvecArg{p_1},\ldots,\unitvecArg{p_{\nbasisflux}}]^T\in\RR{\nbasisflux\times\ndof}$
 is the sampling matrix and $\unitvecArg{p_i}$ is the $p_i$th column of the identity matrix $\identity{\ndof}\in\RR{\ndof\times\ndof}$.  Therefore, Eq.~\eqref{eq:DEIM_approx} becomes
 \begin{equation}\label{eq:DEIM_f}
   \flux \approx \obliqueprojector_{\DEIM} \flux,
 \end{equation}
 where $\obliqueprojector_{\DEIM} := \basisfluxmat
 (\samplemat\basisfluxmat)^{-1}\samplemat \in \RR{\nbasisflux\times\nbasisflux}$ 
 is the DEIM oblique projection matrix. 
 The DEIM approach does not construct the sampling matrix $\samplematNT$.
 Instead, it maintains the sampling indices
 $\{p_1,\ldots,p_{\nbasisflux}\}$ and
 corresponding rows of $\basisfluxmat$ and $\flux$. This enables DEIM to achieve
 a speed-up when it is applied to nonlinear problems. 

 The original DEIM paper \cite{chaturantabut2010nonlinear} constructs the
 nonlinear term basis $\{ \basisfluxArg{1},\ldots,\basisfluxArg{\nbasisflux} \}$
 by applying another POD on the nonlinear term snapshots\footnote{the nonlinear
 term snapshots are $\{ \fluxArg{1},\ldots,\fluxArg{\ntimedof}  \}$ with the
 backward Euler time integrator and $\{ \fluxArg{0},\ldots,\fluxArg{\ntimedof-1}
 \}$ with the forward Euler time integrator} obtained from the FOM simulation at
 every time step.  This implies that DEIM requires two separate SVDs and storage
 for two different snapshots (i.e., solution state and nonlinear term
 snapshots).  Section~\ref{sec:SNS} discusses how to avoid the collection of
 nonlinear term snapshots and an extra SVD without losing the quality of the
 hyper-reduction.

 The sampling indices (i.e., $\samplematNT$) can be found either by a row pivoted
 LU decomposition \cite{chaturantabut2010nonlinear} or the strong column pivoted
 rank-revealing QR (sRRQR)
 decomposition \cite{drmac2016new, drmac2018discrete}. Depending on the
 algorithm of selecting the sampling indices, the DEIM projection 
 error bound is determined. For example, the row pivoted LU decomposition in
 \cite{chaturantabut2010nonlinear} results in the following error bound:
 \begin{equation}\label{eq:errorbound_DEIM}
   \|\flux - \obliqueprojector \flux \|_2 \leq \conditionnumber
   \|(\identity{\ndof}-\basisfluxmat\basisfluxmat^T)\flux  \|_2,
 \end{equation}
 where $\|\cdot\|_2$ denotes $\ell_2$ norm of a vector and 
 $\conditionnumber$ is the condition number of
 $(\samplemat\basisfluxmat)^{-1}$ and it is bounded by
 \begin{equation}\label{eq:crudeBound}
   \conditionnumber \leq
   (1+\sqrt{2\ndof})^{\nbasisflux-1}\|\basisfluxArg{1}\|_\infty^{-1}.
 \end{equation}
 On the other hand, the sRRQR factorization in \cite{drmac2018discrete} 
 reveals tighter bound than \eqref{eq:crudeBound}:
 \begin{equation}\label{eq:tighterBound}
   \conditionnumber \leq
   \sqrt{1+\tuningparam^2\nbasisflux(\ndof-\nbasisflux)}
 \end{equation}
 where $\tuningparam$ is a tuning parameter in the sRRQR factorization
 (i.e., $f$ in Algorithm 4 of \cite{gu1996efficient}).

 \subsection{GNAT}\label{sec:GNAT}
    In contrast to DEIM, the GNAT method takes the Least-Squares
    Petrov--Galerkin (LSPG) approach.  The LSPG method projects a fully
    discretized solution space onto a trial subspace. That is, it discretizes
    Eq.~\eqref{eq:fom} in time domain and replaces $\solArg{n}$ with
    $\solapproxArg{n} \defeq \solArg{0}+\basismatspace\redsolapproxArg{n}$ for
    $n\innat{\ntimedof}$ in residual functions defined in Section~\ref{sec:FOM}
    and Appendix \ref{sec:appendix}. 
    Here, we consider only implicit time integrators because the LSPG projection is
    equivalent to the Galerkin projection when an explicit time integrator is
    used as shown in Section 5.1 in \cite{carlberg2017galerkin}.
    The residual functions for implicit time
    integrators are defined in \eqref{eq:residual_BE}, \eqref{eq:residual_AM2}, 
    and \eqref{eq:residual_BDF2} for various time integrators.  
    For example, the residual function with the backward Euler time integrator\footnote{Although the backward Euler time integrator is used extensively in the paper as an illustration purpose, 
    many other time integrators introduced in Appendix~\ref{sec:appendix} can be
    applied to all the ROM methods dealt in the paper 
    in a straight forward way.} after the trial subspace projection becomes 
    \begin{align}\label{eq:trialsub_residual_BE} 
    \begin{split}
      \resRedApproxArg{n}_{\BE}(\redsolapproxArg{n};\redsolapproxArg{n-1},\param)
&\defeq 
      \resn_{\BE}(\solArg{0}+\basismatspace\redsolapproxArg{n};
      \solArg{0}+\basismatspace\redsolapproxArg{n-1},\param) 
\\ &= \mass\basismatspace(\redsolapproxArg{n} - \redsolapproxArg{n-1}) -\dt\flux(\solArg{0}+\basismatspace\redsolapproxArg{n},t;\param).
    \end{split}
    \end{align}
 The basis matrix $\basismatspace$ can be found by the POD as in the DEIM
 approach.  Note that Eq.~\eqref{eq:trialsub_residual_BE} is an over-determined
 system.  To close the system and solve for the unknown generalized coordinates,
 $\redsolapproxArg{n}$, the LSPG ROM takes the squared norm of the residual vector
 function and minimize it at every time step:
\begin{align} \label{eq:spOpt1}
  \begin{split} 
    \redsolapproxArg{n} = \argmin{\reddummy\in\RR{\nbasisspace}} \quad&
    \frac{1}{2} \left \|\resn_{\BE}(\reddummy;\redsolapproxArg{n-1},\param)
    \right \|_2^2.
  \end{split} 
\end{align}
The Gauss--Newton method with the starting point $\redsolapproxArg{n-1}$ is
applied to solve the minimization problem~\eqref{eq:spOpt1} in GNAT.  However,
as in the DEIM approach, a hyper-reduction is required for a speed-up due to the
presence of the nonlinear residual vector function.  The GNAT method
approximates the nonlinear residual term with gappy POD
\cite{everson1995karhunen}, whose procedure is similar to DEIM, as 
 \begin{equation}\label{eq:GNAT_approx}
   \res \approx \basismatres\redres,
 \end{equation}
 where $\basismatres \defeq
 [\basisresvecArg{1},\ldots,\basisresvecArg{\nbasisres} ] \in
 \RR{\ndof\times\nbasisres}$, $\nbasisspace \leq \nbasisres \ll \ndof$, denotes
 the residual basis
 matrix and $\redres \in \RR{\nbasisres}$ denotes the generalized coordinates of
 the nonlinear residual term.  In contrast to DEIM, the GNAT method solves
 the following least-squares problem to obtain the generalized coordinates
 $\redresi{n}$ at $n$-th time step:
\begin{align} \label{eq:GNAT_least-squares}
  \begin{split} 
    \redresi{n} = \argmin{\reddummy\in\RR{\nbasisres}} \quad&
    \frac{1}{2} \left \|\samplemat(\res - \basismatres\reddummy)
    \right \|_2^2.
  \end{split} 
\end{align}
 where $\samplemat\defeq[\unitvecArg{p_1},\ldots,\unitvecArg{p_{\nressample}}]^T
 \in\RR{\nressample\times\ndof}$, $\nbasisspace \leq \nbasisres \leq \nressample
 \ll \ndof$, is the sampling matrix and $\unitvecArg{p_i}$ is
 the $p_i$th column of the identity matrix
 $\identity{\ndof}\in\RR{\ndof\times\ndof}$.  The solution to
 Eq.~\eqref{eq:GNAT_least-squares} is given as
 \begin{equation}\label{eq:GNAT-generalizedcoordinates}
   \redresi{n} = (\samplemat\basismatres)^\dagger\samplemat\res,
 \end{equation}
 where the Moore--Penrose inverse of a matrix $\weightmat \in \RR{\nressample
 \times \nbasisres}$ with full column rank is defined as $\weightmat^{\dagger} := (\weightmat^T\weightmat)^{-1}\weightmat^T$. Therefore, Eq.~\eqref{eq:GNAT_approx} becomes
 \begin{equation}\label{eq:GNAT_r}
   \res \approx \obliqueprojector_{\GNAT} \res,
 \end{equation}
 where $\obliqueprojector_{\GNAT}:=
 \basismatres (\samplemat\basismatres)^\dagger\samplemat$ is the GNAT oblique
 projection matrix. Note that it has a similar structure to
 $\obliqueprojector_{\DEIM}$. The GNAT projection matrix has a pseudo-inverse
 instead of the inverse becuase it allows $\nbasisres<\nressample$.
 The GNAT method does not construct the sampling matrix $\samplematNT$.  Instead,
 it maintains the sampling indices $\{p_1,\ldots,p_{\nbasisflux}\}$ and
 corresponding rows of $\basismatres$ and $\res$. This enables GNAT to achieve a
 speed-up when it is applied to nonlinear problems.

 The sampling indices (i.e., $\samplematNT$) can be determined by Algorithm 3 of
 \cite{carlberg2013gnat} for computational fluid dynamics problems and Algorithm
 5 of \cite{carlberg2011efficient} for other problems. 
 These two algorithms take greedy procedure to minimize the error 
 in the gappy reconstruction of the POD basis vectors $\basismatres$. 
 The major difference between these sampling algorithms and the
 ones for the DEIM method is that these algorithms for the GNAT method allows
 oversampling (i.e., $\nressample > \nbasisres$), 
 resulting in solving least-squares problems in
 the greedy procedure. These selection algorithms can be viewed as the extension
 of Algorithm 1 in \cite{chaturantabut2010nonlinear} (i.e., a row pivoted LU
 decomposition) to the oversampling case. The nonlinaer residual term projection
 error associated with these sampling algorithms is presented in Appendix D of
 \cite{carlberg2013gnat}. That is,
 \begin{equation}\label{eq:traditionalGNATprojection_error}
   \|\res - \obliqueprojector_{\GNAT}\res\|_2
   \leq \| \righttrianglemat^{-1} \|_2 
   \|\res - \basismatres\basismatres^T\res \|_2
 \end{equation}
 where $\righttrianglemat$ is the triangular factor from the QR factorization of
 $\samplemat\basismatres$ 
 (i.e.,
 $\samplemat\basismatres=\orthogonalmat\righttrianglemat$).

 We present another sampling selection algorithm that was not considered in any
 GNAT papers (e.g., \cite{carlberg2011efficient,carlberg2013gnat}). It is to use
 the sRRQR factorization developed originally in \cite{gu1996efficient} and
 further utilized in the W-DEIM method of \cite{drmac2018discrete} for the case
 of $\nressample \geq \nbasisres$. That is, applying Algorithm 4 of
 \cite{gu1996efficient} to $\basismatres^T$ produces an index selection operator
 $\samplematNT$ whose projection error satisfies
 \begin{equation}\label{eq:betterGNATerror}
   \|\res - \obliqueprojector_{\GNAT}\res\|_2
   \leq \sqrt{1+\tuningparam^2\nbasisres(\ndof-\nbasisres)}
   \|\res - \basismatres\basismatres^T\res \|_2.
 \end{equation}
 This error bound can be obtained by setting the identity matrix as a weight
 matrix in Theorem 4.8 of \cite{drmac2018discrete}.

 Finally, the GNAT method minimizes the following least-squares problem at every
 time step, for example, with the backward Euler time integrator:
  \begin{align} \label{eq:GNAT-spOpt1}
    \begin{split} 
      \redsolapproxArg{n} = \argmin{\reddummy\in\RR{\nbasisspace}} \quad&
      \frac{1}{2} \left \|\ (\samplemat\basismatres)^\dagger\samplemat
      \resn_{\BE}(\reddummy;\redsolapproxArg{n-1},\param)
      \right \|_2^2.
    \end{split} 
  \end{align}
The GNAT method applies another POD to a nonlinear residual term snapshots to
construct $\basismatres$.  The original GNAT paper \cite{carlberg2013gnat}
collects residual snapshots at every Newton iteration from the LSPG
simulations\footnote{We denote the LSPG simulation to be the ROM simulation
without hyper-reduction. That is, LSPG solves Eq.~\eqref{eq:spOpt1} without
any hyper-reduction if the backward Euler time integrator is used.} for several reasons: 
\begin{enumerate}
  \item The GNAT method takes LSPG as a reference model (i.e., denoted as Model
    II in \cite{carlberg2013gnat}).  Its ultimate goal is to achieve the accuracy of Model II. 
  \item The residual snapshots taken every time step (i.e., at the end of Newton iterations at every time step) of the FOM are small in magnitude.
  \item The residual snapshots taken from every Newton step gives information
    about the path that the Newton iterations in Model II take. GNAT tries to
    mimic the Newton path that LSPG takes.
\end{enumerate}
Some disadvantages of the original GNAT approach include:
\begin{enumerate}
  \item The GNAT method requires more storage than DEIM to store residual snapshots from every Newton iteration (cf., The DEIM approach stores only one nonlinear term snapshot per each time step).
  \item The GNAT method requires more expensive SVD for nonlinear residual basis
    construction than the DEIM approach because the number of nonlinear residual
    snapshots in the GNAT method is larger than the number of nonlinear term snapshots in DEIM.
  \item The GNAT method requires the simulations of Model II which are
    computationally expensive. For a parametric global ROM, it is
    computationally expensive because it requires as many Model II simulations as there are training points in a given parameter domain.
\end{enumerate}
Section~\ref{sec:SNS} discusses how to avoid the collection of nonlinear term
snapshots and the extra SVD without losing the quality of the hyper-reduction.

\subsection{Space--time GNAT}\label{sec:STGNAT}
The ST-GNAT method takes the space--time LSPG (ST-LSPG) approach.  Given a time integrator, one can rewrite a fully discretized system of equations to Eq.~\eqref{eq:fom} in a space--time form. For example, if the backward Euler time integrator\footnote{ Here we only consider the backward Euler time integrator for the simplicity of illustration. However, the extension to other linear multistep implicit time integrators is straight forward.  } is used for time domain discretization, then the corresponding space--time formulation to Eq.~\eqref{eq:backwardEuler} becomes
    \begin{equation}\label{eq:st_fom_BE}
      \Abe \solst = \dt\fluxst + \extrast,   
    \end{equation}
    where
    \begin{equation} \label{eq:st_fom_term_definition}
      \Abe = \bmat{\mass &       &        &
\\         -\mass & \mass &        &
\\                &\ddots & \ddots &
\\                &       & -\mass & \mass}, 
\quad
      \solst = \pmat{\solArg{1} \\ \solArg{2} \\ \vdots \\ \solArg{\ntimedof}},
\quad
      \fluxst = \pmat{\fluxArg{1} \\ \fluxArg{2} 
      \\ \vdots \\ \fluxArg{\ntimedof}}, 
\quad
      \extrast = \pmat{\mass\solinit \\ \zerobold 
      \\ \vdots \\ \zerobold}. 
    \end{equation} 
    Note that $\solst(\param)$ denotes the parameterized space--time state implicitly defined as the solution to the problem~\eqref{eq:st_fom_BE} with $\solst:\paramDomain\rightarrow \RR{\nspacedof\ntimedof}$ and $\solst(\param)\in\RR{\nspacedof\ntimedof}$. Further, $\fluxst(\solst;\param)\in\RR{\nspacedof\ntimedof}$ denotes the space--time nonlinear term that is nonlinear in at least its first argument with $\fluxst:\RR{\nspacedof\ntimedof}\times\paramDomain\rightarrow \RR{\nspacedof\ntimedof}$.  The space--time residual function $\resst\in\RR{\nspacedof\ntimedof}$  corresponding to Eq.~\eqref{eq:st_fom_BE} is defined as
    \begin{align}\label{eq:res_st_fom_BE}
      \begin{split}
        \resst_{\text{BE}}(\solst;\param) &:= \Abe \solst(\param) 
        - \dt\fluxst(\solst;\param) - \extrast(\param).
\\                           &=\zerobold 
      \end{split}
    \end{align}

    To reduce both the spatial and temporal dimensions of the full-order model, the ST-LSPG method enforces the approximated numerical solution $\solapproxST\in\RR{\nspacedof\ntimedof}$ to reside in an affine `space--time trial subspace'
\begin{equation}\label{eq:spacetimeLSPGcontinuous}
\solapproxST\in\spacetimesubspace
\subseteq\RR{\nspacedof\ntimedof},
\end{equation}
where  
$\spacetimesubspace\defeq\ones_{\ntimedof}\otimes\solref
+\Span{\stbasisvec{i}}_{i=1}^\nbasisst\subseteq\RR{\nspacedof\ntimedof}$
with
$\dim(\spacetimesubspace)=\nbasisst\ll\nspacedof\ntimedof$ and
$\ones_{\ntimedof} \in \RR{\ntimedof}$ whose elements are all one.
Here, $\otimes$ denotes the Kronecker product of two matrices and the product of
two matrices $\Abold\in\RR{I\times J}$ and $\Bbold\in\RR{K\times L}$ is
defined as
  \begin{equation}\label{eq:kron_prod}
    \Abold\otimes\Bbold = 
    \bmat{a_{11}\Bbold & \cdots  &  a_{1J}\Bbold 
  \\ \vdots & \ddots & \vdots
  \\ a_{I1}\Bbold & \cdots & a_{IJ}\Bbold}.
  \end{equation}
Further, $\stbasisvec{i}\in\RR{\nspacedof\ntimedof}$ denotes a space--time basis vector.  Thus, the ST-LSPG method approximates the numerical solution as
\begin{align}\label{eq:spacetimeExpansion}
\begin{split}
  \solst(\param) \approx\solapproxST(\param)&= \ones_{\ntimedof}\otimes\solref +
\sum_{i=1}^{\nbasisst}
\stbasisvec{i}\redsolapproxSTArg{i} 
\end{split}
\end{align}
where 
$\redsolapproxSTArg{i}\in\RR{}$, $i\innat{\nbasisst}$
denotes the generalized coordinate of the ST-LSPG solution. 

A space--time residual vector function can now be defined with the generalized
coordinates as an argument.  Replacing $\solst$ in Eq.~\eqref{eq:res_st_fom_BE}
with $\solapproxST$ gives the residual vector function 
\begin{equation}\label{eq:st_ROM_BE2}
  \resRedst_{\text{BE}}(\redsolapproxST;\param) 
  := \Abe \stbasismat\redsolapproxST 
  - \dt\fluxst(\solinitst+\stbasismat\redsolapproxST;\param) 
   - \extrast(\param) + \Abe\solinitst(\param),   
\end{equation}
where $\solinitst := \ones_{\ntimedof}\otimes\solinit$ and $\stbasismat\in\RR{\nspacedof\ntimedof\times\nbasisst}$ denotes a space--time basis matrix that is defined as $\stbasismat := \bmat{\stbasisvec{1} & \cdots & \stbasisvec{\nbasisst}}$ and $\redsolapproxST\in\RR{\nbasisst}$ denotes the generalized coordinate vector that is defined as 
\begin{equation}\label{eq:def_stgencoord}
\redsolapproxST := \bmat{\redsolapproxSTArg{1}&\cdots&\redsolapproxSTArg{\nbasisst} }^T.
\end{equation}
Note that $\extrast(\param) + \Abe\solinitst(\param)$ vanishes.
Eq.~\eqref{eq:st_ROM_BE2} becomes
\begin{equation}\label{eq:st_ROM_BE}
  \resRedst_{\text{BE}}(\redsolapproxST;\param) 
  := \Abe \stbasismat\redsolapproxST 
  - \dt\fluxst(\solinitst+\stbasismat\redsolapproxST;\param).
\end{equation}
Reduced space--time residual vector functions for other time integrators can be defined similarly. We denote $\resRedst$ as a reduced spate--time residual vector function with a generic time integrator.

The space--time basis matrix $\stbasismat$ can be obtained by tensor product of spatial and temporal basis vectors.  The spatial basis vectors can be obtained via POD as in the DEIM and GNAT approaches.  The temporal basis vectors can be obtained via the following three tensor decompositions described in \cite{choi2019space}:
\begin{itemize}
  \item Fixed temporal subspace via T-HOSVD
  \item Fixed temporal subspace via ST-HOSVD
  \item Tailored temporal subspace via ST-HOSVD
\end{itemize}
The ST-HOSVD method is a more efficient version of T-HOSVD.  Thus, we will not
consider T-HOSVD. The tailored temporal subspace via ST-HOSVD has a LL1 form
that has appeared, for example, in \cite{sorber2013optimization,de2011blind}.
Therefore, we will refer to it as the LL1 decomposition.

Because of the reduction in spatial and temporal dimension, the space--time
residual vector $\resRedst$ cannot achieve zero most likely.  Thus, the ST-LSPG
method minimizes the square norm of $\resRedst$ and computes the ST-LSPG solution:
\begin{align} \label{eq:t-lspgReducedLargerSimplify}
\begin{split}
 \redsolapproxST(\param) &=
	\underset{\solRedDummyOpt\in\RR{\nbasisst}}{\arg\min}
    \left \|\resRedst(\solRedDummyOpt;\param)
		\right \|_2^2.
\end{split}
\end{align}
The ST-LSPG ROM solves Eq.~\eqref{eq:t-lspgReducedLargerSimplify} without any hyper-reduction. As in the DEIM and GNAT approaches, a hyper-reduction is required for a considerable speed-up due to the presence of the space--time nonlinear residual vector.  Therefore, the \STGNAT\ method approximates the space--time nonlinear residual terms with gappy POD \cite{everson1995karhunen}, which in turn requires construction of a space--time residual basis.  Similar to the GNAT method, the ST-GNAT method approximates the space--time nonlinear residual term as
 \begin{equation}\label{eq:STGNAT_approx}
   \resRedst \approx \basismatrest\redres,
 \end{equation}
 where $\basismatrest \defeq [\basismatrestvec{1},\ldots,\basismatrestvec{\nbasisres} ] \in \RR{\ndof\ntimedof\times\nbasisres}$, $\nbasisst \leq \nbasisres \ll \ndof\ntimedof$, denotes the space--time residual basis matrix and $\redres \in \RR{\nbasisres}$ denotes the generalized coordinates of the nonlinear residual term. The \STGNAT\ solves the following space--time least-squares problem to obtain the generalized coordinates, $\redres$:
\begin{align} \label{eq:STGNAT_least-squares}
  \begin{split} 
    \redres = \argmin{\reddummy\in\RR{\nbasisres}} \quad&
    \frac{1}{2} \left \|\samplematst(\resRedst - \basismatrest\reddummy)
    \right \|_2^2.
  \end{split} 
\end{align}
 where
 $\samplematst\defeq[\unitvecArg{p_1},\ldots,\unitvecArg{p_{\nressample}}]^T
 \in\RR{\nressample\times\ndof\ntimedof}$, $\nbasisst \leq \nbasisres \leq
 \nressample \ll \ndof\ntimedof$, is the sampling matrix and $\unitvecArg{p_i}$
 is the $p_i$th column of the identity matrix
 $\identity{\ndof\ntimedof}\in\RR{\ndof\ntimedof\times\ndof\ntimedof}$.  The
 solution to Eq.~\eqref{eq:STGNAT_least-squares} is given by
 \begin{equation}\label{eq:STGNAT-generalizedcoordinates}
   \redres = (\samplematst\basismatrest)^\dagger\samplematst\resRedst.
 \end{equation}
 Therefore, Eq.~\eqref{eq:STGNAT_approx} becomes
 \begin{equation}\label{eq:STGNAT_r}
   \resRedst \approx\obliqueprojector_{\STGNAT} \resRedst,
 \end{equation}
 where $\obliqueprojector_{\STGNAT} := \basismatrest
 (\samplematst\basismatrest)^\dagger\samplematst$ is the $\STGNAT$ oblique projection
 matrix. Note that $\obliqueprojector_{\STGNAT}$ has the same structure as
 $\obliqueprojector_{\GNAT}$.  
 The ST-GNAT method does not construct the sampling matrix $\samplematstNT$.
 Instead, it maintains the sampling indices $\{p_1,\ldots,p_{\nbasisflux}\}$ and
 corresponding rows of $\basismatrest$ and $\resRedst$. This enables the ST-GNAT
 to achieve a speed-up when it is applied to nonlinear problems.

 Section 5.3 in \cite{choi2019space} discusses three different options to determine 
 the sampling indices (i.e., $\samplematstNT$). However, all these three options
 are simple variations of Algorithm 3 in \cite{carlberg2013gnat} and Algorithm 5
 in \cite{carlberg2011efficient}. They all
 minimize the error in the gappy reconstruction of the POD basis vectors for
 nonlinear space and time residuals. Therefore, the space--time nonlinear
 residual projection error due to \eqref{eq:STGNAT_r} is similar to the ones in
 the GNAT method. That is, 
 \begin{equation}\label{eq:STGNATprojection_error}
   \|\resRedst - \basismatrest (\samplematst\basismatrest)^\dagger\samplematst\resRedst\|_2
   \leq \| \righttrianglemat^{-1} \|_2 
   \|\resRedst - \basismatrest\basismatrest^T\resRedst \|_2,
 \end{equation}
 where $\righttrianglemat$ is a triangle matrix from QR factorization of
 $\samplematst\basismatrest$ (i.e.,
 $\samplematst\basismatrest=\orthogonalmat\righttrianglemat$). On the other
 hand, one can also apply the sRRQR factorization in Algorithm 4 with a tuning
 parameter $\tuningparam$ of
 \cite{gu1996efficient} to $\basismatrest^T$ to obtain $\samplematstNT$ that is
 associated with a tighter error bound for
 the projection error:
 \begin{equation}\label{eq:betterSTGNATprojection_error}
   \|\resRedst - \basismatrest (\samplematst\basismatrest)^\dagger\samplematst\resRedst\|_2
   \leq \sqrt{1+\tuningparam^2\nbasisres(\ndof-\nbasisres)}
   \|\resRedst - \basismatrest\basismatrest^T\resRedst \|_2.
 \end{equation}
 This error bound can be obtained by setting the identity matrix as a weight
 matrix in Theorem 4.8 of \cite{drmac2018discrete}.

 Finally, the ST-GNAT solves the following least-squares problem at every time
 step, for example, with the backward Euler time integrator:
  \begin{align} \label{eq:STGNAT-spOpt1}
    \begin{split} 
      \redsolapproxST(\param) = \argmin{\reddummy\in\RR{\nbasisspace}} \quad&
      \frac{1}{2} \left \|\ (\samplematst\basismatrest)^\dagger\samplematst
      \resRedst_{\BE}(\reddummy;\param)
      \right \|_2^2.
    \end{split} 
  \end{align}

The original ST-GNAT paper introduces three different ways of collecting space--time residual snapshots, that are in turn used for the space--time residual basis construction (see Section 5.2 in \cite{choi2019space}).  Below is a list of the approaches introduced in \cite{choi2019space} and explains advantages and disadvantages of each:
\begin{enumerate}
  \item {\it ST-LSPG ROM training iterations.} This approach takes the
    space--time residual snapshot from every Newton iteration of the ST-LSPG
    simulations at training points in $\paramDomainTrain$.  This case leads to
    the number of residual snapshots, $\nrestrain =
    \sum_{\param\in\paramDomainTrain}(\kmax(\param)+1)$, where $\kmax(\param)$
    is the number of Newton iterations taken to solve the ST-LSPG simulation for
    the training point $\param$.  This approach is not realistic for a
    large-scale problem because it requires $\card{\paramDomainTrain}$ training
    simulations of the computationally expensive ST-LSPG ROM, where
    $\card{\paramDomainTrain}$ denotes the cardinality of the set
    $\paramDomainTrain$. Furthermore, it
    requires the extra SVD on the residual snapshots, which is not necessary for
  our proposed method.  \item {\it Projection of FOM training solutions.} This
    approach takes the following steps:
    \begin{enumerate}
      \item take FOM state solution at every Newton iteration.
      \item re-arrange them in the space--time form (i.e., $\solst$ in Eq.~\eqref{eq:st_fom_term_definition}).\footnote{
        Note that these are FOM solutions from time marching algorithms in which
        each time step results in different number of Newton iterations if
        implicit time integrators are used.  Some time steps take a smaller
        number of Newton iterations than other time steps. However, in order to
        re-arrange each Newton iterate state in the space--time form, we must
        have the same number of Newton iterations at each time step.  Therefore,
        for the time steps that have converged with a smaller number of Newton
        iterations than other time steps, we pad the solution state vectors of
        the Newton iterations beyond the convergence with the converged
        solution.  This only applies to an implicit time integrator because an explicit time integrator does not require any Newton solve.
        }
      \item project them onto the space--time subspace, $\spatiotemporalSubspace$ (i.e., $\solapproxST = \stbasismat(\stbasismat^T\stbasismat)^{-1}\stbasismat^T(\solst-\solinitst)$).
      \item compute the corresponding space--time residual (e.g., $\resst(\solapproxST;\param)$ in Eq.~\eqref{eq:res_st_fom_BE} in the case of the backward Euler time integrator).
      \item use those residuals as residual snapshots.
    \end{enumerate}
    This approach simply requires $\nrestrain$ projections and evaluations of
    the space--time residual.  However, it requires the extra SVD on the residual snapshots, which is not necessary for our proposed method.
  \item {\it Random samples.} This approach generates a random space--time
    solution state samples (e.g., via Latin hypercube sampling or random
    sampling from uniform distribution) and computes the corresponding
    space--time residual (e.g., $\resst(\solapproxST;\param)$ in
    Eq.~\eqref{eq:res_st_fom_BE} in the case of the backward Euler time
    integrator).  This approach simply requires $\nrestrain$ random sample
    generations and evaluations of the space--time residual.  However, random
    samples are hardly correlated with actual data. Therefore, it is likely to
    generate poor space--time residual subspace.  Furthermore, it requires the extra SVD on the residual snapshots, which is not necessary for our proposed method.
\end{enumerate}

\section{Solution-based Nonlinear Subspace (\methodAcronym) method}\label{sec:SNS}
  Finally, we state our proposed method that avoids collecting nonlinear term
  snapshots and additional POD for the DEIM and GNAT approaches or additional
  tensor decomposition for the ST-GNAT method.  We propose to use solution
  snapshots to construct nonlinear term basis in the DEIM, GNAT, and ST-GNAT
  approaches. A justification for using the solution snapshots comes from the
  subspace inclusion relation between the subspace spanned by the solution
  snapshots and the subspace spanned by the nonlinear term snapshots as shown in
  Eqs.~\eqref{eq:forwardEuler_total_inclusion} and
  \eqref{eq:backwardEuler_total_inclusion} with the forward and backward Euler
  time integrators.\footnote{Subspace inclusion relations for other time
  integrators are shown in Appendix~\ref{sec:appendix}.}
  \subsection{DEIM-SNS}\label{sec:DEIM-SNS}
    We are going back to Eq.~\eqref{eq:rom} and replace the nonlinear term with the approximation in Eq.~\eqref{eq:DEIM_f}:
    \begin{equation}\label{eq:DEIM}
      \mass\basismatspace\frac{d\redsolapprox}{dt} =
      \basisfluxmat (\samplemat\basisfluxmat)^{-1}\samplemat
      \flux(\solArg{0}+\basismatspace\redsolapprox,t; \param).
    \end{equation}
    Eq.~\eqref{eq:DEIM} is an over-determined system, so it is likely that
    it will not have a solution. 
    However, if there is a solution, then necessary conditions
    for Eq.~\eqref{eq:DEIM} to have a non-trivial solution (i.e., $\redsolapprox
    \neq \zerobold$) are $\flux(\solArg{0},t;\param) \neq \zerobold$ and 
    \begin{equation}\label{eq:DEIMnecessary}
      \Range{\mass\basismatspace} \cap \Range{\basisfluxmat} \neq \{\zero\}. 
    \end{equation}
    The second condition says that the intersection of
    $\Range{\mass\basismatspace}$ and $\Range{\basisfluxmat}$ needs to be
    non-trivial if there is a non-trivial solution 
    to Eq.~\eqref{eq:DEIM}. Typically, we build $\basismatspace$ first, using the POD
    approach explained in Section~\ref{sec:DEIM} and
    $\Range{\mass\basismatspace}$ is set by $\basismatspace$. 
    Therefore, the intersection of $\Range{\mass\basismatspace}$ and 
    $\Range{\basisfluxmat}$ can be controlled by the choice of $\basisfluxmat$
    we made. The larger the intersection of those two range spaces are, 
    the more likely it is that there is a solution to Eq.~\eqref{eq:DEIM}. 
    Given $\basismatspace$, the largest subspace intersection it
    can be is $\Range{\mass\basismatspace}$, i.e.,
    \begin{equation}\label{eq:conformingCondition}
      \Range{\mass\basismatspace} \cap \Range{\basisfluxmat} =
      \Range{\mass\basismatspace}.
    \end{equation}
    We call 
    this condition as {\it \bf the conforming subspace condition}. The conforming
    subspace condition leads to two obvious choices for $\basisfluxmat$:
    \begin{itemize}
      \item The first choice is to ensure 
        $\Range{\mass\basismatspace} = \Range{\basisfluxmat}$.
        If $\basisfluxmat = \mass\basismatspace$, then the range space of
      the left and right-hand sides of Eq.~\eqref{eq:DEIM} are the same. This leads 
      Eq.~\eqref{eq:DEIM} to become
      \begin{equation}\label{eq:DEIM-SNS}
        \mass\basismatspace\frac{d\redsolapprox}{dt} =
        \mass\basismatspace (\samplemat\mass\basismatspace)^{-1}\samplemat
        \flux(\solArg{0}+\basismatspace\redsolapprox,t; \param).
      \end{equation}
      Because Eq.~\eqref{eq:DEIM-SNS} is an over-determined system and 
        unlikely to have a solution, 
      applying the Galerkin projection to Eq.~\eqref{eq:DEIM-SNS} becomes:
      \begin{equation}\label{eq:DEIM-SNS-Galerkin}
        \basismatspace^T\mass\basismatspace\frac{d\redsolapprox}{dt} =
        \basismatspace^T\mass\basismatspace (\samplemat\mass\basismatspace)^{-1}\samplemat
        \flux(\solArg{0}+\basismatspace\redsolapprox,t; \param).
      \end{equation}
      Assuming $\basismatspace^T\mass\basismatspace$ is invertible,
      Eq.~\eqref{eq:DEIM-SNS-Galerkin} becomes:
      \begin{equation}\label{eq:DEIM-SNS-Galerkin-reduced}
        \frac{d\redsolapprox}{dt} =
        (\samplemat\mass\basismatspace)^{-1}\samplemat
        \flux(\solArg{0}+\basismatspace\redsolapprox,t; \param).
      \end{equation}
      For the special case of $\mass$ being an identity matrix,
      Eq.~\eqref{eq:DEIM-SNS-Galerkin-reduced} becomes:
      \begin{equation}\label{eq:DEIM-SNS-Galerkin-reduced-Imass}
        \frac{d\redsolapprox}{dt} =
        (\samplemat\basismatspace)^{-1}\samplemat
        \flux(\solArg{0}+\basismatspace\redsolapprox,t; \param).
      \end{equation}

      \item  The second choice is to ensure
        $\Range{\mass\basismatspace} \subset \Range{\basisfluxmat}$.
        This can be achieved by taking 
      an extended solution basis, 
      $\basismatspaceext \in \RR{\nspacedof\times\nbasisspaceext}$ with
      $\nbasisspace < \nbasisspaceext \ll \nspacedof$ and 
      $\Range{\basismatspace} \subset \Range{\basismatspaceext}$.
      Then we set $\basisfluxmat = \mass\basismatspaceext$, which leads to
      $\Range{\mass\basismatspace} \subset \Range{\mass\basismatspaceext}$.
      The obvious choice for $\basismatspaceext$ is to take a larger truncation
      of the left singular matrix from SVD of the solution snapshot matrix
      than $\basismatspace$. Note that this particular choice of $\basismatspaceext$ 
      results in the first $\nbasisspace$ columns of $\basismatspaceext$ being the same 
      as $\basismatspace$ 
      (i.e., $\basismatspaceext = \bmat{\basismatspace & \basismatspaceextra}$
      where $\basismatspaceextra\in\RR{\nspacedof\times\nbasisspaceextra}$ with
      $\nbasisspaceextra = \nbasisspaceext -\nbasisspace$ and
      $\basismatspace^T\basismatspaceextra = \zerobold \in\RR{\nbasisspace\times\nbasisspaceextra}$).
      By setting $\basisfluxmat = \mass\basismatspaceext$,
      Eq.~\eqref{eq:DEIM} becomes 
      \begin{equation}\label{eq:DEIM-SNS2}
        \mass\basismatspace\frac{d\redsolapprox}{dt} =
        \mass\basismatspaceext (\samplemat\mass\basismatspaceext)^{-1}\samplemat
        \flux(\solArg{0}+\basismatspace\redsolapprox,t; \param).
      \end{equation}
      Because Eq.~\eqref{eq:DEIM-SNS2} is unlikely to have a solution, 
      applying the Galerkin projection to Eq.~\eqref{eq:DEIM-SNS2} becomes:
      \begin{align} \label{eq:DEIM-SNS-Galerkin2} 
        \basismatspace^T\mass\basismatspace\frac{d\redsolapprox}{dt} &=
        \basismatspace^T\mass\basismatspaceext (\samplemat\mass\basismatspaceext)^{-1}\samplemat
        \flux(\solArg{0}+\basismatspace\redsolapprox,t; \param) \\
        &= \basismatspace^T\mass\bmat{\basismatspace & \basismatspaceextra} 
        (\samplemat\mass\basismatspaceext)^{-1}\samplemat
        \flux(\solArg{0}+\basismatspace\redsolapprox,t; \param).
      \end{align}
      Assuming $\basismatspace^T\mass\basismatspace$ is invertible,
      Eq.~\eqref{eq:DEIM-SNS-Galerkin2} becomes:
      \begin{equation}\label{eq:DEIM-SNS-Galerkin-reduced2}
        \frac{d\redsolapprox}{dt} =
        \bmat{\identity{\nbasisspace} & (\basismatspace^T\mass\basismatspace)^{-1}
        (\basismatspace^T\mass\basismatspaceextra) }
        (\samplemat\mass\basismatspaceext)^{-1}\samplemat
        \flux(\solArg{0}+\basismatspace\redsolapprox,t; \param).
      \end{equation}
      For the special case of $\mass$ being an identity matrix, 
      Eq.~\eqref{eq:DEIM-SNS-Galerkin-reduced2} becomes
      \begin{equation}\label{eq:DEIM-SNS-Galerkin-reduced-special}
        \frac{d\redsolapprox}{dt} =
        \bmat{\identity{\nbasisspace} & \zerobold} 
        (\samplemat\basismatspaceext)^{-1}\samplemat
        \flux(\solArg{0}+\basismatspace\redsolapprox,t; \param).
      \end{equation}
    \end{itemize}
    The DEIM-SNS approach solves either Eq.~\eqref{eq:DEIM-SNS-Galerkin-reduced} or
    Eq.~\eqref{eq:DEIM-SNS-Galerkin-reduced2} depending on the choice of
    $\basisfluxmat$ above.
    Applying a time integrator to 
    Eq.~\eqref{eq:DEIM-SNS-Galerkin-reduced} or 
    Eq.~\eqref{eq:DEIM-SNS-Galerkin-reduced-special}
    leads to a reduced
    O$\Delta$E, whose solution, $\redsolapprox_\star$, 
    can be lifted to find the full order approximate solution via
    $\solapproxfom = \solArg{0}+\basismatspace\redsolapprox_\star$.

    Additionally, the subspace inclusion
    relations\footnote{
      Eq.~\eqref{eq:forwardEuler_total_inclusion} 
    of the forward Euler time integrator, 
    Eq.~\eqref{eq:AF2_total_inclusion}
    of the Adams--Moulton time integrator, and  
    Eq.~\eqref{eq:RK2_inclusion} of the midpoint Runge-Kutta method show the
    subspace inclusion relations
    for explicit time integrators.
    }
    show that the subspace spanned by the solution snapshots
    include the subspace spanned by the nonlinear term snapshots.
    This fact further motivates the use of solution snapshots to build a
    nonlinear term basis.
    Indeed, numerical experiments show that
    the solution accuracy obtained by the DEIM-SNS approach is comparable to
    the one obtained by the traditional DEIM approach. 
    For example, see Figs.~\ref{fig:forwardEuler_DEIM_nldiff},
    and \ref{fig:backwardEuler_DEIM_nldiff}
    in Section~\ref{sec:DEIMvsDEIMSNS}.

    The obvious advantage of the DEIM-SNS approach over DEIM is that
    no additional SVD or eigenvalue decomposition is required, which
    can save the computational cost of the offline phase.

    {\remark\label{remark:DEIM-SNS} 
    Although the numerical experiments show that the two choices of
    $\basisfluxmat$ above give promising results, the error analysis in 
    Section~\ref{sec:erroranalysis} shows that the nonlinear term projection error bound
    increases by the condition number of $\mass$ (see Theorem
    \ref{theorem:DEIM-SNSbound}). This is mainly because the orthogonality of
    $\basisfluxmat$ is lost when $\basisfluxmat =
    \mass\basismatspace$ or $\mass\basismatspaceext$. This issue
    is resolved by orthogonalizing $\basisfluxmat$ (e.g., apply economic 
    QR factorization $\basisfluxmat = \orthogonalmat\righttrianglemat$)
    before using it as nonlinear
    term basis. Then apply, for example, Algorithm 4 of \cite{gu1996efficient}
    to the transpose of orthogonalized one (e.g., $\orthogonalmat^T$) to
    generate a sampling matrix $\samplematNT$. This procedure eliminates the
    condition number of $\mass$ in the error bound because 
    \eqref{eq:tighterBound} is valid.}

    {\remark\label{remark:weightedInnerProduct}
    Inspired by the weighted inner product space introduced for DEIM in
    \cite{drmac2018discrete}, another oblique DEIM-SNS projection is possible.
    With the weight matrix $\Wbold =
    \mass^{-T}\mass^{-1}$, the selection operator $\Sbold^T =
    \samplemat\mass^T$, and the basis $\hat{\Ubold} = \mass\basismatspace$ according to
    Section 4.2 of \cite{drmac2018discrete}, the weighted oblique DEIM-SNS projection
    can be defined as:
    \begin{align}\label{eq:weightedObliqueProjection}
      \obliqueprojector_{\SNS} &=
      \hat{\Ubold}(\Sbold^T\Wbold\hat{\Ubold})^{\dagger}\Sbold^T\Wbold \\
      &=
      \mass\basismatspace(\samplemat\basismatspace)^{\dagger}\samplemat\mass^{-1}.
    \end{align}
    }



  \subsection{GNAT-SNS}\label{sec:GNAT-SNS}
  The GNAT method needs to build a nonlinear residual term basis, $\basismatres$,
  as described in Section~\ref{sec:GNAT}. The nonlinear residual term is nothing
  more than linear combinations of the nonlinear term and the time derivative
  approximation as in Eq.~\eqref{eq:trialsub_residual_BE}.
  Thus, the subspace spanned by the nonlinear term residual snapshots
  are included in the subspace spanned by the solution snapshots.
  This motivates the use of the solution snapshots 
  for the construction of a nonlinear residual term basis. 
  Therefore, the same type of the nonlinear term basis in the DEIM-SNS approach
  can be used to construct the nonlinear residual term basis in the GNAT-SNS method:
  \begin{itemize}
    \item The first choice is to set $\basismatres = \mass\basismatspace$.
      Thus, for example, the GNAT-SNS method solves the following
      least-squares problem with the backward Euler time integrator:
      \begin{align} \label{eq:GNAT-SNS-least-squares}
        \begin{split} 
          \redsolapproxArg{n} = \argmin{\reddummy\in\RR{\nbasisspace}} \quad&
          \frac{1}{2} \left \|\mass\basismatspace(\samplemat\mass\basismatspace)^\dagger
          \samplemat
          (\mass\basismatspace(\reddummy - \redsolapproxArg{n-1}) 
          -\dt\flux(\solArg{0}+\basismatspace\redsolapproxArg{n},t;\param))
          \right \|_2^2,
        \end{split} 
      \end{align}
      which can be manipulated to
      \begin{align} \label{eq:GNAT-SNS-least-squares-simplified}
        \begin{split} 
          \redsolapproxArg{n} = \argmin{\reddummy\in\RR{\nbasisspace}} \quad&
          \frac{1}{2} \left \|\mass\basismatspace(\reddummy - \redsolapproxArg{n-1}) 
          -\dt\mass\basismatspace(\samplemat\mass\basismatspace)^\dagger
                    \samplemat
          \flux(\solArg{0}+\basismatspace\redsolapproxArg{n},t;\param)
          \right \|_2^2.
        \end{split} 
      \end{align}
      Note that the terms in $\ell_2$ norm in
      Eq.~\eqref{eq:GNAT-SNS-least-squares-simplified} is very similar to
      the discretized DEIM residual
      before Galerkin projection (i.e., applying the backward Euler 
      time integrator to Eq.~\eqref{eq:DEIM-SNS} gives the terms in $\ell_2$
      norm in Eq.~\eqref{eq:GNAT-SNS-least-squares-simplified}).
      They are only different by the fact that one is an inverse and the other one
      is the Moore--Penrose inverse. In fact, if $\nressample = \nbasisres$, then
      Eq.~\eqref{eq:GNAT-SNS-least-squares-simplified} is equivalent to 
      applying the backward Euler time integrator to Eq.~\eqref{eq:DEIM-SNS} and
      minimize the $\ell_2$ norm of the corresponding residual.

      For the special case of $\mass$ being an identity matrix, 
      Eq.~\eqref{eq:GNAT-SNS-least-squares-simplified} becomes
      \begin{align} \label{eq:GNAT-SNS-least-squares-simplified2}
        \begin{split} 
          \redsolapproxArg{n} = \argmin{\reddummy\in\RR{\nbasisspace}} \quad&
          \frac{1}{2} \left \|\reddummy - \redsolapproxArg{n-1} 
          -\dt(\samplemat\basismatspace)^\dagger
                    \samplemat
          \flux(\solArg{0}+\basismatspace\redsolapproxArg{n},t;\param)
          \right \|_2^2.
        \end{split} 
      \end{align}
      
    \item The second choice is to set $\basismatres = \mass\basismatspaceext$
      where $\basismatspaceext = \bmat{\basismatspace & \basismatspaceextra}$ 
      as in Section~\ref{sec:DEIM-SNS}.
      This leads to the following least-squares problem, for example, 
      using the backward Euler time integrator:
      \begin{align} \label{eq:GNAT-SNS-least-squares2}
        \begin{split} 
          \redsolapproxArg{n} = \argmin{\reddummy\in\RR{\nbasisspace}} \quad&
          \frac{1}{2} \left \|\mass\basismatspaceext(\samplemat\mass\basismatspaceext)^\dagger
          \samplemat
          (\mass\basismatspace(\reddummy - \redsolapproxArg{n-1}) 
          -\dt\flux(\solArg{0}+\basismatspace\redsolapproxArg{n},t;\param))
          \right \|_2^2.
        \end{split} 
      \end{align}
      For the special case of $\mass$ being an identity matrix, 
      Eq.~\eqref{eq:GNAT-SNS-least-squares2} becomes
      \begin{align} \label{eq:GNAT-SNS-least-squares-reduced2}
        \begin{split} 
          \redsolapproxArg{n} = \argmin{\reddummy\in\RR{\nbasisspace}} \quad&
          \frac{1}{2} \left \|(\samplemat\basismatspaceext)^\dagger
          \samplemat
          (\basismatspace(\reddummy - \redsolapproxArg{n-1}) 
          -\dt\flux(\solArg{0}+\basismatspace\redsolapproxArg{n},t;\param))
          \right \|_2^2.
        \end{split} 
      \end{align}
  \end{itemize}
  The GNAT-SNS method solves either Eq.~\eqref{eq:GNAT-SNS-least-squares-simplified} or
  Eq.~\eqref{eq:GNAT-SNS-least-squares2} depending on the choice of
  $\basismatres$ above.
  For the special case of $\mass$ being an identity matrix, 
  the GNAT-SNS method solves either Eq.~\eqref{eq:GNAT-SNS-least-squares-simplified2} or
  Eq.~\eqref{eq:GNAT-SNS-least-squares-reduced2} depending on the choice of
  $\basismatres$.
  The reduced solution $\redsolapproxArg{n}$ 
  can be lifted to find the full order approximate solution via
  $\solapproxArg{n} = \solArg{0}+\basismatspace\redsolapproxArg{n}$.

   { \remark\label{remark:GNAT-SNS} 
    Similar to Remark~\ref{remark:DEIM-SNS}, 
    the error analysis in Section~\ref{sec:erroranalysis} shows 
    that the nonlinear residual projection error bound, 
    regarding the two choices of $\basismatres$ above,
    increases by the condition number of $\mass$ (see Theorem
    \ref{theorem:GNAT-SNSbound}). This is mainly because the orthogonality of
    $\basismatres$ is lost when $\basismatres =
    \mass\basismatspace$ or $\mass\basismatspaceext$. This issue
    is resolved by orthogonalizing $\basismatres$ (e.g., apply economic 
    QR factorization $\basismatres = \orthogonalmat\righttrianglemat$)
    before using it as nonlinear
    residual basis. Then apply, for example, Algorithm 4 of \cite{gu1996efficient}
    to the transpose of orthogonalized one (e.g., $\orthogonalmat^T$) to
    generate a sampling matrix $\samplematNT$. This procedure eliminates the
    condition number of $\mass$ in the error bound because 
    \eqref{eq:betterGNATerror} is valid.}

  \subsection{\STGNAT-SNS}\label{sec:ST-GNAT-SNS}
  The ST-GNAT method needs to build a space--time nonlinear residual term basis, 
  $\basismatrest$, as described in Section~\ref{sec:STGNAT}. 
  We are going back to Eq.~\eqref{eq:st_ROM_BE} to find 
  {\it \bf the conforming subspace condition} for ST-GNAT-SNS. 
  In order to increase the chance of
  making the space--time residual function
  defined in Eq.~\eqref{eq:st_ROM_BE} zero, 
  the following conforming subspace condition can be made:
  \begin{equation}\label{eq:STGNAT-SNS_conforming}
    \Range{\Abe\stbasismat} \cap \Range{\basismatrest} = \Range{\Abe\stbasismat}.
  \end{equation}
  Therefore, we propose the following
  bases of the space--time nonlinear residual term:
  \begin{itemize}
    \item The first choice is to set $\basismatrest = \Abe\stbasismat$. Thus,
      for example, the ST-GNAT-SNS method solves the following least-squares problem
      with the backward Euler time integrator:
      \begin{align} \label{eq:STGNAT-SNS-spOpt1}
        \begin{split} 
          \redsolapproxST(\param) = \argmin{\reddummy\in\RR{\nbasisspace}} \quad&
          \frac{1}{2} \left \|\ (\samplematst\Abe\stbasismat)^\dagger\samplematst
          (\Abe \stbasismat\redsolapproxST 
          - \dt\fluxst(\Abe\solinitst+\stbasismat\redsolapproxST;\param) 
          )
          \right \|_2^2,
        \end{split} 
      \end{align}
      which can be rewritten as
      \begin{align} \label{eq:STGNAT-SNS-spOpt2}
        \begin{split} 
          \redsolapproxST(\param) = \argmin{\reddummy\in\RR{\nbasisspace}} \quad&
          \frac{1}{2} \left \|\ 
          \redsolapproxST 
          - \dt(\samplematst\Abe\stbasismat)^\dagger\samplematst
          \fluxst(\Abe\solinitst+\stbasismat\redsolapproxST;\param) 
          \right \|_2^2.
        \end{split} 
      \end{align}
      This is what the ST-GNAT-SNS method solves if $\basismatrest =
      \Abe\stbasismat$.
    \item The second choice is to set $\basismatrest = \Abe\stbasismatext$ where
      $\stbasismatext\in\RR{\ndof\ntimedof\times\nbasisspaceext}$ with
      $\nspacedof\ntimedof \gg \nressample \geq \nbasisspaceext > \nbasisst$ and 
      $\Range{\stbasismat} \subset \Range{\stbasismatext}$.
      The obvious choice for $\stbasismatext$ is to take a larger truncation
      of the factor matrices from the tensor decomposition (e.g., 
      ST-HOSVD and LL1) of the solution snapshot tensor 
      than $\stbasismat$. Note that this particular choice of $\stbasismatext$ 
      results in the first $\nbasisst$ columns of $\stbasismatext$ being the same 
      as $\stbasismat$ 
      (i.e., $\stbasismatext = \bmat{\stbasismat & \stbasismatextra}$
      where $\stbasismatextra\in\RR{\nspacedof\ntimedof\times\nbasisspaceextra}$ with
      $\nbasisspaceextra = \nbasisspaceext -\nbasisst$). 
      In this case, the ST-GNAT-SNS method solves the following least-squares problem:
      \begin{align} \label{eq:STGNAT-SNS-spOpt3}
        \begin{split} 
          \redsolapproxST(\param) = \argmin{\reddummy\in\RR{\nbasisspace}} \quad&
          \frac{1}{2} \left \|\ (\samplematst\Abe\stbasismatext)^\dagger\samplematst
          (\Abe \stbasismat\redsolapproxST 
          - \dt\fluxst(\Abe\solinitst+\stbasismat\redsolapproxST;\param) 
          )   
          \right \|_2^2,
        \end{split} 
      \end{align}
  \end{itemize}
  In addition to the choices above, we propose the following two choices 
  for the special case of $\mass$ being an identity matrix:
  \begin{itemize}
    \item The first choice is to set $\basismatrest = \stbasismat$. Thus,
      for example, the ST-GNAT-SNS method solves the following least-squares problem
      with the backward Euler time integrator:
      \begin{align} \label{eq:STGNAT-SNS-spOpt4}
        \begin{split} 
          \redsolapproxST(\param) = \argmin{\reddummy\in\RR{\nbasisspace}} \quad&
          \frac{1}{2} \left \|\ (\samplematst\stbasismat)^\dagger\samplematst
          (\Abe \stbasismat\redsolapproxST 
          - \dt\fluxst(\Abe\solinitst+\stbasismat\redsolapproxST;\param) 
          )
          \right \|_2^2,
        \end{split} 
      \end{align}
    \item The second choice is to set $\basismatrest =
      \stbasismatext$.
      In this case, the ST-GNAT-SNS method solves the following least-squares problem:
      \begin{align} \label{eq:STGNAT-SNS-spOpt5}
        \begin{split} 
          \redsolapproxST(\param) = \argmin{\reddummy\in\RR{\nbasisspace}} \quad&
          \frac{1}{2} \left \|\ (\samplematst\stbasismatext)^\dagger\samplematst
          (\Abe \stbasismat\redsolapproxST 
          - \dt\fluxst(\Abe\solinitst+\stbasismat\redsolapproxST;\param) 
          )   
          \right \|_2^2,
        \end{split} 
      \end{align}
  \end{itemize}

  The space--time generalized coordinates, $\redsolapproxST$, can be lifted
  to the approximate full space--time solution $\solapproxST(\param)$ 
  via Eq.~\eqref{eq:spacetimeExpansion}. 

   { \remark\label{remark:STGNAT-SNS} 
    Similar to Remarks~\ref{remark:DEIM-SNS} and \ref{remark:GNAT-SNS}, 
    the error analysis in Section~\ref{sec:erroranalysis} shows 
    that the nonlinear space--time residual projection error bound, 
    regarding the first two choices of $\basismatrest$ above,
    increases by the condition number of $\Abe$ (see Theorem
    \ref{theorem:STGNAT-SNSbound}). This is mainly because the orthogonality of
    $\basismatrest$ is lost when $\basismatrest =
    \Abe\stbasismat$ or $\Abe\stbasismatext$. 
    This issue is resolved by orthogonalizing $\basismatrest$ (e.g., apply economic 
    QR factorization $\basismatrest = \orthogonalmat\righttrianglemat$)\footnote{
    This might be challenging because of the size of $\basismatrest$. 
    However, a parallel QR factorization can be used to compute QR
    efficiently for a large matrix, such as
    \cite{buttari2008parallel,elmroth2000high}}
    before using it as nonlinear residual basis. 
    Then apply, for example, Algorithm 4 of \cite{gu1996efficient}
    to the transpose of orthogonalized one (e.g., $\orthogonalmat^T$) to
    generate a sampling matrix $\samplematstNT$. This procedure eliminates the
    condition number of $\Abe$ in the error bound because 
    \eqref{eq:betterSTGNATprojection_error} is valid.}

\section{Error analysis}\label{sec:erroranalysis}
    Section~\ref{sec:SNS} introduced the SNS method. If $\mass = \identity{}$,
    then all the same error analysis presented in Section~\ref{sec:ROMs} holds
    by replacing $\basisfluxmat$, $\basismatres$ with $\basismatspace$ or
    $\basismatspaceext$ for DEIM-SNS and GNAT-SNS and 
    $\basismatrest$ with $\stbasismat$ or $\stbasismatext$ for ST-GNAT-SNS.
    However, if $\mass \neq
    \identity{}$ is a general non-singular matrix, 
    then the error analysis has to be revisited. The error analysis below
    is inspired by \cite{chaturantabut2010nonlinear, drmac2016new,
    drmac2018discrete}.

    {\lemma\label{lemma:projectionbound} 
    Let $\dummy\in\RR{\nbig}$ is an arbitrary vector;
    $\nonsingularmat\in\RR{\nbig\times\nbig}$ is a non-singular matrix; 
    $\orthobasis\in\RR{\nbig\times\nsmall}$
    denotes a full rank matrix; $\samplematNT\in\RR{\nbig\times\nsample}$,
    $\nbig > \nsample \geq \nsmall$
    denotes a sampling matrix, defined in Section~\ref{sec:ROMs}. Let
    $\projection\in\RR{\nbig\times\nbig}$ 
    be an oblique projection matrix, defined as $\projection :=
    \nonsingularmat\orthobasis(\samplemat\nonsingularmat\orthobasis)^{\dagger}\samplemat$;
    let $\projections\in\RR{\nbig\times\nbig}$
    be another oblique projection matrix
    onto $\range{\nonsingularmat\orthobasis}$
    (i.e., $\projections := 
    \nonsingularmat\orthobasis(\nonsingularmat\orthobasis)^{\dagger}$).
    Then, a projection error bound is given by
    \begin{equation}\label{eq:projectionerror}
      \| (\identity{\nbig} - \projection)\dummy \|_2 \leq 
      \|(\samplemat\orthogonalmat)^{\dagger}\|_2
      \|(\identity{\nbig}-\projections)\dummy\|_2,
    \end{equation}
    where $\orthogonalmat$ is obtained from a QR factorization of
    $\nonsingularmat\orthobasis$, i.e., 
    $\nonsingularmat\orthobasis\defeq\orthogonalmat\righttrianglemat$.

    \proof Note that $\projection\projections = \projections$ is
    true becuase $\samplemat\nonsingularmat\orthobasis$ has full column rank. Thus,
    $(\identity{\nbig}-\projection)\projections = \zerobold$. This leads to
    \begin{equation*}
      (\identity{\nbig} - \projection)\dummy = (\identity{\nbig} -
      \projection)(\identity{\nbig}-\projections)\dummy.
    \end{equation*}
    Because $\projection \neq \zerobold$, $\projection \neq \identity{\nbig}$,
    it holds that $\|\projection\|_2 = \|\identity{\nbig}-\projection\|_2$.
    Note also that 
    $\|\projection\|_2 \leq
    \|(\samplemat\orthogonalmat)^{\dagger}\|_2$ is true because
    $\projection =
    \orthogonalmat(\samplemat\orthogonalmat)^{\dagger}\samplemat$,
    which gives \eqref{eq:projectionerror}.
    \hspace*{1em}\endproof}

  {\theorem\label{theorem:DEIM-SNSbound} Let the DEIM-SNS projection matrix be
  $\projection =
  \basisfluxmat(\samplemat\basisfluxmat)^{-1}\samplemat\in\RR{\ndof\times\ndof}$
  and $\projections = \basisfluxmat(\basisfluxmat)^{\dagger}$.
  If the sampling matrix $\samplematNT\in\RR{\ndof\times\nsample}$ is
  constructed using sRRQR factorization of $\orthogonalmat$ or
  $\orthogonalmat_e$, where $\orthogonalmat$ and $\orthogonalmat_e$ are obtained from QR
  factorizations of $\mass\basismatspace$ and $\mass\basismatspaceext$,
  respectively, 
  with a tuning parameter, $\tuningparam$, as in \cite{drmac2018discrete}
  (i.e., $\nsample = \nbasisflux$) or oversampled (i.e.,
  $\ndof \geq \nsample > \nbasisflux$), then
  the DEIM-SNS method with either $\basisfluxmat=\mass\basismatspace$ or 
  $\basisfluxmat=\mass\basismatspaceext$ has a nonlinear term projection error
  bound for a given vector $\flux\in\RR{\ndof}$:
  \begin{equation}\label{eq:DEIM-SNSbound}
    \|(\identity{\ndof}-\projection)\flux\|_2 \leq 
    \sqrt{1+\tuningparam^2\nbasisflux(\ndof-\nbasisflux)}
    \|(\identity{\ndof}-\projections)\flux\|_2.
  \end{equation}
  \proof 
  Combining Lemmas \ref{lemma:projectionbound} with 
  $\nonsingularmat = \mass$, $\orthobasis = \basismatspace$ or
  $\basismatspaceext$ in Section~\ref{sec:DEIM-SNS} and $\nbig=\ndof$, and
  $\nsmall=\nbasisflux$ in Section~\ref{sec:DEIM} 
  and applying sPPQR factorization of
  \cite{drmac2018discrete} to either $\orthogonalmat$ or $\orthogonalmat_e$ 
  and noting that
  $\|(\samplemat\orthogonalmat_e)^{\dagger}\|_2 \leq
  \|(\samplemat\orthogonalmat)^{-1}\|_2$ 
  proves the bound.
    \hspace*{1em}\endproof}

  {\theorem\label{theorem:GNAT-SNSbound} Let the GNAT-SNS projection matrix be
  $\projection =
  \basismatres(\samplemat\basismatres)^{-1}\samplemat\in\RR{\ndof\times\ndof}$
  and $\projections = \basismatres(\basismatres)^{\dagger}$.
  set the sampling matrix $\samplematNT$, 
  If the sampling matrix $\samplematNT\in\RR{\ndof\times\nsample}$ is
  constructed using sRRQR factorization of $\orthogonalmat$ or
  $\orthogonalmat_e$, where $\orthogonalmat$ and $\orthogonalmat_e$ are obtained from QR
  factorizations of $\mass\basismatspace$ and $\mass\basismatspaceext$,
  respectively, 
  with a tuning parameter, $\tuningparam$, as in \cite{drmac2018discrete}
  (i.e., $\nsample = \nbasisres$) or oversampled (i.e.,
  $\ndof \geq \nsample > \nbasisres$), then
  the GNAT-SNS method with either $\basismatres=\mass\basismatspace$ or 
  $\basismatres=\mass\basismatspaceext$ has a nonlinear residual term projection error
  bound for a given vector $\res\in\RR{\ndof}$:
  \begin{equation}\label{eq:GNAT-SNSbound}
    \|(\identity{\ndof}-\projection)\res\|_2 \leq 
    \sqrt{1+\tuningparam^2\nbasisres(\ndof-\nbasisres)} 
    \|(\identity{\ndof}-\projections)\res\|_2.
  \end{equation}
  \proof 
  Combining Lemmas \ref{lemma:projectionbound} with 
  $\nonsingularmat = \mass$, $\orthobasis = \basismatspace$ or
  $\basismatspaceext$ in Section~\ref{sec:GNAT-SNS} and 
  $\nbig=\ndof$, and $\nsmall=\nbasisres$ 
  in Section~\ref{sec:GNAT} 
  and applying sPPQR factorization of
  \cite{drmac2018discrete} to either $\orthogonalmat$ or $\orthogonalmat_e$
  and noting that $\|(\samplemat\orthogonalmat_e)^{\dagger}\|_2 \leq
  \|(\samplemat\orthogonalmat)^{-1}\|_2$
  proves the bound.
    \hspace*{1em}\endproof}

  {\theorem\label{theorem:STGNAT-SNSbound} Let the ST-GNAT-SNS projection matrix be
  $\projection =
  \basismatrest(\samplematst\basismatrest)^{-1}\samplematst\in\RR{\ndof\ntimedof\times\ndof\ntimedof}$
  and $\projections = \basismatrest(\basismatrest)^{\dagger}$.
  If the sampling matrix $\samplematst\in\RR{\ndof\ntimedof\times\nsample}$ is
  constructed using sRRQR factorization of $\orthogonalmatst$ or
  $\orthogonalmatst_e$, where $\orthogonalmatst$ and $\orthogonalmatst_e$ are
  obtained from QR
  factorizations of $\Abe\stbasismat$ and $\Abe\stbasismatext$,
  respectively, 
  with a tuning parameter, $\tuningparam$, as in \cite{drmac2018discrete}
  (i.e., $\nsample =
  \nbasisres$) or oversampled (i.e., $\ndof\ntimedof \geq \nsample > \nbasisres$), then
  the ST-GNAT-SNS method with either $\basismatrest=\Abe\stbasismat$ or 
  $\basismatrest=\Abe\stbasismatext$ has 
  a nonlinear space--time residual term projection error
  bound for a given vector $\resst\in\RR{\ndof\ntimedof}$:
  \begin{equation}\label{eq:ST-GNAT-SNSbound}
    \|(\identity{\ndof\ntimedof}-\projection)\resst\|_2 \leq 
    \sqrt{1+\tuningparam^2\nbasisres(\ndof\ntimedof-\nbasisres)}
    \|(\identity{\ndof\ntimedof}-\projections)\resst\|_2.
  \end{equation}
  \proof 
  Combining Lemmas \ref{lemma:projectionbound} with
  $\nonsingularmat = \Abe$, $\orthobasis = \stbasismat$ or
  $\stbasismatext$ in Section~\ref{sec:ST-GNAT-SNS} and 
  $\nbig=\ndof\ntimedof$, and $\nsmall=\nbasisres$
  in Section~\ref{sec:STGNAT} 
  and applying sPPQR factorization to either $\orthogonalmatst$ or
  $\orthogonalmatst_e$
  and noting that $\|(\samplematst\orthogonalmatst_e)^{\dagger}\|_2 \leq
  \|(\samplematst\orthogonalmatst)^{-1}\|_2$
  proves the bound.
    \hspace*{1em}\endproof}

\section{Numerical Results}\label{sec:numericresults}
In this section, we demonstrate that the SNS methods
reduce the offline computational time 
without losing accuracy of the DEIM, GNAT, and ST-GNAT approaches.
The focus of the numerical experiments is not to show the accuracy and speed-up
of all the model order reduction techniques considered in the paper.
For the benefits of the DEIM, GNAT, and ST-GNAT methods in terms of accuracy and
speed-up, we refer to their original papers:
\cite{chaturantabut2010nonlinear,carlberg2013gnat,choi2019space}.
For the DEIM method, which does not require to collect the nonlinear term
snapshots from other simulations than the corresponding FOM, the only offline
computational time reduction comes from the fact that the SNS methods require
only one compression of the only solution snapshots instead of two compressions
that are requried in the DEIM method. Therefore, the compression computational
time reduction due to the DEIM-SNS method is about a factor of two (e.g., see
Figs.~\ref{fig:forwardEuler_DEIM_nldiff}(c) and
\ref{fig:backwardEuler_DEIM_nldiff}(c)).
For the GNAT and ST-GNAT methods, on the other hand, the offline computational
time reduction due to the SNS methods is large because the GNAT and ST-GNAT
require to collect nonlinear residual term snapshots from their corresponding
LSPG and ST-LSPG simulations for the best performance, but the SNS methods do
not require that. Therefore, the numerical experiments show that 
the offline computational time reduction is from a factor
of three to a hundred (e.g., see
Figs.~\ref{fig:explicit_GNAT_burgers_offlinecost},
\ref{fig:implicit_GNAT_burgers_offlinecost},
\ref{fig:implicit_STGNAT_burgers_offlinetime}, \ref{fig:GNAT_euler}(b), and
\ref{fig:STGNAT_euler_offlinetime}).

We consider three different problems:
a 2D nonlinear diffusion problem is solved in
Section~\ref{sec:diffusion},
a parameterized 1D Burgers' equation is solved in
Section~\ref{sec:burgersProblem}, and 
a parameterized quasi-1D Euler equation is solved in
Section~\ref{sec:eulerProblem}.
The performance of the DEIM and DEIM-SNS approaches are compared
in Section~\ref{sec:diffusion}.
The performance of the GNAT and 
GNAT-SNS methods are compared in Sections~\ref{sec:burgersProblem}
and \ref{sec:eulerProblem}.
The performance of the ST-GNAT and ST-GNAT-SNS methods are
compared in Sections~\ref{sec:burgersProblem} and \ref{sec:eulerProblem}.

The following greedy algorithms for
constructing sample indices are used for each method: 
\begin{itemize}
  \item Algorithm 1 in \cite{chaturantabut2010nonlinear} with the DEIM and DEIM-SNS
    approaches.
  \item Algorithm 3 in \cite{carlberg2013gnat} with the GNAT and GNAT-SNS methods.
  \item Algorithm 1 and 2 in \cite{choi2019space} with the
    ST-GNAT and ST-GNAT-SNS methods.
\end{itemize}
Procedure identifier 1 in Table 1 of \cite{carlberg2013gnat}
is used for residual snapshot-collection procedures for GNAT.
The accuracy of any ROM solution
$\solapproxFuncOnlyParam$ is assessed from its
mean squared state-space error:
\begin{equation} 
  \mseError  = \left.
  \sqrt{\sum_{n=1}^\ntimedof\|\solapproxFuncArg{n}-\solFuncArg{n}\|_2^2}
  \middle/ 
  \sqrt{\sum_{n=1}^\ntimedof\|\solFuncArg{n}\|_2^2} \right..
\end{equation}
We measure the computational offline cost in terms of the wall time.  All
timings with the GNAT, ST-GNAT, GNAT-SNS, and ST-GNAT-SNS methods are obtained by
performing calculations on an Intel Core i7 CPU @ $2.5$ GHz, $16$ GB $1\,600$ MHz
DDR3 using the modified version of \verb+MORTestbed+\cite{zahr2010comparison} in MATLAB.  
All the timings with the DEIM and DEIM-SNS approaches are obtained by performing
calculations on the Quartz cluster at Lawrence Livermore National Laboratory
using MFEM-based reduced order model house code \cite{mfem-library}. Quartz has
$2\,634$ nodes, each with an Intel Xeon E5-2695 with $36$ cores operating at
$2.1$ Ghz and $128$ Gb RAM memory. The time units of
Figures~\ref{fig:explicit_GNAT_burgers_offlinecost},
\ref{fig:implicit_GNAT_burgers_offlinecost},
\ref{fig:implicit_STGNAT_burgers_offlinetime}, and
\ref{fig:STGNAT_euler_offlinetime} are seconds.


\subsection{Nonlinear diffusion equation}\label{sec:diffusion}
\begin{figure}[htbp]
  \centering
  \subfigure[initial temperature and mesh]{
  \includegraphics[width=0.2\textwidth]{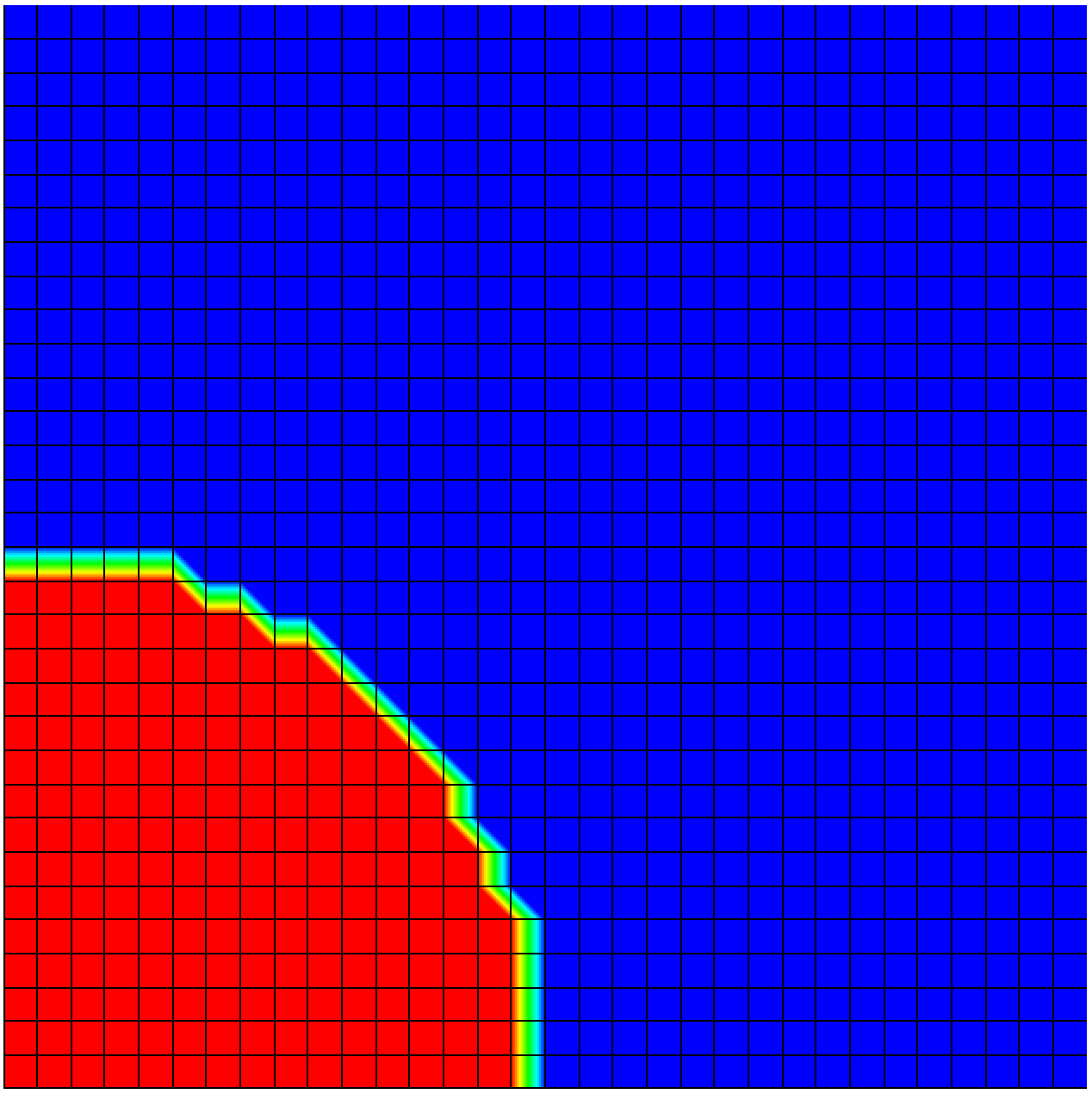}
}~~~~~~~
\subfigure[temperature legend]{
  \includegraphics[width=0.1\textwidth]{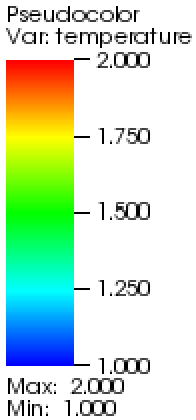}
}
  \caption{Initial temperature distribution, mesh, and legend}
  \label{fig:initial_temperature}
\end{figure}

We now consider a parameterized 2D nonlinear diffusion equation associated 
with the problem of time dependent nonlinear heat conduction.
The problem corresponds to the following initial boundary value problem 
on the unit square and $t\in[0,T]$: 
 \begin{align}\label{eq:diff_eq}
  \begin{split}
      \frac{\partial \temperature}{\partial t} = 
           \nabla \cdot(\kappa + \alpha \temperature)\nabla \temperature,
      \quad \forall (x, y)\in D = [0,1] m\times[0,1] m,
      \quad \forall t\in[0,\totaltime]\\
  \end{split}
 \end{align}
 where $\temperature((x,y),t)\in H^1(D)$ denotes the space and time dependent temperature
 function
 with $\temperature:D\times [0,1] \rightarrow \RR{}$ implicitly defined as the
 solution to Eq.~\eqref{eq:diff_eq}.
The diffusivity depends linearly on $\temperature$ with coefficients,
$\kappa = 0.5$ $m^2/s$ and $\alpha = 0.01$ $m^2/(s\cdot K)$.
Zero temperature gradient boundaries are employed and the simulation is initialized 
by a step function defined on a quarter circle given by: 
\begin{equation}
  \frac{\partial \temperature}{\partial n} = 0, \quad \text{on} 
  \quad \Gamma = \left\{(x,y) |   
      x\in\left\{ 0, 1 \right\},
      y\in\left\{ 0, 1 \right\}
  \right\}
\end{equation}
\begin{equation}
    \temperature(x,y,0;\mu) = \left\{
                  \begin{array}{ll}
                    2 \quad \text{if } x^{2}+y^{2} \leq 0.5^{2} \\
                    1 \quad \text{if otherwise}
                  \end{array}
                  \right. 
\end{equation}
After applying a linear finite-element spatial discretization with $\nspacedof = 1089$
($32$ elements at each side. See Fig.~\ref{fig:initial_temperature} for the mesh).
Eq.~\eqref{eq:diff_eq} lead to an initial-value ODE problem consistent with
Eq.~\eqref{eq:fom} with $\mass$ being a volume matrix.
For time discretization, 
the forward and backward Euler schemes are applied with a uniform time step.
The solution basis dimension is set $\nbasisspace = 20$.

\subsubsection{DEIM versus DEIM-SNS}\label{sec:DEIMvsDEIMSNS}
\begin{figure}[htbp]
    \centering
    \subfigure[FOM temperature]{
    \includegraphics[width=0.2\textwidth]{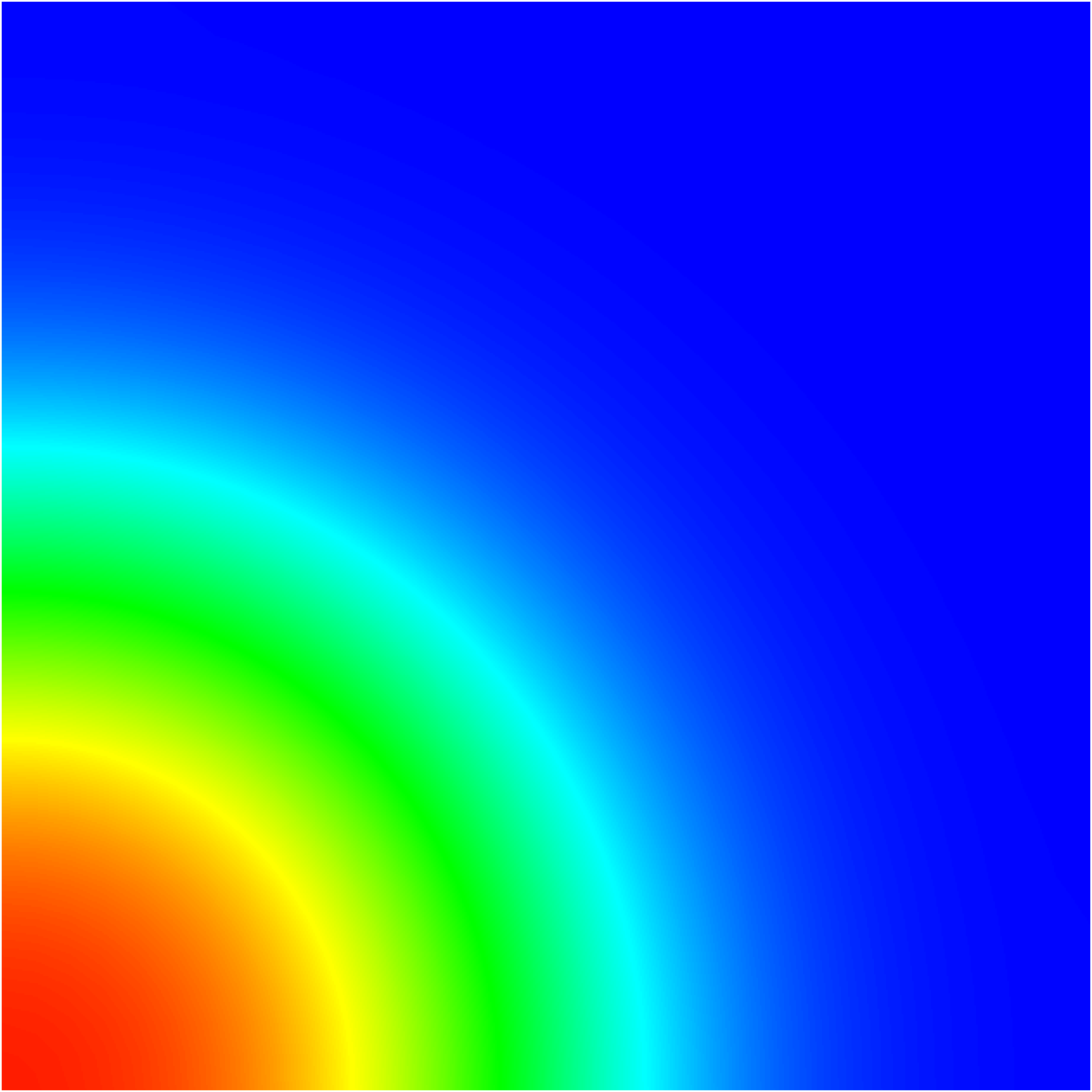}
  }~~~~~~~
  \subfigure[DEIM temperature]{
    \includegraphics[width=0.2\textwidth]{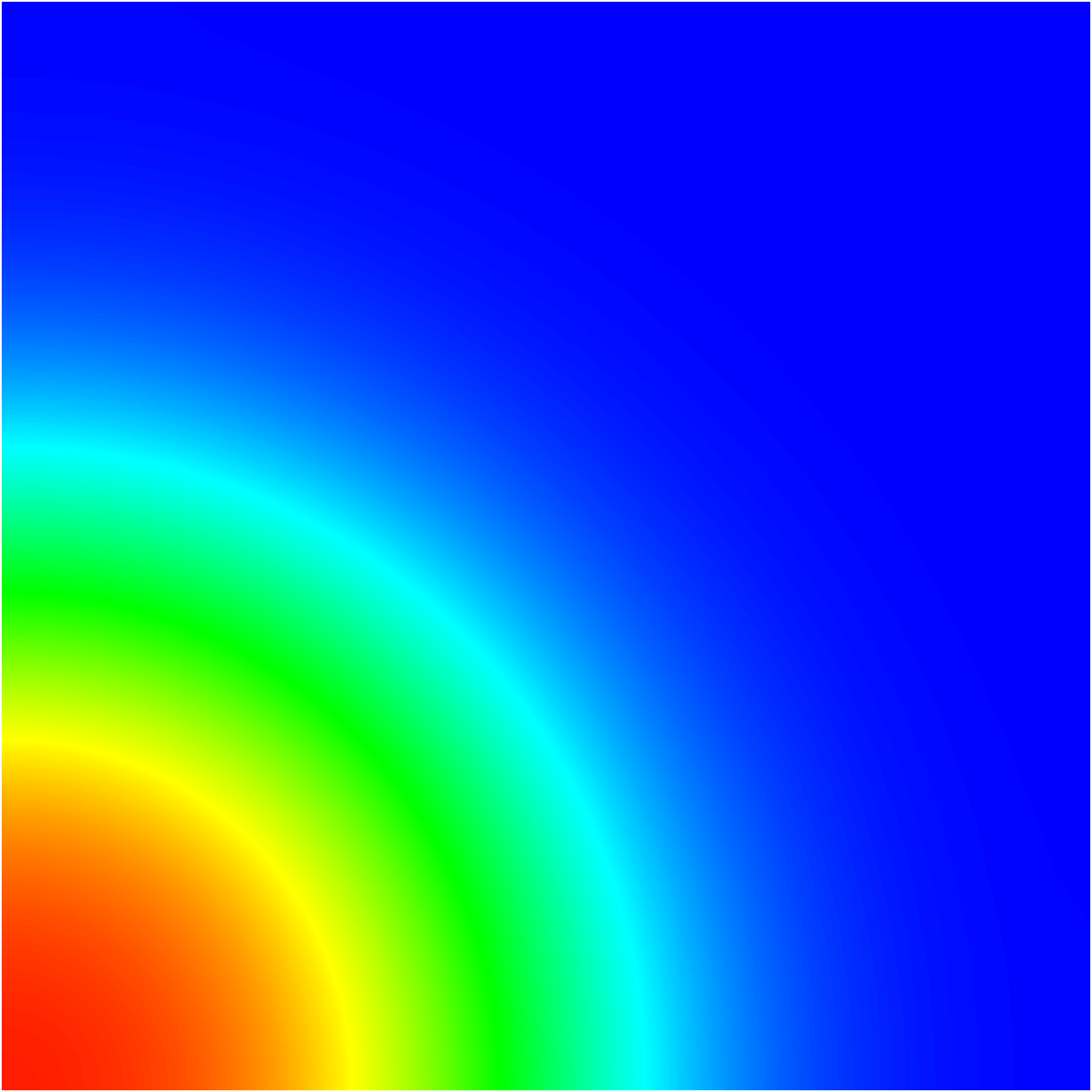}
  }~~~~~~~
  \subfigure[DEIM-SNS temperature]{
    \includegraphics[width=0.2\textwidth]{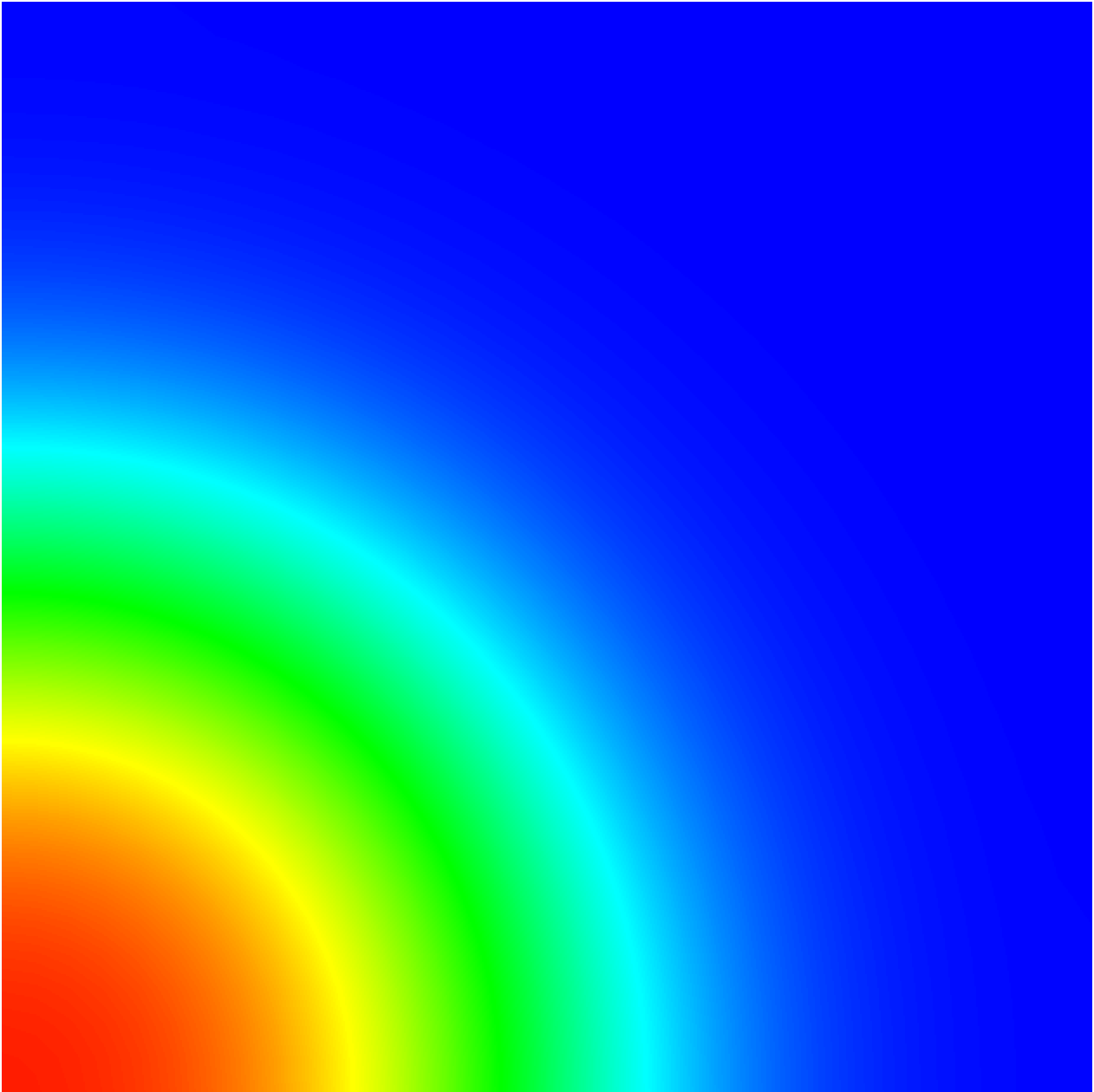}
  }
    \caption{Temperature distribution using the forward Euler time integrator.
    For the DEIM and DEIM-SNS approaches, $\nbasisspace = 20$ and $\nbasisres = 20$
  are used.}
    \label{fig:modelsolutions_DEIM_explicit}
  \end{figure}

\begin{figure}[htbp]
    \centering
    \subfigure[relative error]{
    \includegraphics[width=0.3\textwidth]{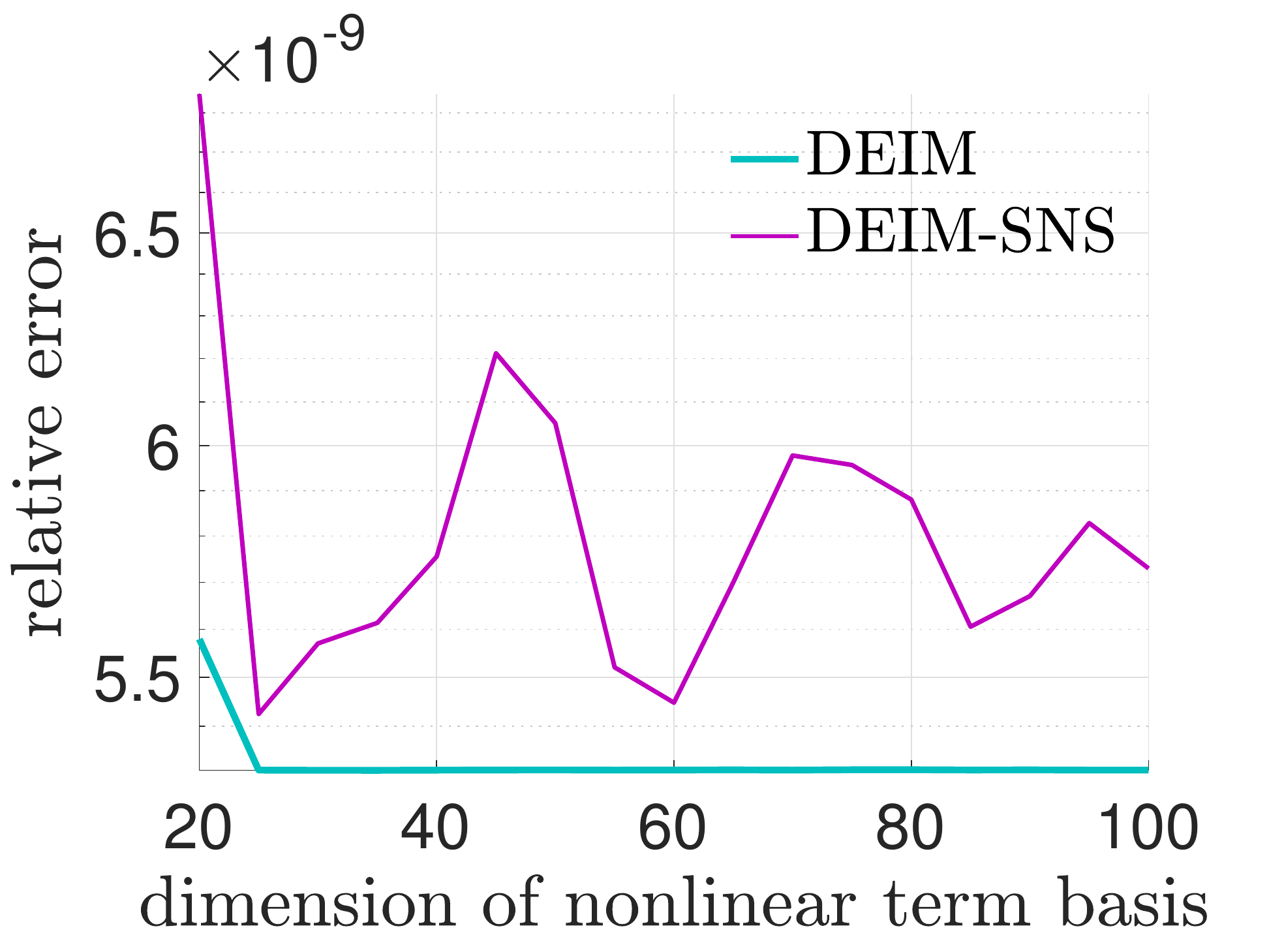}
  }~~~~~~~
  \subfigure[offline time]{
    \includegraphics[width=0.3\textwidth]{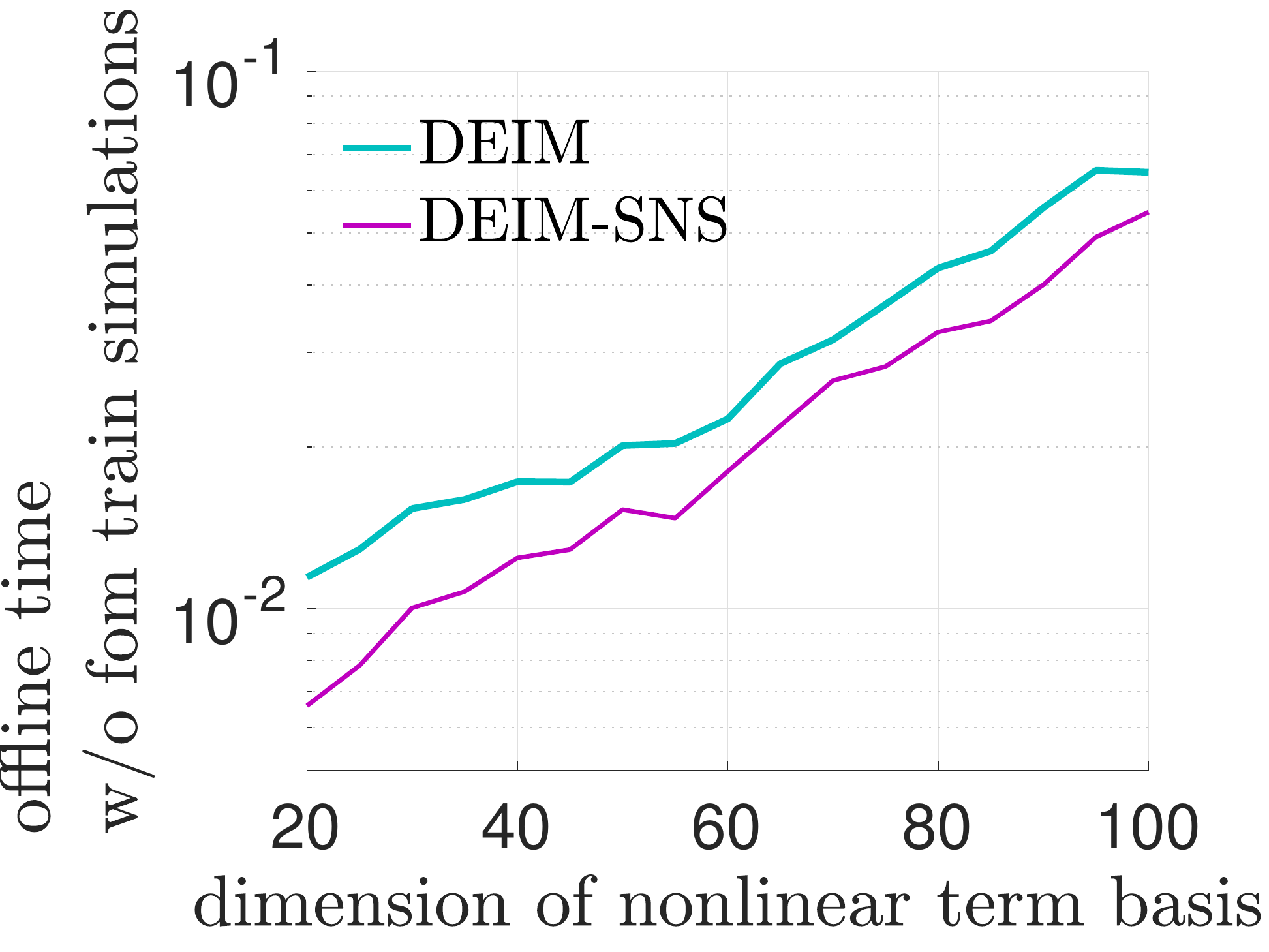}
  }~~~~~~~
  \subfigure[SVD time]{
    \includegraphics[width=0.3\textwidth]{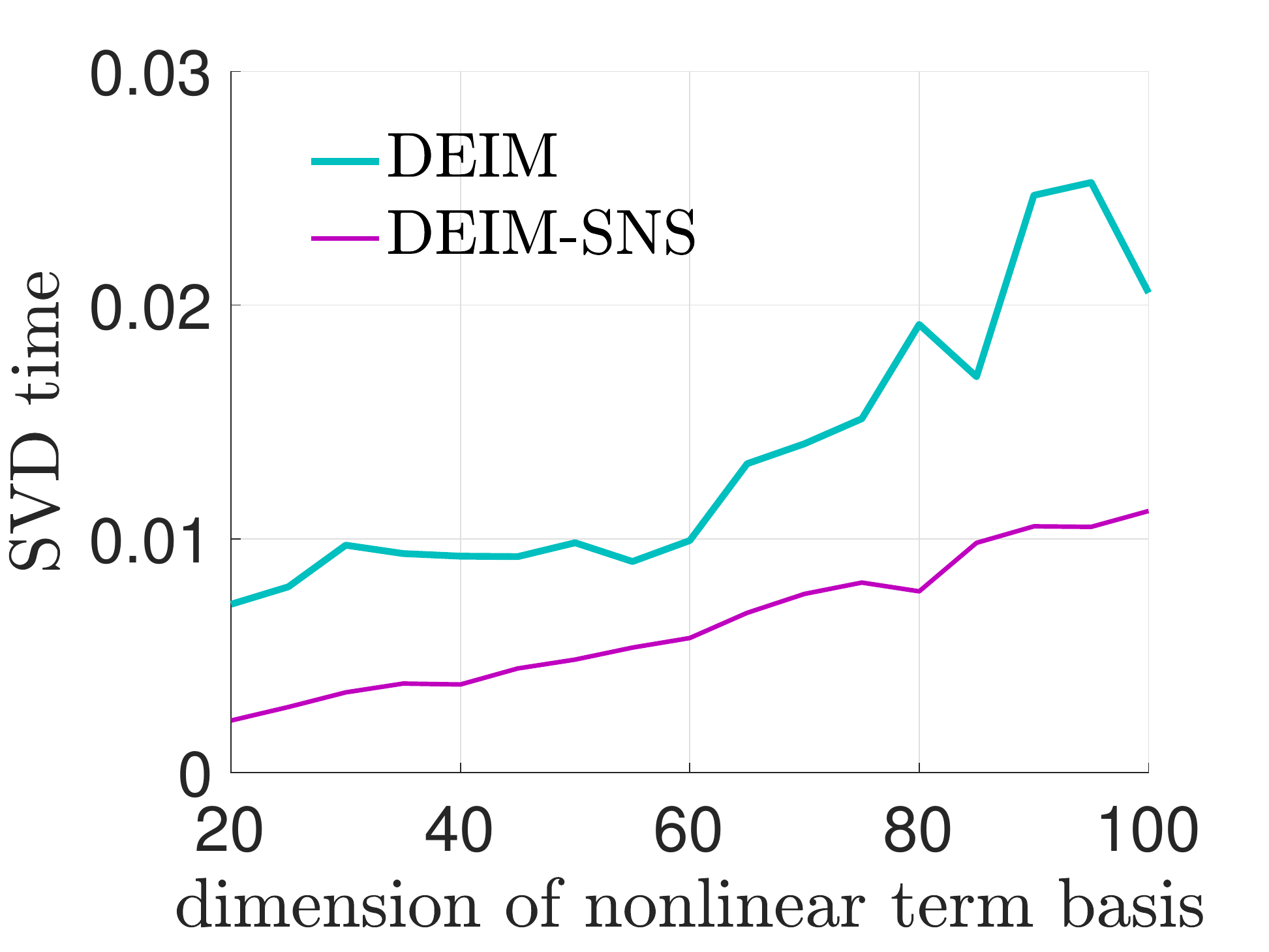}
  }~~~~~~~
    \caption{Relative errors and offline time with the forward Euler time
    integrator, a number of time steps ($\ntimedof = 100$).}
    \label{fig:forwardEuler_DEIM_nldiff}
  \end{figure}

\begin{figure}[htbp]
    \centering
    \subfigure[relative error]{
    \includegraphics[width=0.3\textwidth]{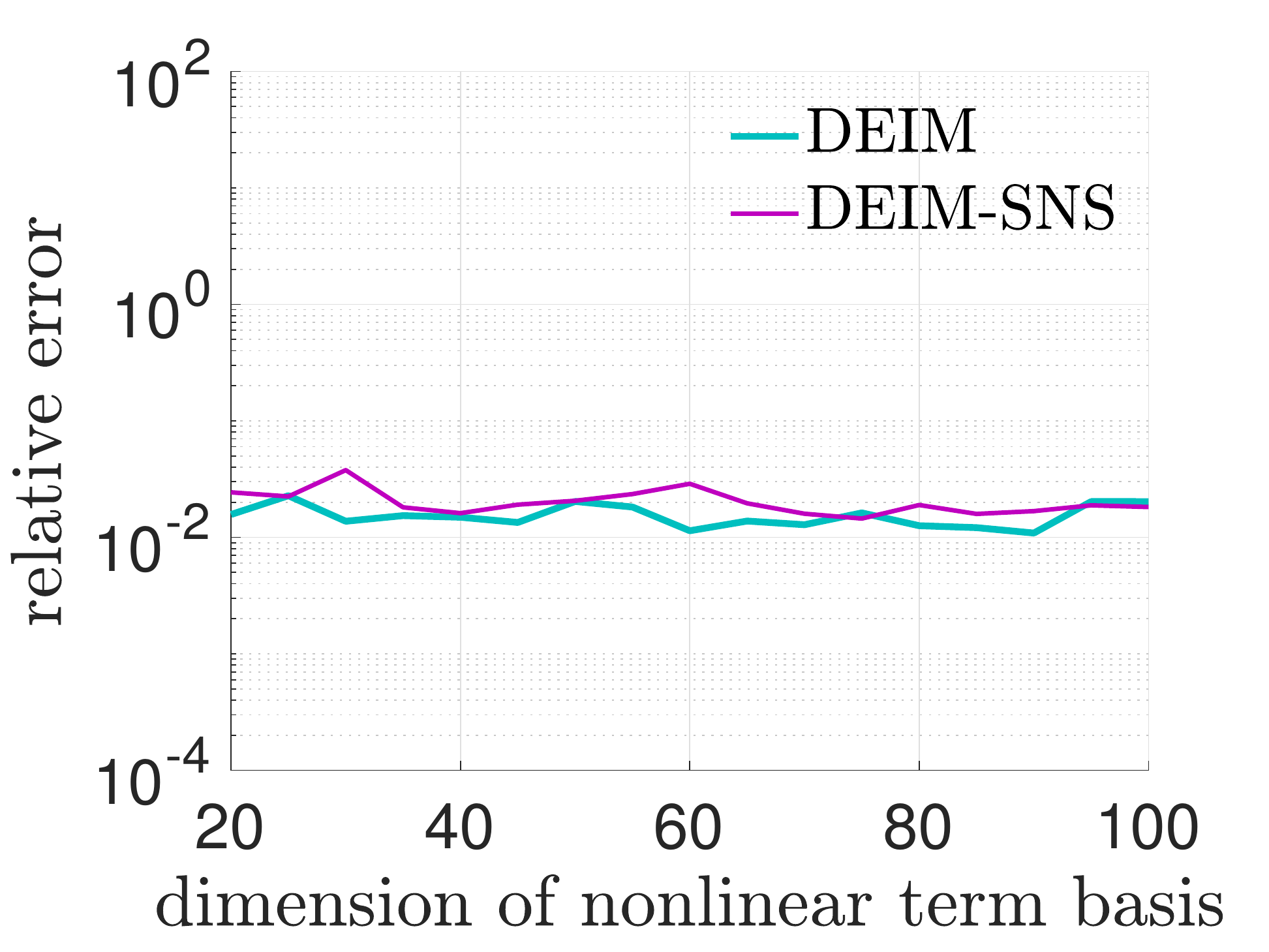}
  }~~~~~~~
  \subfigure[offline time]{
    \includegraphics[width=0.3\textwidth]{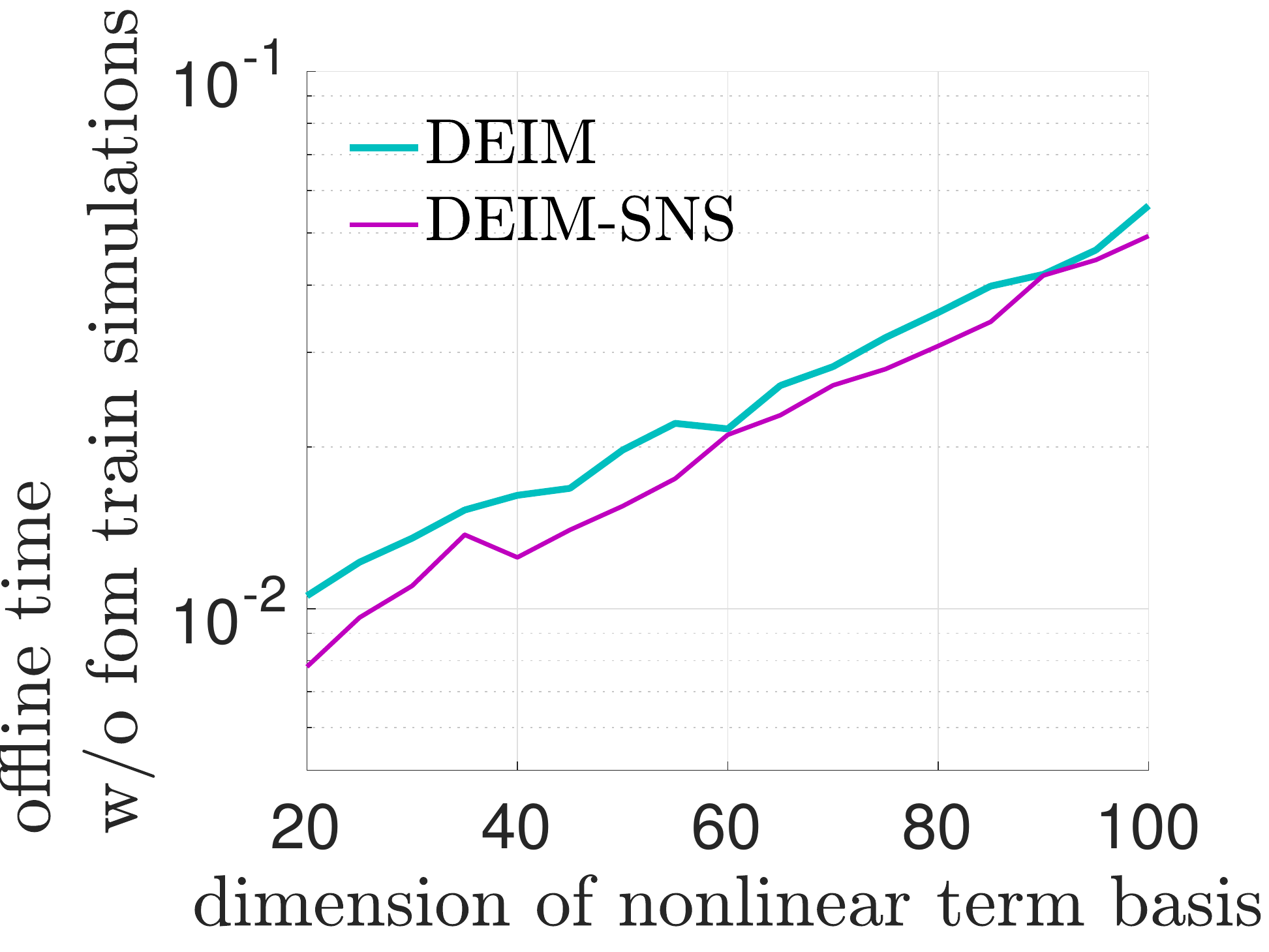}
  }~~~~~~~
  \subfigure[SVD time]{
    \includegraphics[width=0.3\textwidth]{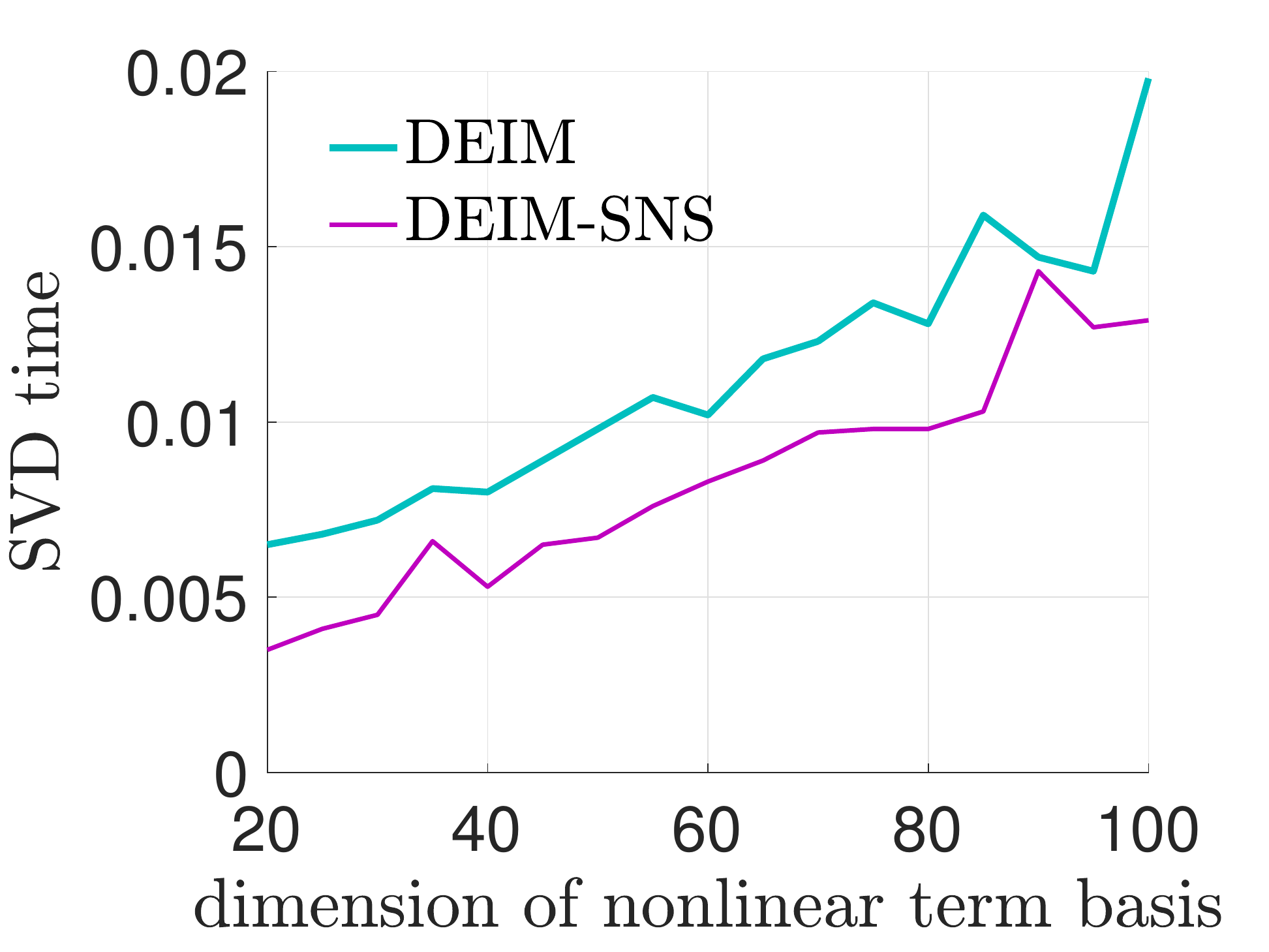}
  }~~~~~~~
    \caption{Relative errors and offline time with the backward Euler time
    integrator, a number of time steps ($\ntimedof = 100$).}
    \label{fig:backwardEuler_DEIM_nldiff}
  \end{figure}

 The DEIM and DEIM-SNS approaches are compared numerically. For this parabolic
 problem, the DEIM and DEIM-SNS approach tries to reproduce the solution of the
 corresponding high fidelity model with the same problem parameter settings. For
 a parametric case, where the DEIM and DEIM-SNS approaches are trained with a number of
 sample points in a parameter space and are used to predict the solution of a
 new parameter point, is considered for the hyperbolic problems in
 Sections~\ref{sec:burgersProblem} and \ref{sec:eulerProblem}.
 Figs.~\ref{fig:modelsolutions_DEIM_explicit}, \ref{fig:forwardEuler_DEIM_nldiff} and 
 \ref{fig:backwardEuler_DEIM_nldiff} are generated by setting $T = 0.01$ s and
 $\dt = 1.0\times10^{-4}$ s, leading to $\ntimedof = 100$.
 The relative error and the offline time are plotted as
 the dimension of nonlinear term basis $\nbasisflux$ increases.
  For DEIM-SNS, $\basisfluxmat = \mass\basismatspace$ is used if $\nbasisflux = \nbasisspace$
  while $\basisfluxmat = \mass\basismatspaceext$ is used for $\nbasisflux >
  \nbasisspace$, where $\mass$ is a volume matrix.


  In Fig.~\ref{fig:forwardEuler_DEIM_nldiff}(a), 
  the relative errors of both the DEIM and DEIM-SNS approaches are plotted as the
  dimension of nonlinear term basis increases from 20 to 100 
  by 5 with the forward Euler time integrator.
  The figure shows that DEIM-SNS is comparable to DEIM 
  in terms of accuracy. Note that the order of relative
  errors both with DEIM and DEIM-SNS is $10^{-9}$ 
  for the whole range of the nonlinear term basis dimension considered here.
  This implies that setting the dimension of the nonlinear term basis
  $\nbasisflux$ as small as the dimension of the solution basis $\nbasisspace$
  is sufficient to achieve a good accuracy.  

  Fig.~\ref{fig:forwardEuler_DEIM_nldiff}(b) shows the offline times 
  required by DEIM and DEIM-SNS. 
  The offline time of the DEIM approach includes the time of two SVDs for the solution 
  and nonlinear term bases construction and the time of constructing sample indices. 
  The offline time of DEIM-SNS includes the time of `one' SVD for
  the solution basis and
  nonlinear term basis and the time of constructing sample indices. 
  Because DEIM-SNS requires only one SVD, 
  the offline time of the DEIM-SNS approach is
  less than the one of the DEIM approach.
  This fact is shown more clearly in Fig.~\ref{fig:forwardEuler_DEIM_nldiff}(c)
  that shows the SVD times only, where DEIM-SNS achieves a speed-up of around
  two with respect to DEIM. 
  This excludes the time of constructing sample indices from
  Fig.~\ref{fig:forwardEuler_DEIM_nldiff}(b).

  In Fig.~\ref{fig:backwardEuler_DEIM_nldiff}(a), the relative errors 
  of both the DEIM and DEIM-SNS approaches are plotted as the
  dimension of the nonlinear term basis increases from 20 to 100 by 5 with the backward Euler time integrator.
  The figure shows that the DEIM-SNS approach is comparable to the DEIM approach in terms of accuracy. 
  The order of relative errors of both DEIM and DEIM-SNS is $10^{-2}$ 
  for the whole range of the nonlinear term basis dimensions considered here.

  Fig.~\ref{fig:backwardEuler_DEIM_nldiff}(b) shows 
  the offline time required by DEIM and DEIM-SNS. 
  The offline time of the DEIM approach includes the time of two SVDs for the solution
  and nonlinear term bases construction and the time of constructing sample indices. 
  The offline time of the DEIM-SNS approach includes the time of `one' SVD for
  the solution and
  nonlinear term bases construction and the time of constructing sample indices. 
  Because DEIM-SNS requires only one SVD, 
  the offline time of DEIM-SNS is less than the one of DEIM.
  This fact is shown more clearly in Fig.~\ref{fig:backwardEuler_DEIM_nldiff}(c)
  that shows the SVD times only.
  This excludes the time of constructing sample indices from
  Fig.~\ref{fig:backwardEuler_DEIM_nldiff}(b). 

\subsection{Parameterized 1D Burgers' equation}\label{sec:burgersProblem}
We first consider the parameterized inviscid Burgers' equation described in
Ref.~\cite{rewienski2003trajectory}, which corresponds to the following initial boundary value
problem for $x\in[0,1]$ m and $t\in[0,\totaltime]$ with
$\totaltime = 0.5$ s:
 \begin{align}\label{eq:burgers_eq} 
\begin{split}
   \frac{\partial w(x,t;\mu)}{\partial t} + \frac{\partial f(w(x,t;\mu))}{\partial x} &= 0.02e^{\mu_2 x}, 
   \quad \forall x\in[0,1], \quad \forall t\in[0,\totaltime]\\
   w(0,t;\param) &= \mu_1, \quad \forall t\in[0,\totaltime] \\
   w(x,0) &= 1, \quad \forall x\in[0,1]
\end{split}
 \end{align} 
 where $w: [0,1] \times [0,\totaltime]\times\paramDomain \rightarrow
 \RR{}$ is a conserved quantity and the $\nparam=2$ parameters comprise
 the left boundary value and source-term coefficient with 
$\param\equiv(\mu_1,\mu_2)\in\paramDomain = [1.2,1.5]\times[0.02,0.025]$.

After applying Godunov's scheme for spatial discretization, Eqs.~\eqref{eq:burgers_eq} lead to a parameterized
initial-value ODE problem consistent with Eq.~\eqref{eq:fom} with $\mass$ being an identity matrix.  For this problem,
all ROMs employ a training set $\paramDomainTrain= \{ 1.2,1.3,1.4,1.5 \}\times\{0.02,0.025\}$ such that $\ntrain=8$ at
which the FOM is solved. Then, the target parameter, $\param = (1.45,0.0201)$,
is pursued.

\subsubsection{GNAT-SNS versus GNAT}\label{sec:burgers-gnat}
  For the spatial ROMs the domain is discretized with $1\,000$ control volumes, 
  for $\nspacedof=1\,000$ spatial degrees of freedom.
  We employ $\ntimedof = 2\,000$, leading to a uniform time step of $\dt = 2.5\times10^{-4}$.
  The solution basis dimension of $\nbasisspace = 100$ is used. 
  The relative error and the offline time are plotted as
  the dimension of nonlinear residual term basis $\nbasisres$ increases.
  For the GNAT-SNS method, $\basismatres = \mass\basismatspace$ is set if
  $\nbasisres = \nbasisspace$,
  while $\basismatres = \mass\basismatspaceext$ is set for $\nbasisres > \nbasisspace$.

  \begin{figure}[htbp]
    \centering
    \subfigure[forward Euler time integrator]{
    \includegraphics[width=0.4\textwidth]{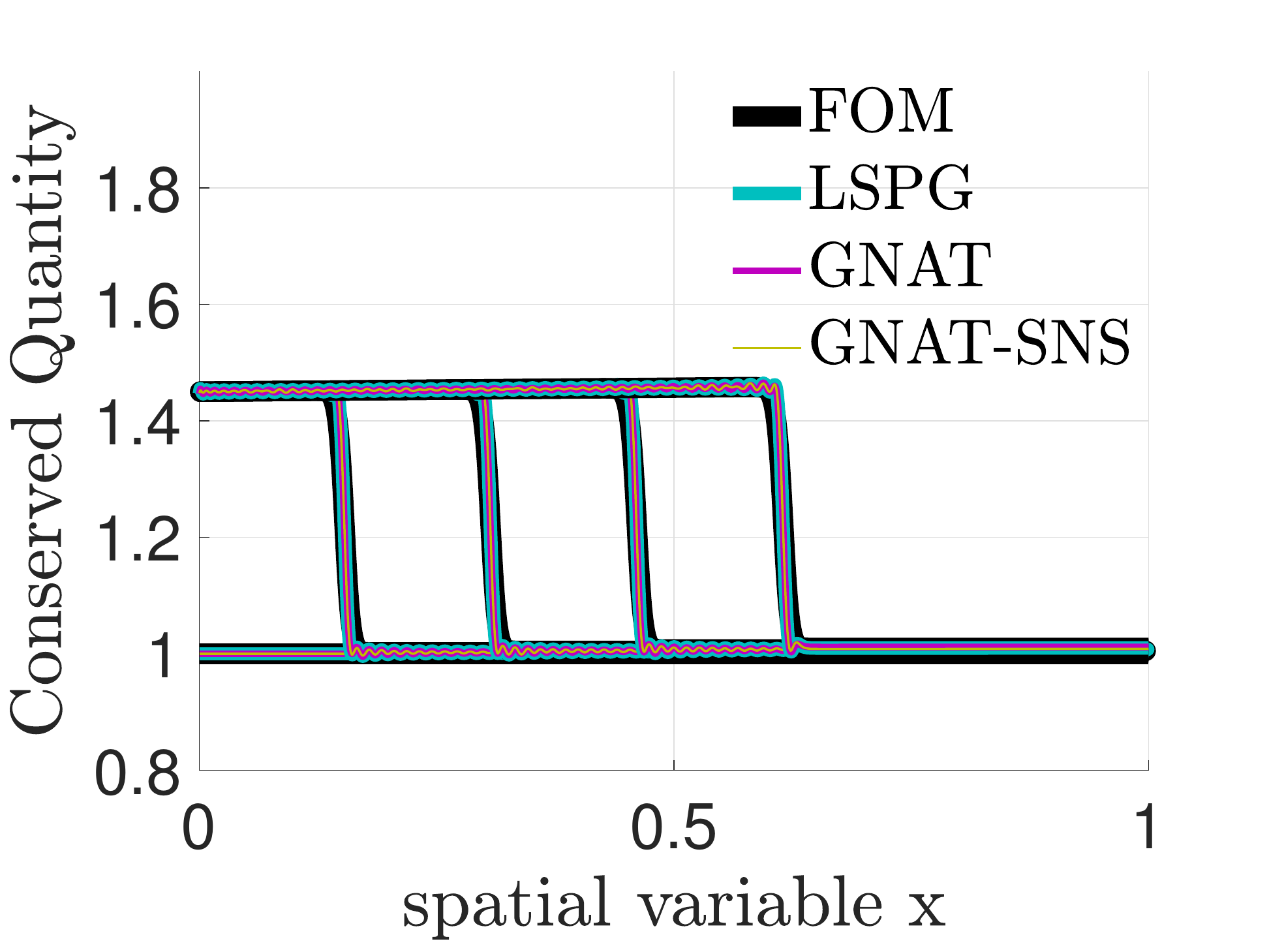}
  }~~~~~~~
  \subfigure[backward Euler time integrator]{
    \includegraphics[width=0.4\textwidth]{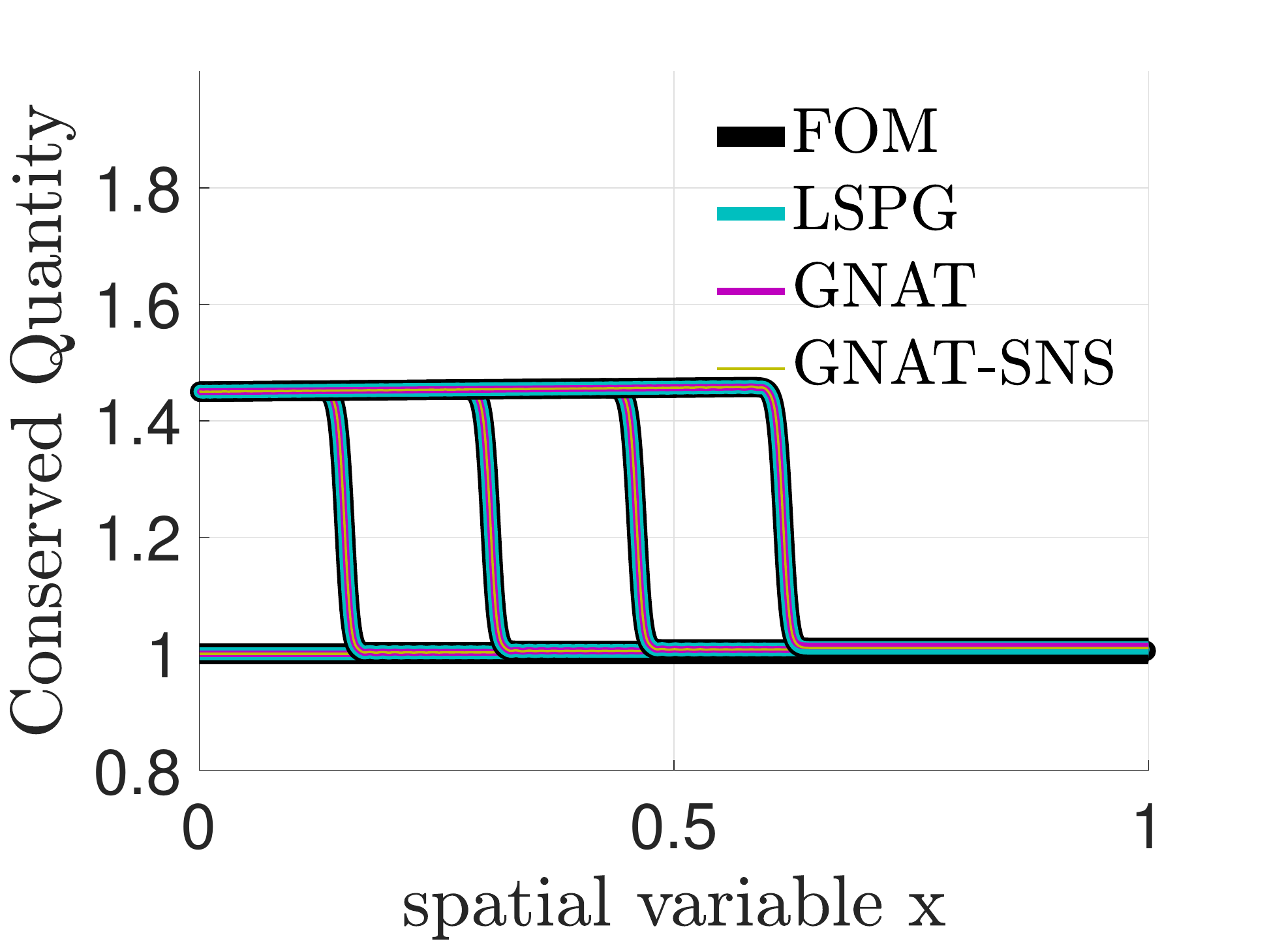}
  }
  \caption{Solution snapshots at $t\in\{0,0.125,0.25,0.375,0.5\}$, $\nbasisres =
  100$, $\nressample = 300$}
    \label{fig:GNAT_burgers_snapshots}
  \end{figure}
  Figs.~\ref{fig:GNAT_burgers_snapshots} compare the solution snapshots of several
  methods: the LSPG, GNAT, and GNAT-SNS methods with the FOM solution snapshots.
  Fig.~\ref{fig:GNAT_burgers_snapshots}(a) is generated with the forward Euler
  time integrator, while 
  Fig.~\ref{fig:GNAT_burgers_snapshots}(b) is generated with the backward Euler 
  time integrator. 
  All the methods are able to generate almost the same solutions as the FOM
  solutions.  

  \begin{figure}[htbp]
    \centering
    \subfigure[forward Euler]{
    \includegraphics[width=0.3\textwidth]{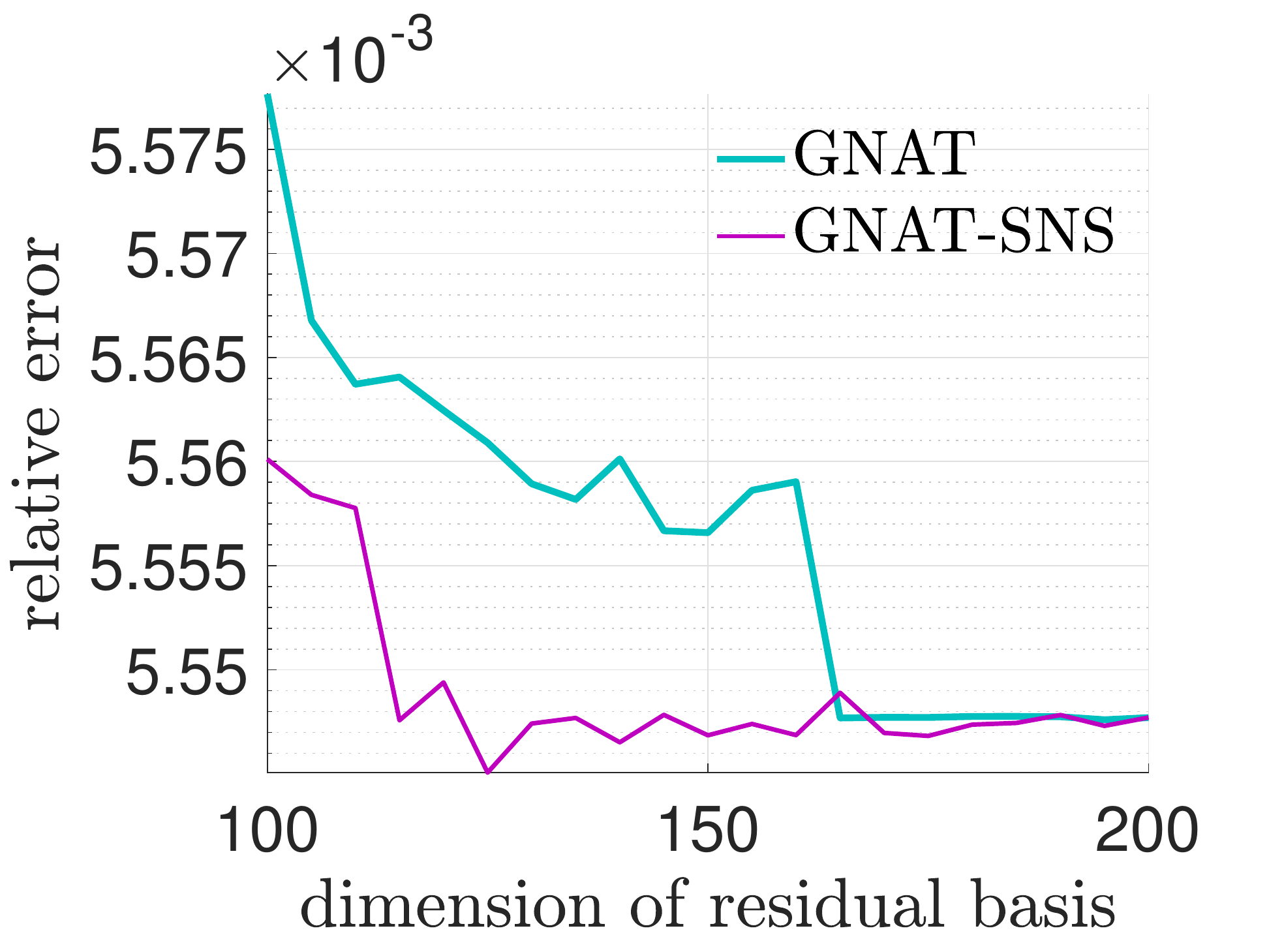}
  }~~~~~~~
  \subfigure[Adams--Bashforth]{
    \includegraphics[width=0.3\textwidth]{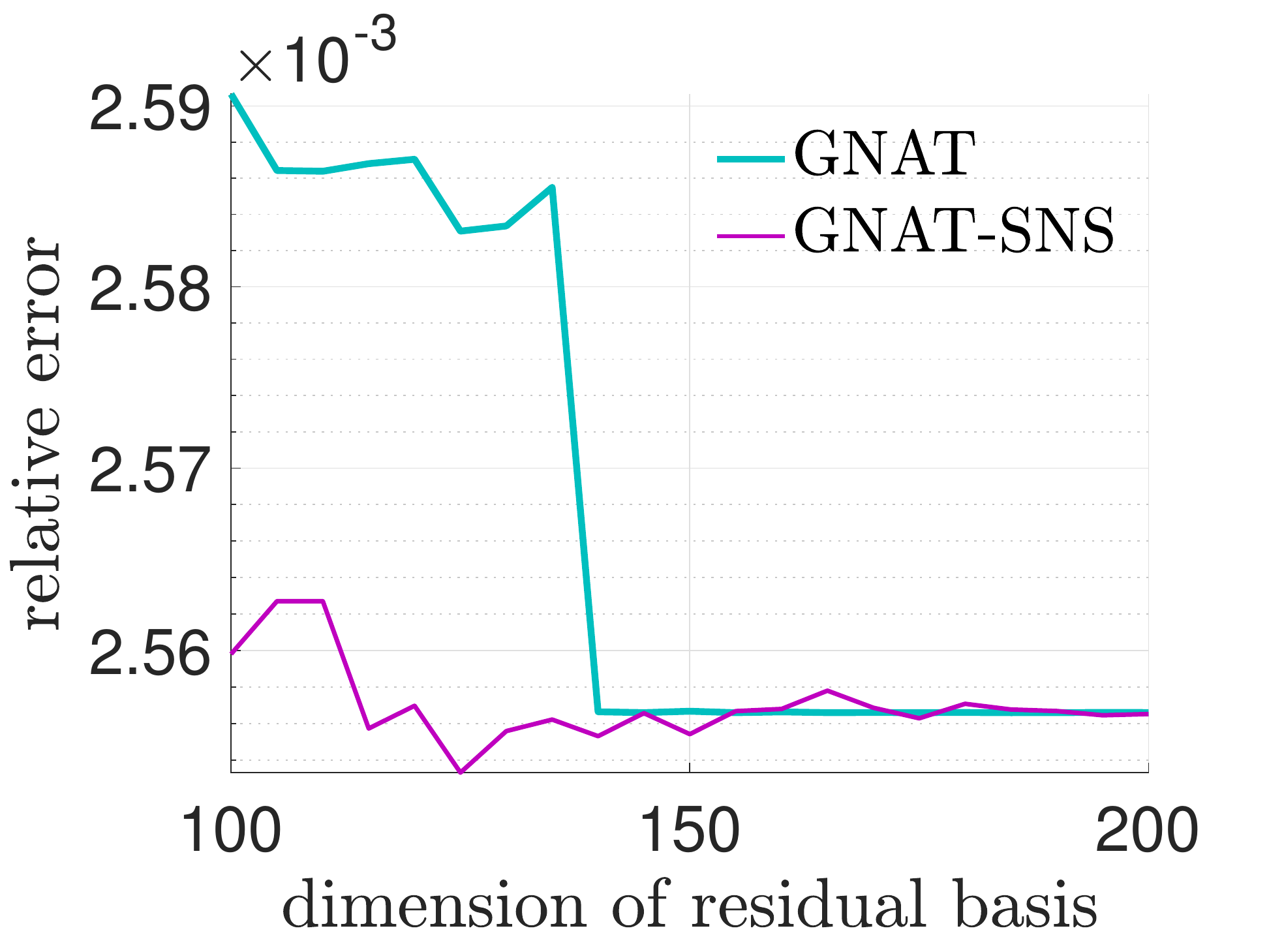}
  }~~~~~~~
  \subfigure[midpoint Runge--Kutta]{
    \includegraphics[width=0.3\textwidth]{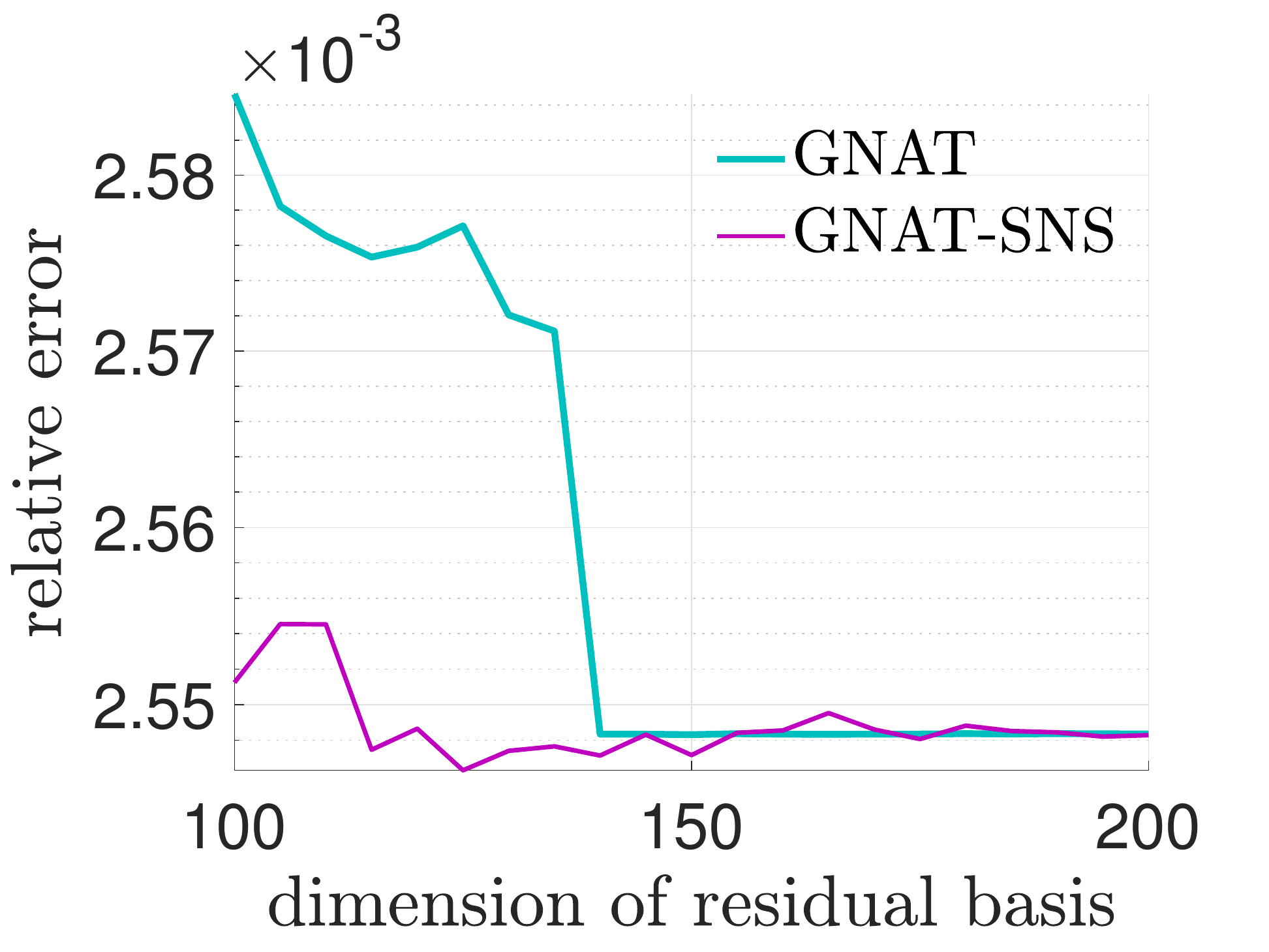}
  }
    \caption{Relative errors with respect to FOM solution for explicit time integrators}
    \label{fig:explicit_GNAT_burgers}
  \end{figure}
  In Figs.~\ref{fig:explicit_GNAT_burgers}, the relative errors of both the GNAT
  and GNAT-SNS methods are plotted as the dimension of the residual basis
  ($\nbasisres$) increases
  from 100 to 200 by 5 with a fixed number of sample indices, $\nressample =
  300$, with the three
  different explicit time integrators: the forward Euler, the Adams--Bashforth,
  and the midpoint Runge--Kutta time integrators.  The figures show that the GNAT-SNS
  method is comparable to the GNAT method in terms of accuracy; the order of
  relative errors are $10^{-3}$ for the whole range of the residual basis
  dimensions considered here.

  \begin{figure}[htbp]
    \centering
    \subfigure[backward Euler]{
    \includegraphics[width=0.3\textwidth]{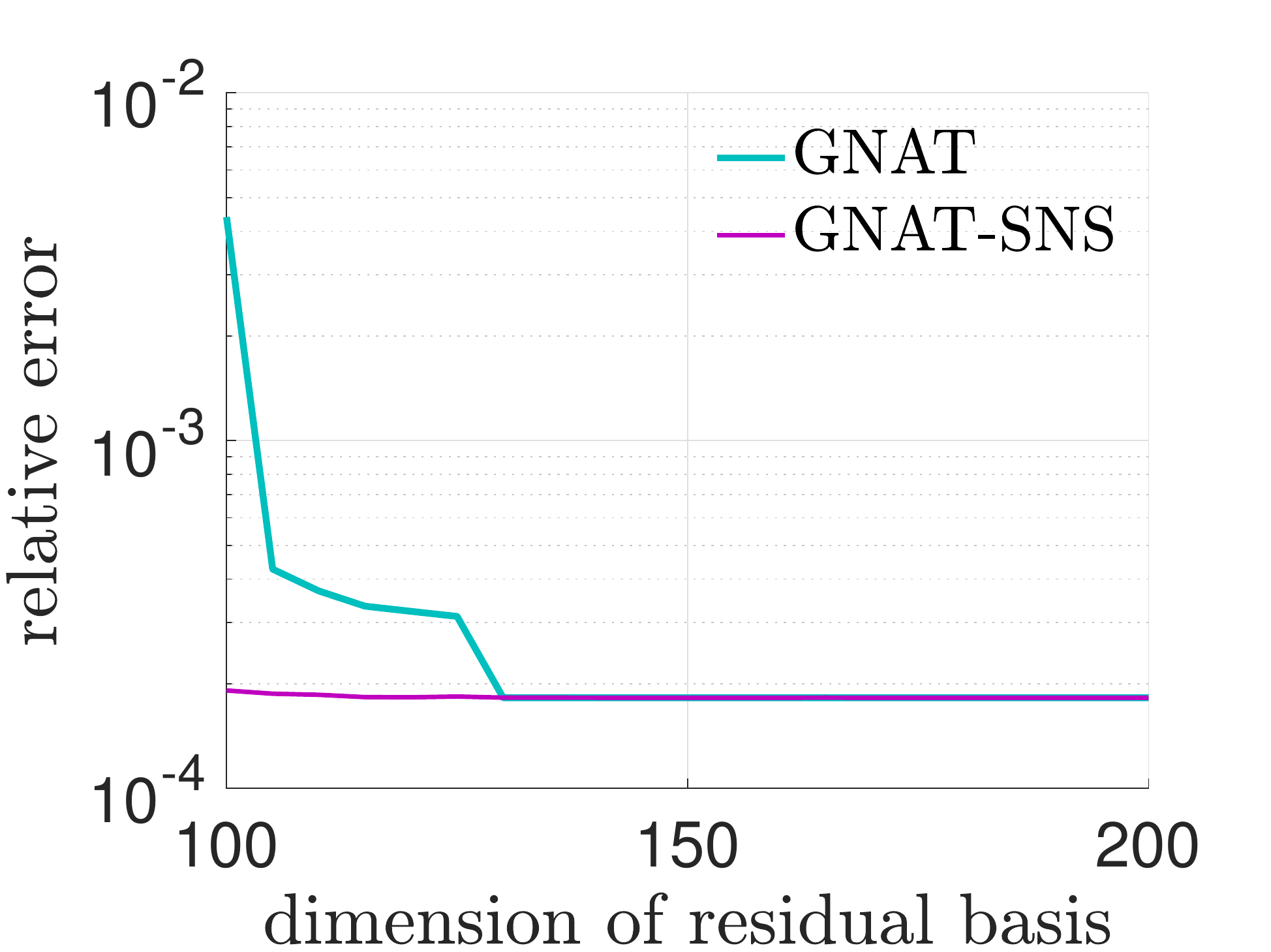}
  }~~~~~~~
  \subfigure[Adams--Moulton]{
    \includegraphics[width=0.3\textwidth]{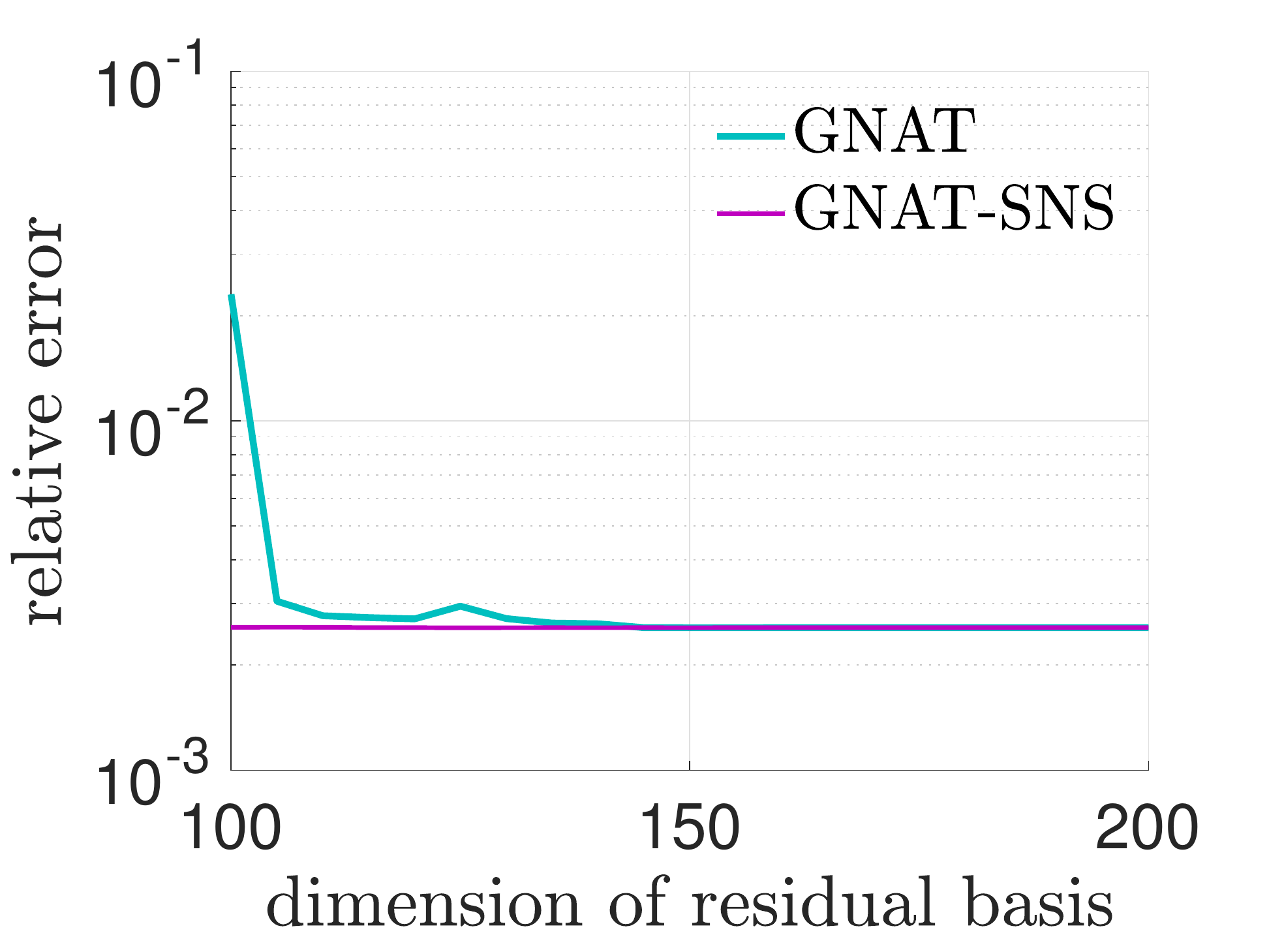}
  }~~~~~~~
  \subfigure[backward differentiation formula]{
    \includegraphics[width=0.3\textwidth]{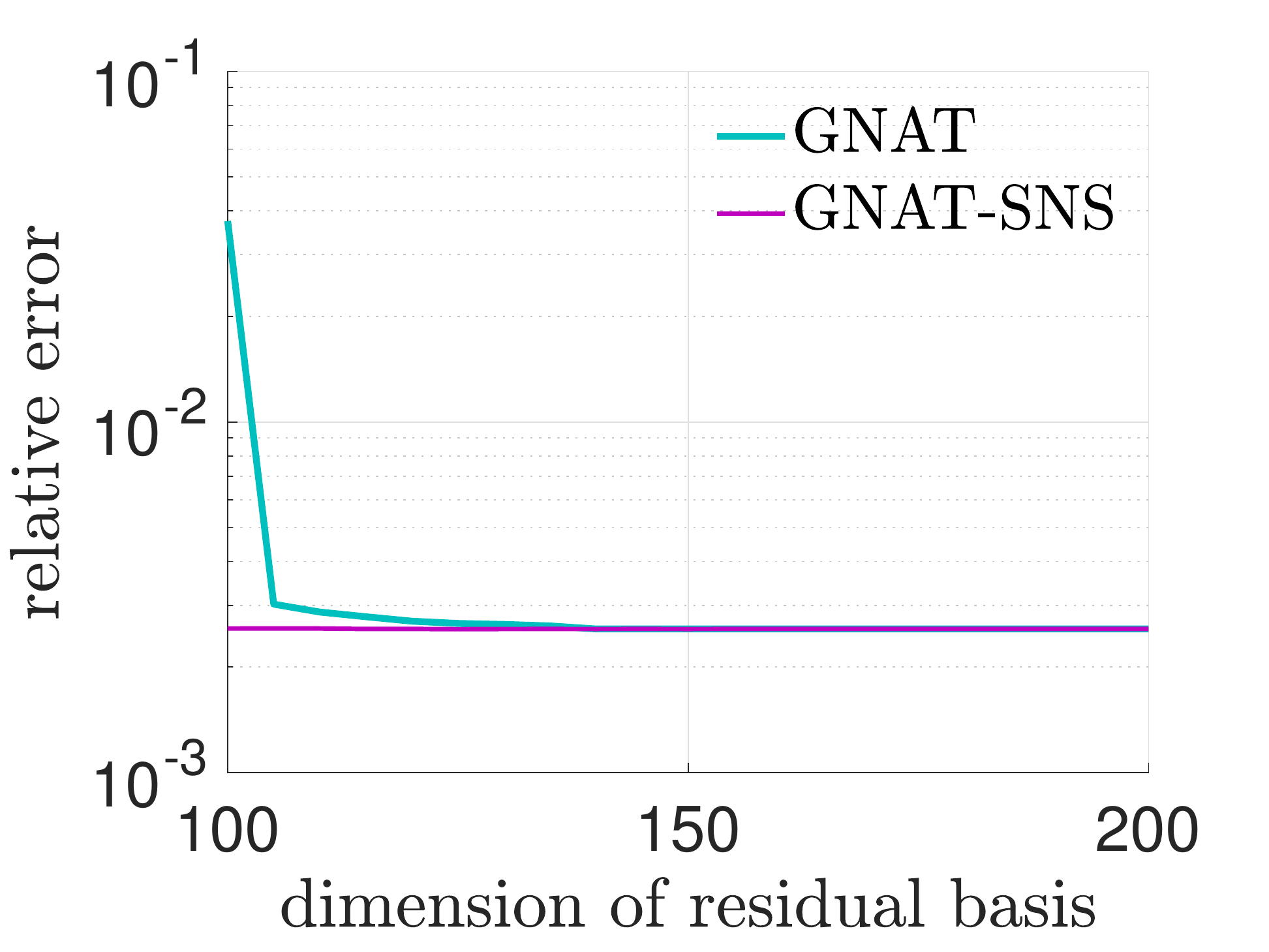}
  }
      \caption{Relative errors with respect to FOM solution for implicit time integrators}
    \label{fig:implicit_GNAT_burgers}
  \end{figure}
  In Figs.~\ref{fig:implicit_GNAT_burgers}, the relative errors of both the
  GNAT and GNAT-SNS methods are plotted as the dimension of residual basis
  ($\nbasisres$)
  increases from 100 to 200 by 5 with a fixed number of sample indices,
  $\nressample = 300$, with the
  three different implicit time integrators: the backward Euler, the
  Adams--Moulton, and the BDF time integrators. The figures show that the GNAT-SNS
  method produces results with better accuracy than the GNAT method when the residual
  basis dimensions are between 100 and 120 (i.e., $\nbasisres \approx
  \nbasisspace$). In fact, the GNAT-SNS achieves an accuracy as good as the LSPG
  can achieve. We are not sure why this is so at this time, but it would be interesting
  to investigate the cause of it in the future work.

  \begin{figure}[htbp]
    \centering
    \subfigure[forward Euler]{
      \includegraphics[width=0.3\textwidth]{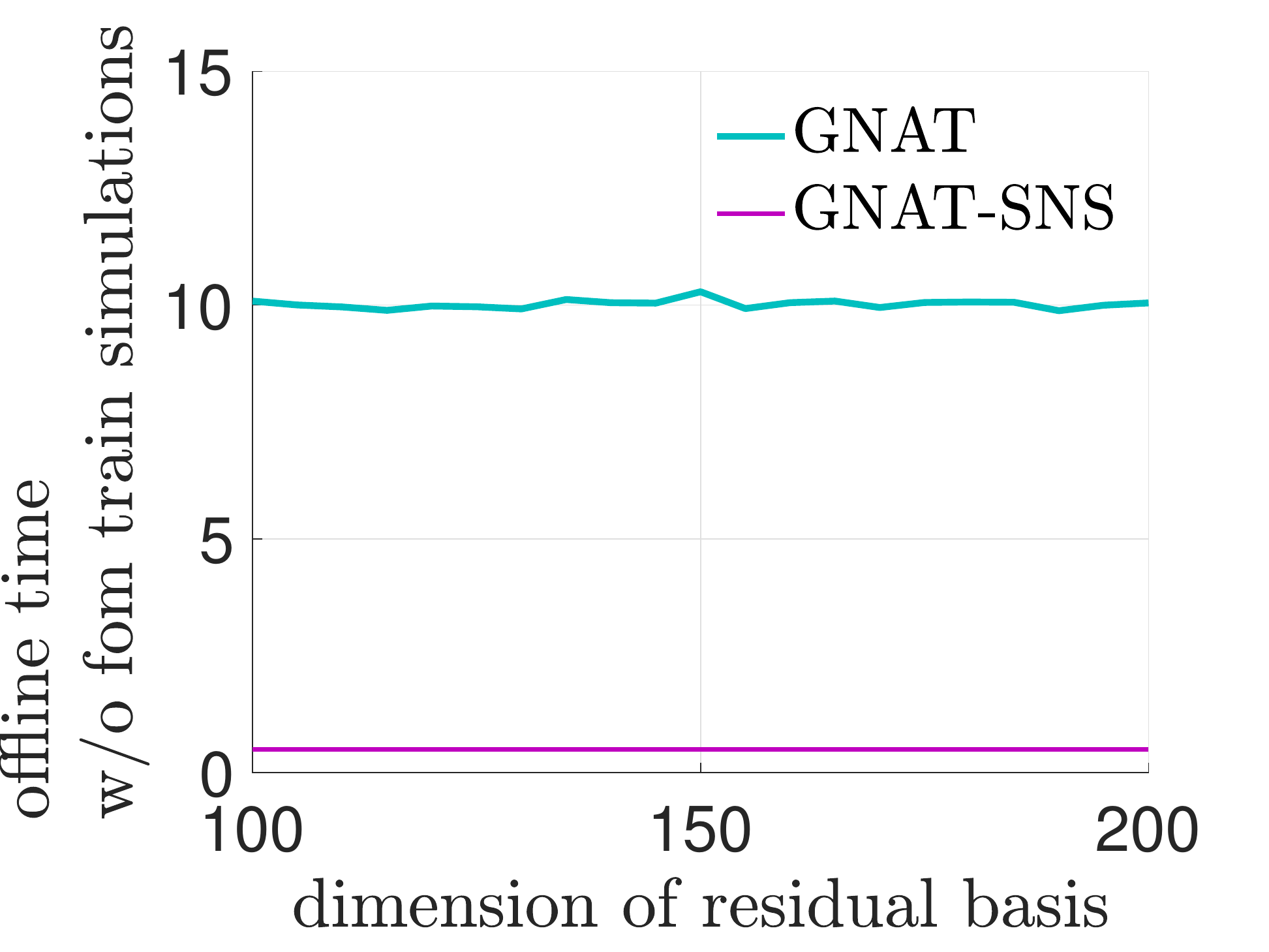}
  }~~~~~~~
  \subfigure[Adams--Bashforth]{
    \includegraphics[width=0.3\textwidth]{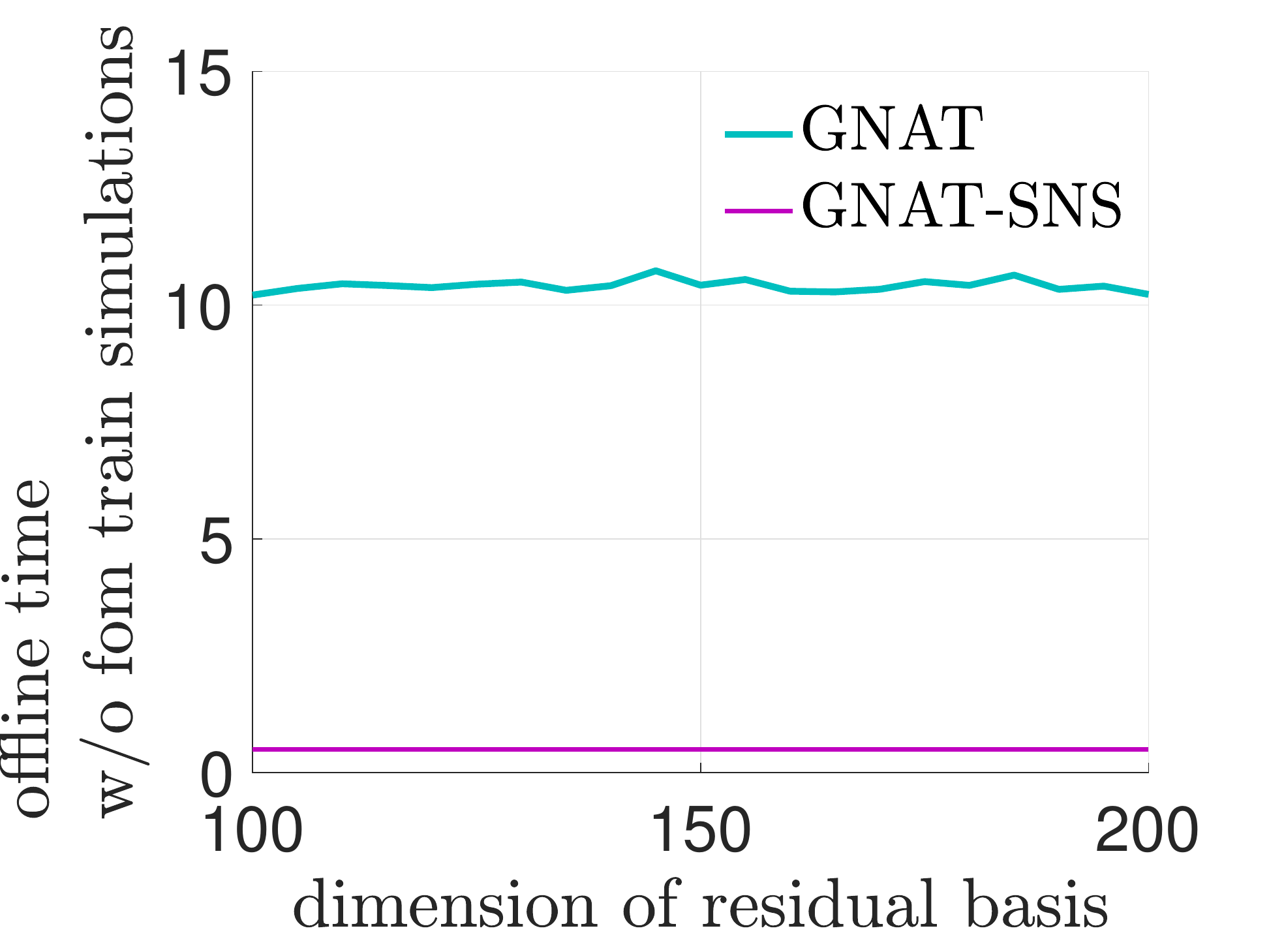}
  }~~~~~~~
  \subfigure[midpoint Runge--Kutta]{
    \includegraphics[width=0.3\textwidth]{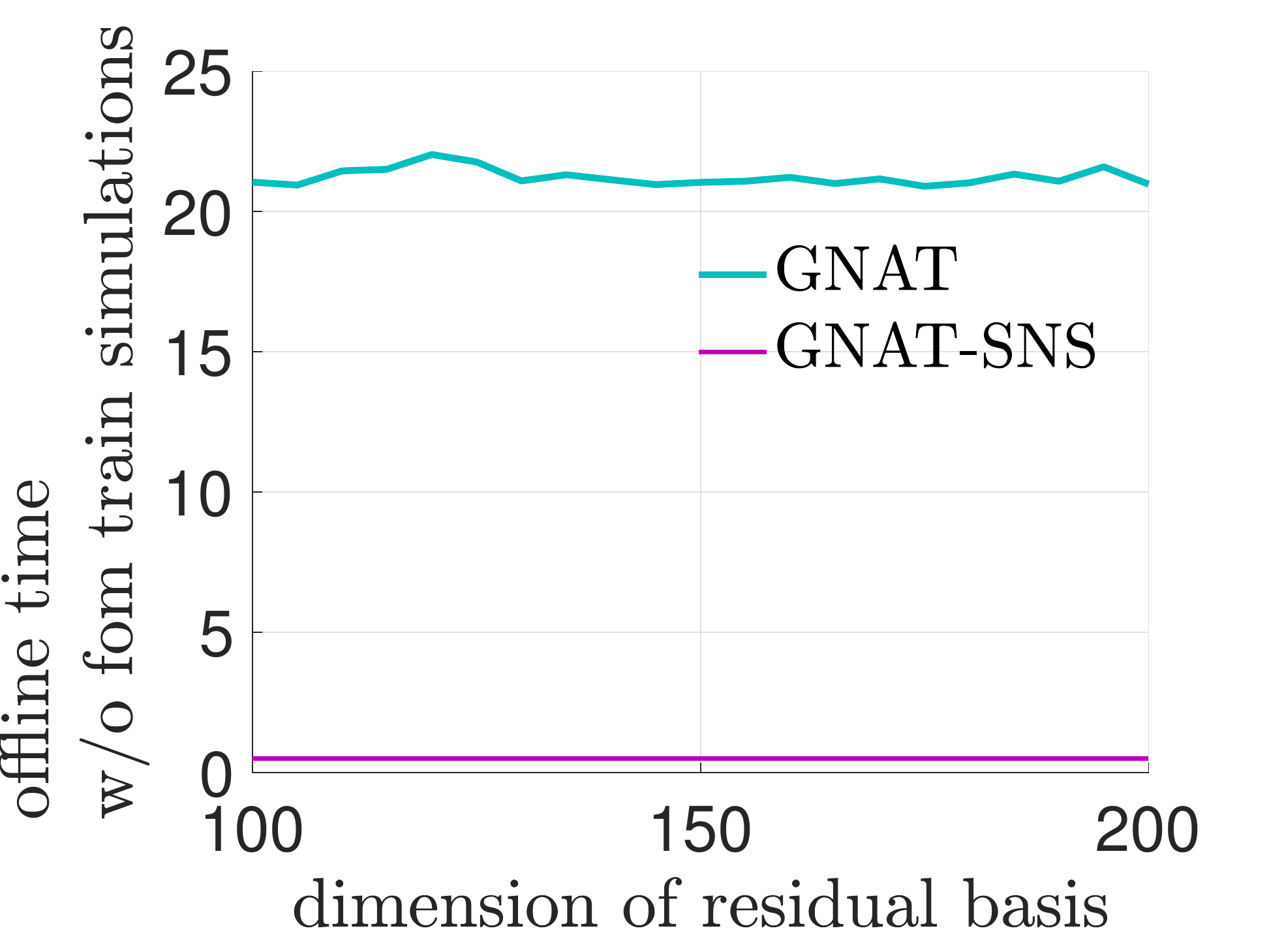}
  }
      \caption{Offline time for explicit time integrators}
    \label{fig:explicit_GNAT_burgers_offlinecost}
  \end{figure}
  Figs.~\ref{fig:explicit_GNAT_burgers_offlinecost} show the offline time of
  the GNAT and GNAT-SNS methods for the three different explicit time integrators.
  The offline time of the GNAT method includes the time of two SVDs for the
  solution and residual bases construction, the time of the training LSPG simulations for collecting
  residual snapshots, and the time of constructing sample indices. The offline
  time of the GNAT-SNS method includes the time of `one' SVD for the solution and
  residual bases and the time of constructing sample indices. Because the
  GNAT-SNS method does not need to solve the training LSPG simulations for collecting
  residual snapshots and it requires only one SVD, the offline time of the
  GNAT-SNS is less than the one of the GNAT method.
  Here, we obtain the offline computational time speed-up, a factor of 20 to 40
  with the SNS methods.
  
  \begin{figure}[htbp]
    \centering
    \subfigure[backward Euler]{
      \includegraphics[width=0.3\textwidth]{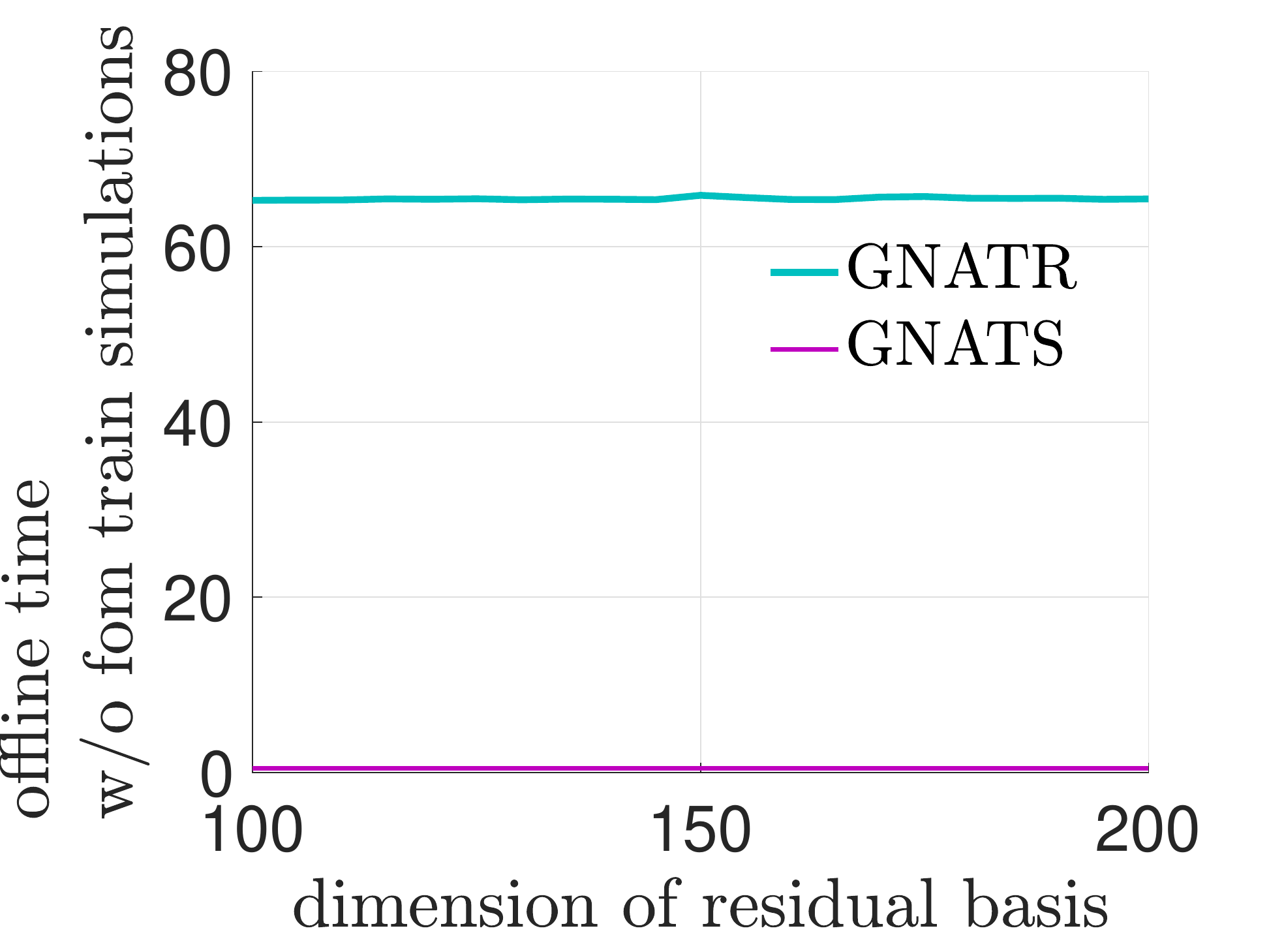}
  }~~~~~~~
  \subfigure[Adams--Moulton]{
    \includegraphics[width=0.3\textwidth]{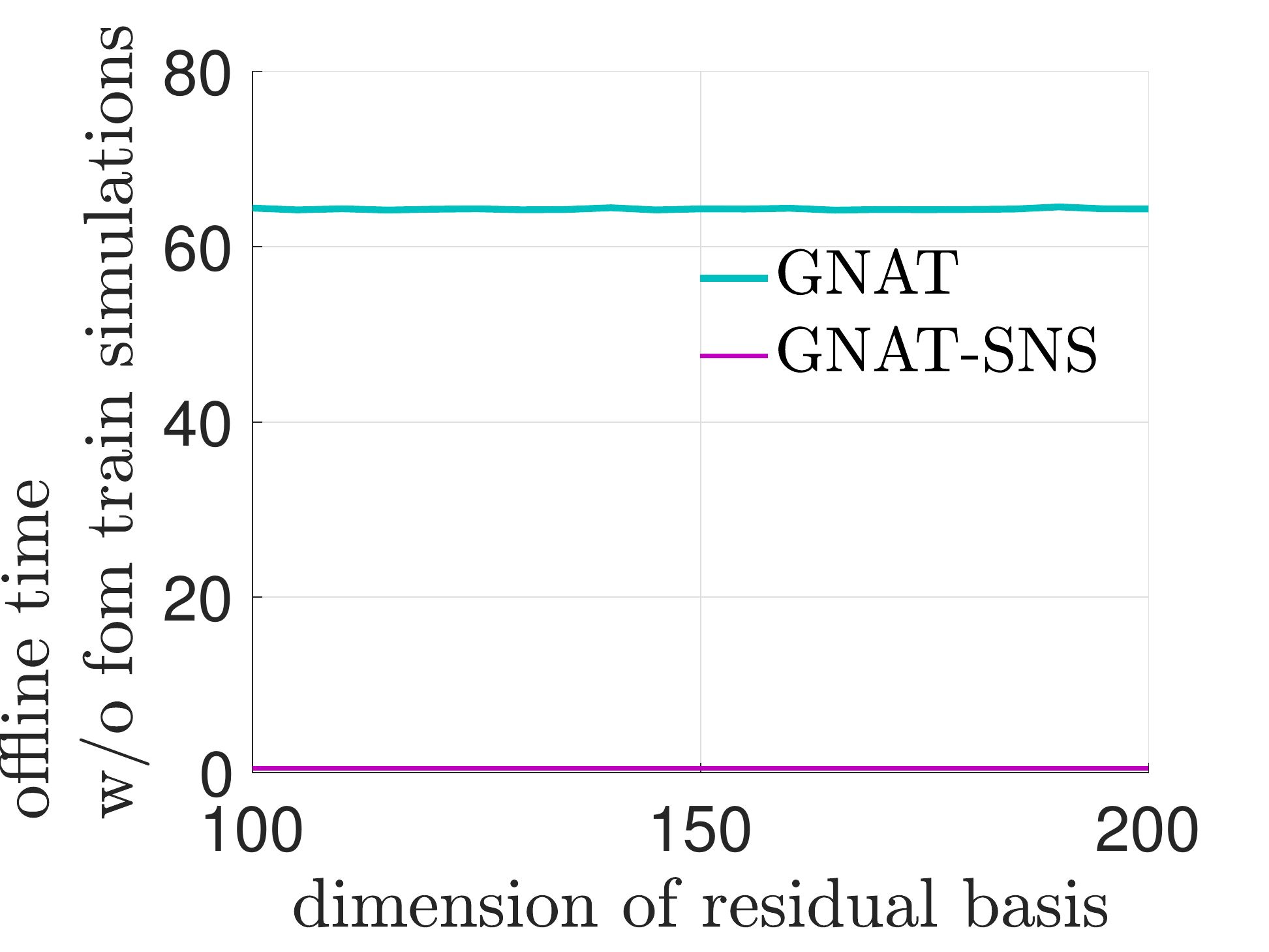}
  }~~~~~~~
  \subfigure[backward differentiation formula]{
    \includegraphics[width=0.3\textwidth]{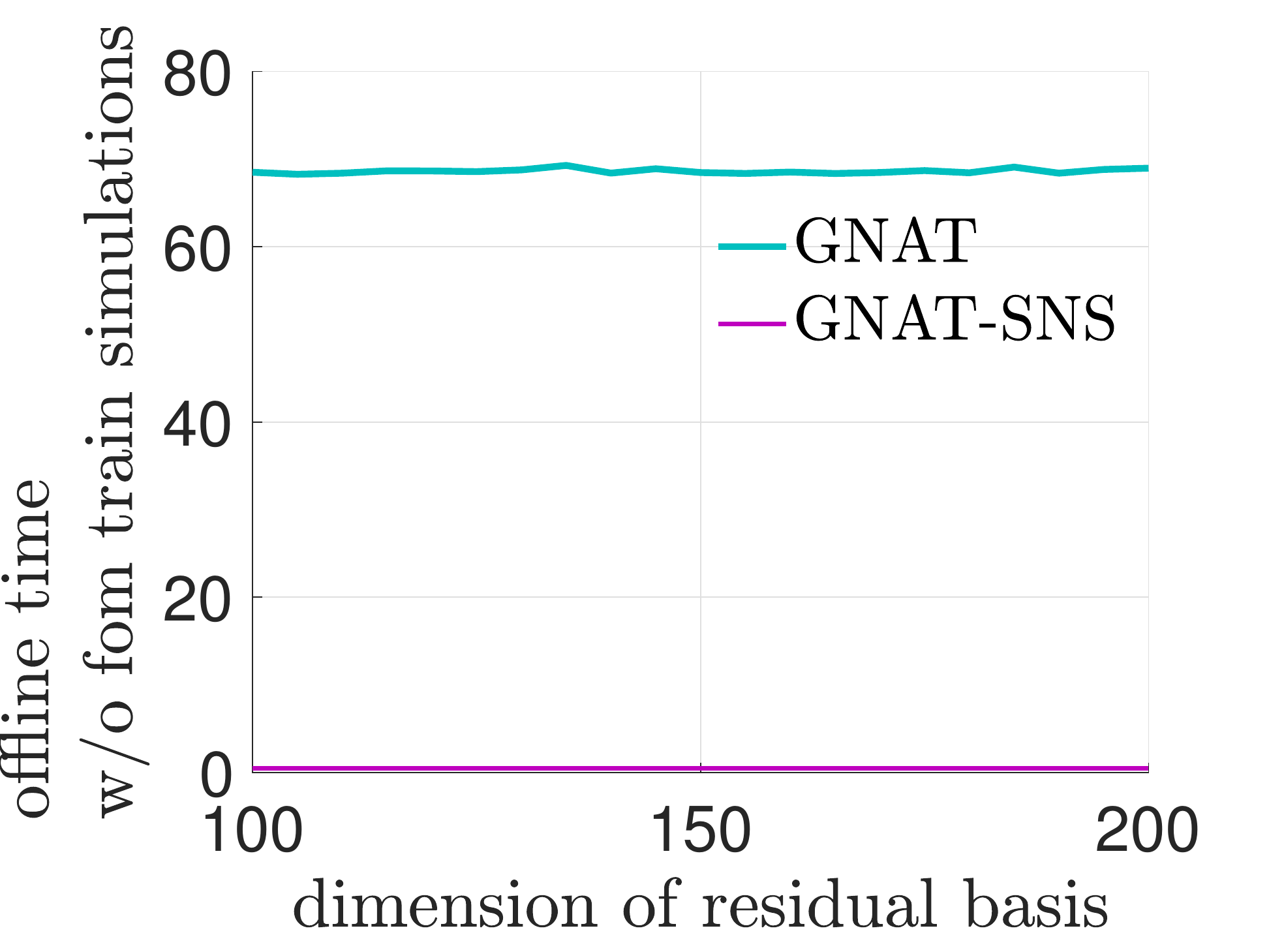}
  }
      \caption{Offline time for implicit time integrators}
    \label{fig:implicit_GNAT_burgers_offlinecost}
  \end{figure}
  Figs.~\ref{fig:implicit_GNAT_burgers_offlinecost} show the offline time of
  the GNAT and GNAT-SNS methods for the three different implicit time integrators.
  Because the GNAT-SNS method does not need to solve the training LSPG simulations for
  collecting residual snapshots and it requires only one SVD, the offline time
  of the GNAT-SNS method is less than the one of the GNAT method.
  Here, we obtain the offline computational time speed-up, an order of 100
  with the SNS methods.

  \subsubsection{ST-GNAT-SNS versus ST-GNAT}\label{sec:burgers-stgnat}
  For the space--time ROMs, the domain is discretized with $100$ control volumes,
  for $\nspacedof=100$ spatial degrees of freedom.  The time discretization used
  in the space--time ROMs is the backward Euler time integrator. We employ $\ntimedof =
  2\,000$, leading to a uniform time step of $\dt = 2.5\times10^{-4}$ s. The
  description for the various ways of collecting the space--time residual
  snapshots are shown in Sections~\ref{sec:STGNAT} and \ref{sec:ST-GNAT-SNS}. We
  use {\it ST-LSPG ROM training iterations} to collect the ST-GNAT residual basis
  snapshots. The relative error and the offline time are plotted as the dimension
  of nonlinear residual term basis $\nbasisres$ increases.  For the ST-GNAT-SNS
  method, we set $\basismatrest = \stbasismat$ if $\nbasisst = \nbasisres$, while
  $\basismatrest = \stbasismatext$ is used if $\nbasisst < \nbasisres$.

  \begin{figure}[htbp]
    \centering
    \subfigure[ST-HOSVD]{
    \includegraphics[width=0.3\textwidth]{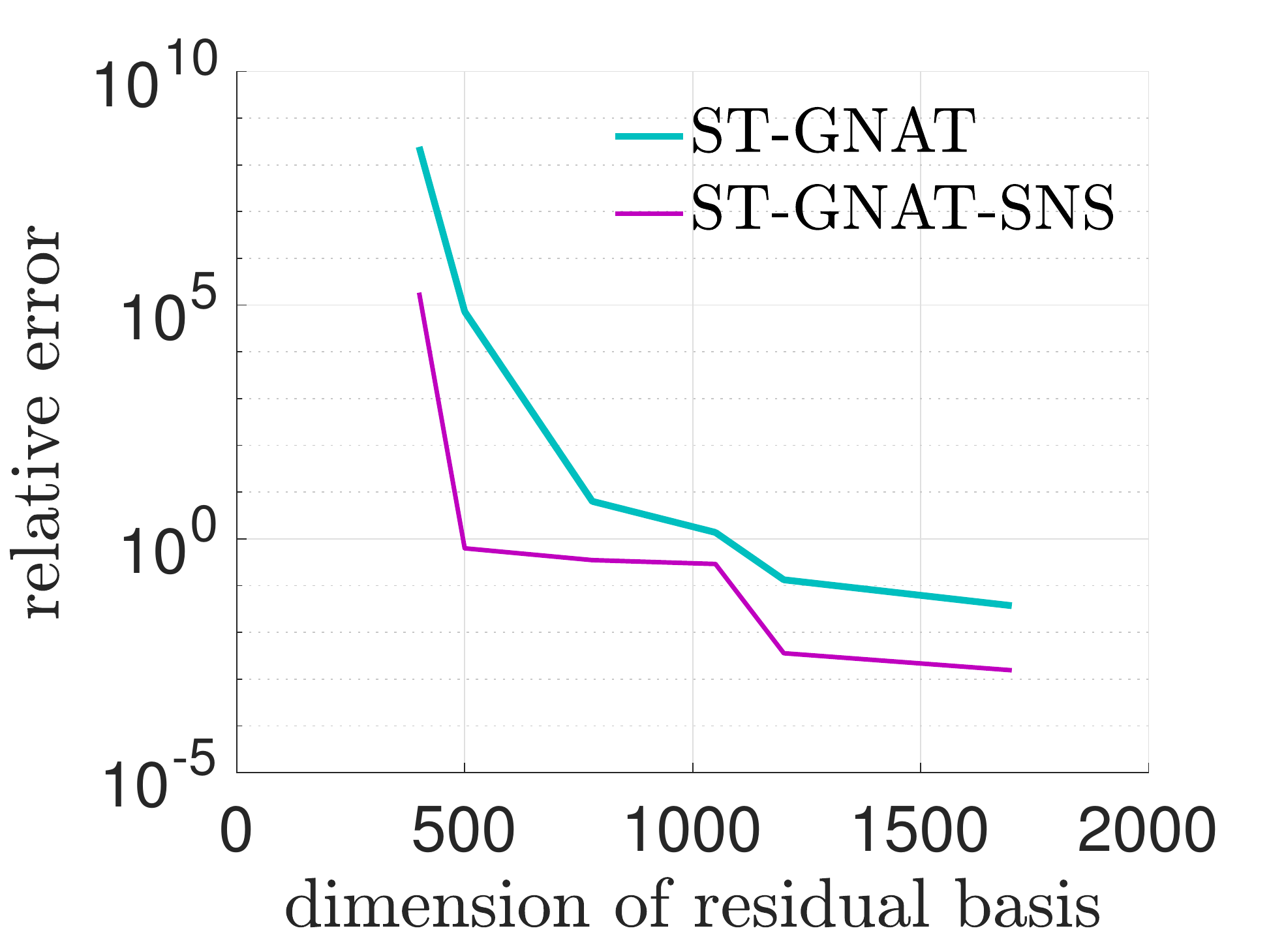}
  }~~~~~~~
  \subfigure[LL1]{
    \includegraphics[width=0.3\textwidth]{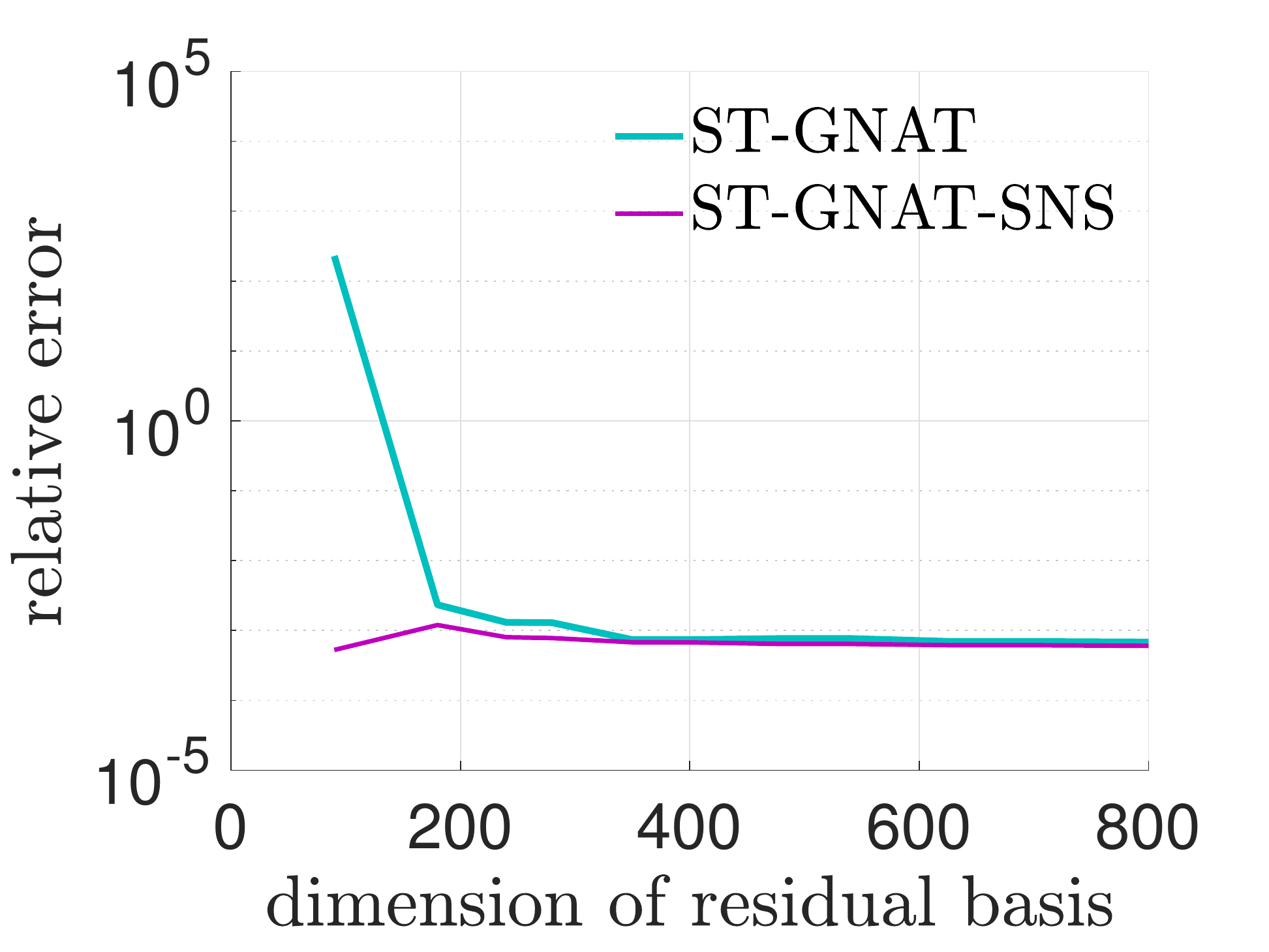}
  }~~~~~~~
      \caption{Relative errors of the space--time ROMs for solving Burgers' equation}
    \label{fig:implicit_STGNAT_burgers}
  \end{figure}
  Figs.~\ref{fig:implicit_STGNAT_burgers} show the relative error for the two
  different space--time basis generation methods, namely two different tensor
  decompositions: ST-HOSVD and LL1.  For each case, the dimension of residual
  basis $\nbasisres$ varies.  For the ST-HOSVD, we use 
  $\nbasisst = 400$ and 
  $\nbasisres\in\{400$, $500$, $780$,
  $1\,050$, $1\,200$, $1\,700\}$.  For the LL1, we use 
  $\nbasisst = 90$ and 
  $\nbasisres\in\{90$,
  $180$, $240$, $280$, $350$, $400$, $480$, $540$, $630$, $700$, $800\}$.
  For all the cases, the ST-GNAT-SNS method is comparable to the ST-GNAT
  method.  

  \begin{figure}[htbp]
    \centering
    \subfigure[ST-HOSVD]{
    \includegraphics[width=0.3\textwidth]{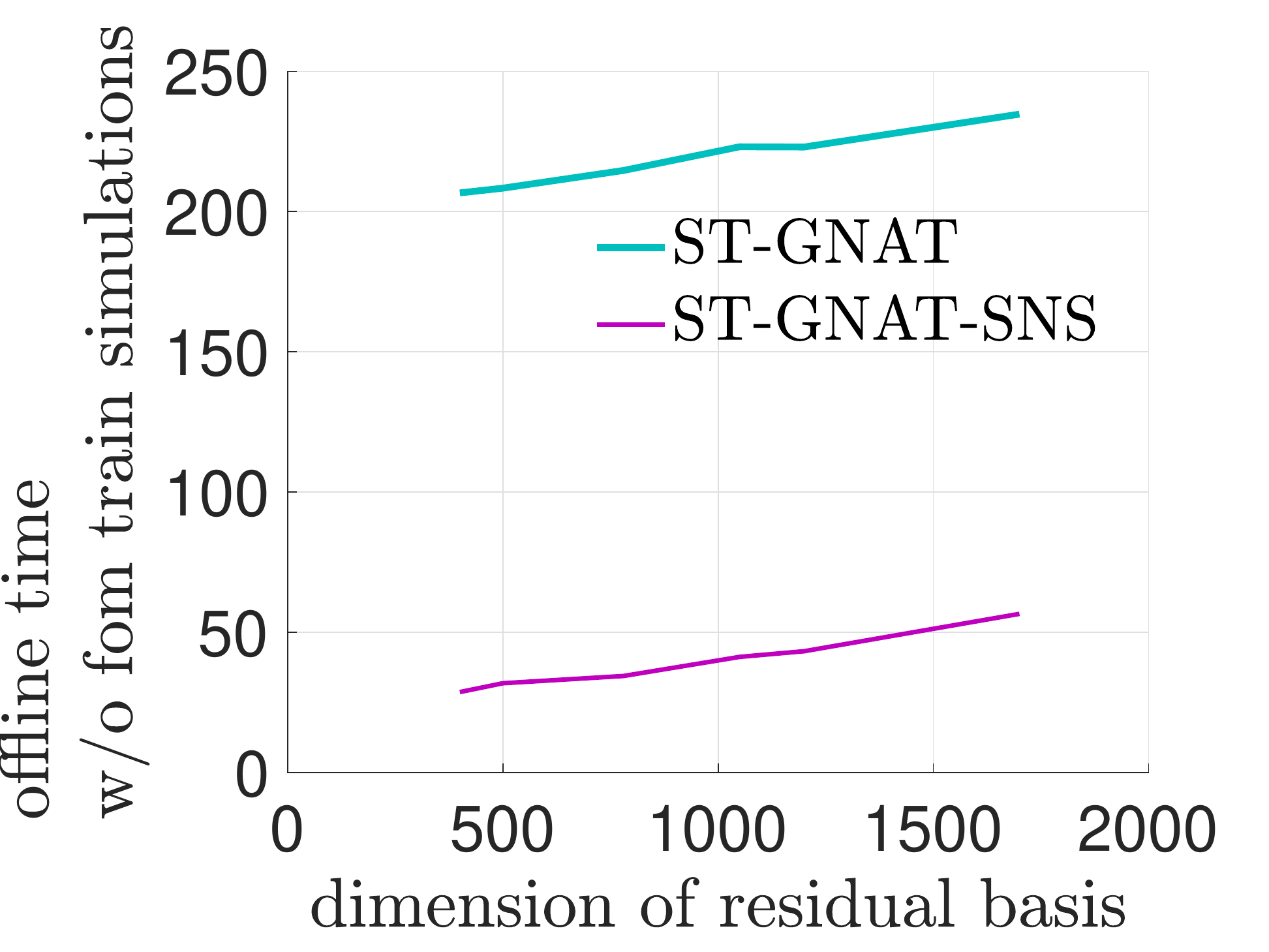}
  }~~~~~~~
  \subfigure[LL1]{
    \includegraphics[width=0.3\textwidth]{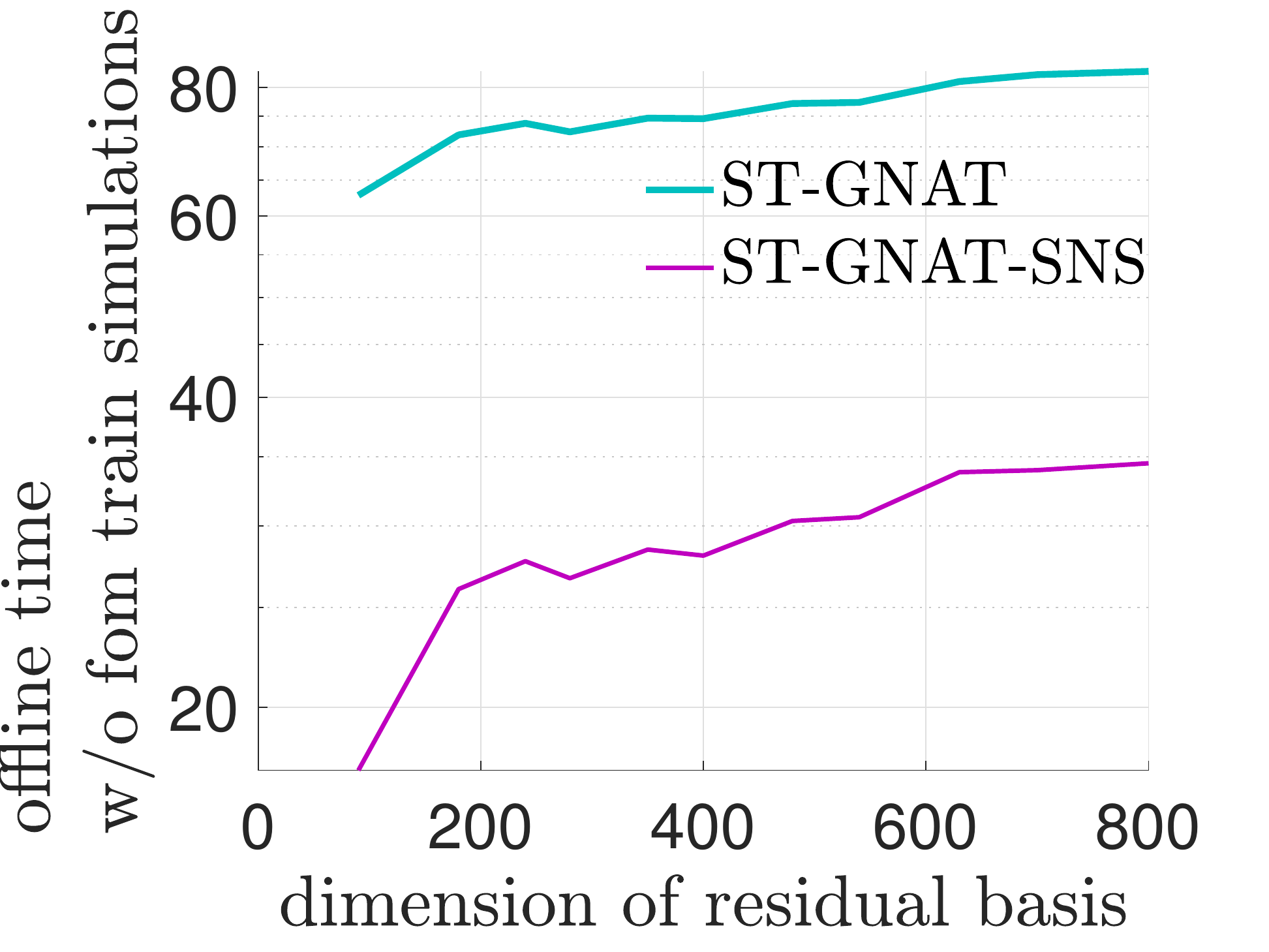}
  }~~~~~~~
      \caption{Offline time of the space--time ROMs for solving Burgers' equation}
    \label{fig:implicit_STGNAT_burgers_offlinetime}
  \end{figure}
  Figs.~\ref{fig:implicit_STGNAT_burgers_offlinetime} show the offline time of
  the ST-GNAT and ST-GNAT-SNS methods. The offline time of the ST-GNAT
  method includes the time of two tensor decompositions (e.g., ST-HOSVD or LL1)
  for the solution and residual bases construction, the time of the training ST-LSPG
  simulations for collecting residual snapshots, and the time of constructing
  sample indices. The offline time of the ST-GNAT-SNS method includes the time
  of `one' tensor decomposition for the solution and residual bases and the
  time of constructing sample indices. Because the ST-GNAT-SNS does not need to
  solve the training ST-LSPG simulations for collecting residual snapshots and
  it only requires one tensor decomposition, the offline time of the ST-GNAT-SNS
  method is less
  than the one of the ST-GNAT method.
  Here, we obtain the offline computational time speed-up, a factor of three to
  seven with the SNS methods.

\subsection{Parameterized quasi-1D Euler equation}\label{sec:eulerProblem}
\begin{figure}[htbp] 
  \centering 
  \includegraphics[width=0.4\textwidth]{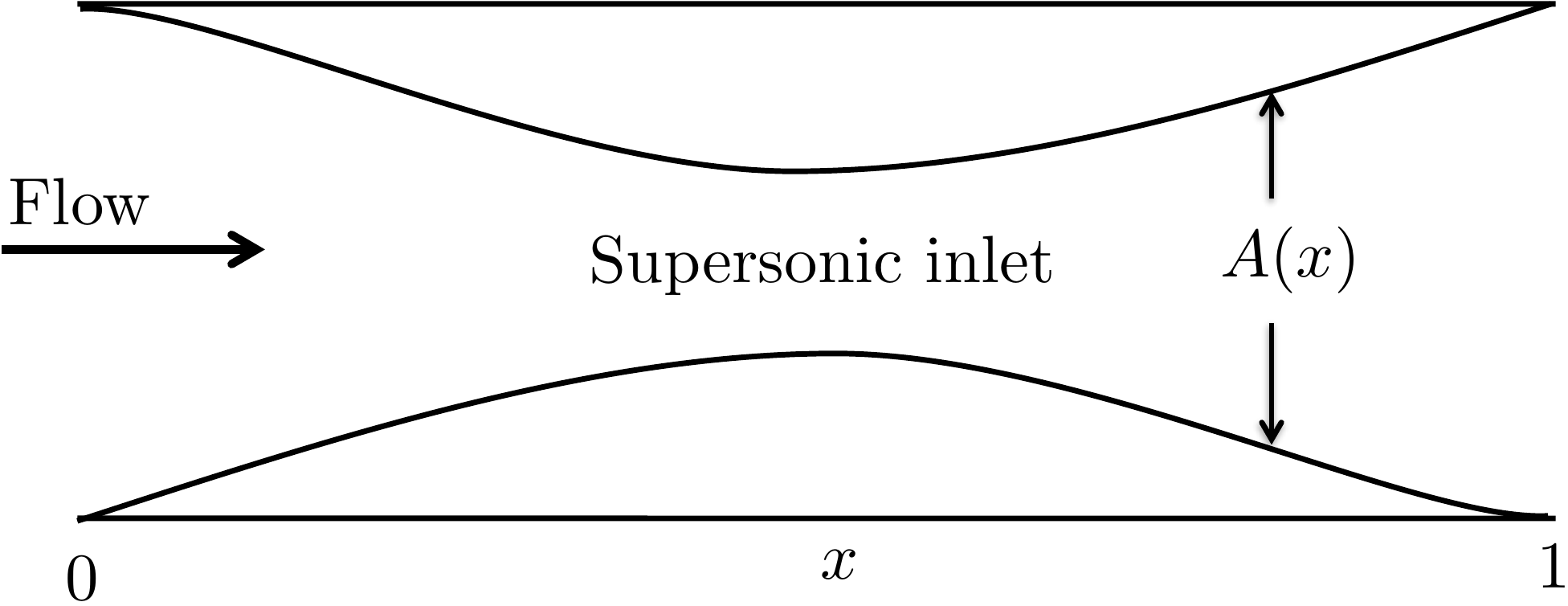} 
  \caption{\textit{Quasi-1D Euler.} Schematic figures of converging-diverging nozzle. } 
  \label{fig:converging-diverging_nozzle} 
\end{figure}

We now consider a parameterized quasi-1D Euler equation associated with
modeling inviscid compressible flow in a one-dimensional converging--diverging nozzle with a
continuously varying cross-sectional area \cite[Chapter 13]{maccormackNote};
Figure \ref{fig:converging-diverging_nozzle} depicts the problem geometry. 

The governing system of nonlinear partial differential equations is
 \begin{equation}\label{eq:euler} 
   \frac{\partial \solw}{\partial t} + \frac{1}{\crossareaSymb}\frac{\partial (\flux(\solw) \crossareaSymb)}{\partial x} 
   = \externalforce(\solw), \quad \forall x\in[0,1]\ \text{m}, \quad \forall t\in[0,\totaltime], 
 \end{equation} 
where $\totaltime=0.6$ s and
 \begin{equation} 
      \solw = \bmat{\densitySymb \\ \densitySymb \velocitySymb \\ \energydensity}, 
\quad \flux(\solw) = \bmat{\densitySymb \velocitySymb \\ \densitySymb \velocitySymb^2 + \pressureSymb \\ 
                    (\energydensity+\pressureSymb)\velocitySymb},
\quad \externalforce(\solw) = \bmat{0 \\ \frac{\pressureSymb}{\crossareaSymb}\frac{\partial \crossareaSymb}{\partial x}\\0},
\quad
      \pressureSymb = (\specificheat-1)\densitySymb\energypermass,
\quad \energypermass = \frac{\energydensity}{\densitySymb} - \frac{\velocitySymb^2}{2},
\quad \crossareaSymb = \crossareaSymb(x).
 \end{equation}
Here,
$\densitySymb$ denotes density,
$\velocitySymb$ denotes velocity,
$\pressureSymb$ denotes pressure, 
$\energypermass$ denotes potential energy per unit mass,
$\energydensity$ denotes total energy density,
$\specificheat$ denotes the specific heat ratio, and
$\crossareaSymb$ denotes the converging--diverging nozzle cross-sectional area.
We employ a
specific heat ratio of $\specificheat=1.3$, 
a specific gas constant of $\specificgasconstant=355.4$
$\text{m}^2/\text{s}^2/\text{K}$,
a total temperature of $\temperatureSymb_t = 300$ K,
and a total pressure of $\pressureSymb_t = 10^6$ $\text{N}/\text{m}^2$.
The cross-sectional area $A(x)$ is determined by a cubic spline interpolation over the points 
$(x,A(x)) \in \{ (0,0.2)$, $(0.25,0.173)$, $(0.5,0.17)$, $(0.75, 0.173)$, $(1,0.2)
\}$, which results in
 \begin{equation}
   A(x) = \left\{  
             \begin{array}{ll}
                 -0.288x^3 + 0.4080x^2 - 0.1920x + 0.2,                   &
								 x\in[0,0.25)\ \text{m} \\
           -0.288(x-0.25)^3 + 0.1920(x-0.25)^2 - 0.0420(x-0.25) + 0.1730, & x\in[0.25,0.5) \ \text{m}\\
            0.288(x-0.5)^3 - 0.0240(x-0.5)^2 + 0.17,                      & x\in[0.5,0.75) \ \text{m}\\
            0.288(x-0.75)^3 + 0.1920(x-0.75)^2 + 0.0420(x-0.75) + 0.1730, & x\in[0.75,1]\ \text{m}.
             \end{array}
          \right.   
 \end{equation} 
A perfect gas is assumed that obeys the ideal gas law (i.e., $\pressureSymb = \densitySymb \specificgasconstant \temperatureSymb$).
The initial flow field is created in several steps.
First, the following isentropic relations are used to generate a zero pressure-gradient flow field:
 \begin{align}
   \begin{split}
     \machSymb(x) &= \frac{\machSymb_m\crossareaSymb_m}{\crossareaSymb(x)}\left ( \frac{1+\frac{\specificheat-1}{2}\machSymb(x)^2}{1+\frac{\gamma-1}{2}\machSymb_m^2} \right
  )^{\frac{\specificheat+1}{2(\specificheat-1)}},\quad
  \pressureSymb(x) = \pressureSymb_t\left (1+\frac{\specificheat-1}{2}\machSymb(x)^2\right )^{\frac{-\specificheat}{\specificheat-1}}
  \\ \temperatureSymb(x) &= \temperatureSymb_t\left
	(1+\frac{\specificheat-1}{2}\machSymb(x)^2\right )^{-1},\quad
  \densitySymb(x) = \frac{\pressureSymb(x)}{R\temperatureSymb(x)},\quad
  \speedofsoundSymb(x) =
	\sqrt{\specificheat\frac{\pressureSymb(x)}{\densitySymb(x)}} ,\quad
  \velocitySymb(x) = \machSymb(x) \speedofsoundSymb(x),
   \end{split}
 \end{align} 
 where a subscript $m$ indicates the flow quantity at $x=0.5$ m, and
 $\machSymb$ denotes the
 Mach number.
Then, a shock is placed at $x=0.85$ m of the flow field. 
   The jump relations for a stationary shock and 
   the perfect gas equation of state are used to derive the velocity across the shock
	 $\velocitySymb_2$, which satisfies
   the quadratic equation
 \begin{equation} \label{eq:quadratic}
   \left ( \frac{1}{2}-\frac{\specificheat}{\specificheat-1} \right ) \velocitySymb_2^2 
   + \frac{\specificheat}{\specificheat-1}\frac{n}{m}\velocitySymb_2 - h = 0.
 \end{equation}  
Here, $m \defeq \densitySymb_2\velocitySymb_2 = \densitySymb_1\velocitySymb_1$, 
      $n\defeq\densitySymb_2\velocitySymb_2^2+\pressureSymb_2 = \densitySymb_1\velocitySymb_1^2+\pressureSymb_1$,
      $h\defeq(\energydensity_2+\pressureSymb_2)/\densitySymb_2 =
			(\energydensity_1+\pressureSymb_1)/\densitySymb_1$, 
		and subscripts $1$ and $2$ denote a flow quantity to the left and to the
		right of
the shock, respectively. The solution $\velocitySymb_2$ is employed
			to Eq.~\eqref{eq:quadratic} that leads to a discontinuity (i.e., shock).
Finally, the exit pressure is increased to a factor $\pexit$ of its
original value in order to generate transient dynamics.  

Applying a finite-volume spatial discretization with 50 equally spaced control
volumes and fully implicit boundary conditions leads to a parameterized system
of nonlinear ODEs consistent with Eq.~\eqref{eq:fom} with $\nspacedof = 150$
spatial degrees of freedom.  The Roe flux difference vector splitting method
is used to compute the flux at each intercell face \cite[Chapter
9]{maccormackNote}.  For time discretization, we apply the backward
Euler scheme and a uniform time step of $\dt = 0.001$ s, leading to $\ntimedof=600$.

For this problem, we use the following two parameters: 
the pressure factor $\paramSymb_1 = \pexit$ and the Mach number at the middle of the nozzle
$\paramSymb_2 = \machSymb_m$.
All ROMs employ a training set at which the FOM is solved of
$\paramDomainTrain = \{1.7+0.01i\}_{i=0}^3\times \{1.7,1.72\}$ such that
$\ntrain = 8$. Then the target parameter, $\param = (1.7225,1.705)$, is pursued.

\subsubsection{GNAT-SNS versus GNAT}\label{sec:euler-gnat}
  The solution basis dimension of $\nbasisspace = 30$ is used. 
  The relative error and the offline time are plotted as
  the dimension of nonlinear residual term basis $\nbasisres$ increases.
  For the GNAT-SNS method, $\basismatres = \mass\basismatspace$ is used if
  $\nbasisres = \nbasisspace$,
  while $\basismatres = \mass\basismatspaceext$ is used if $\nbasisres > \nbasisspace$.

  \begin{figure}[htbp]
    \centering
    \subfigure[$\nbasisres = 30$]{
    \includegraphics[width=0.3\textwidth]{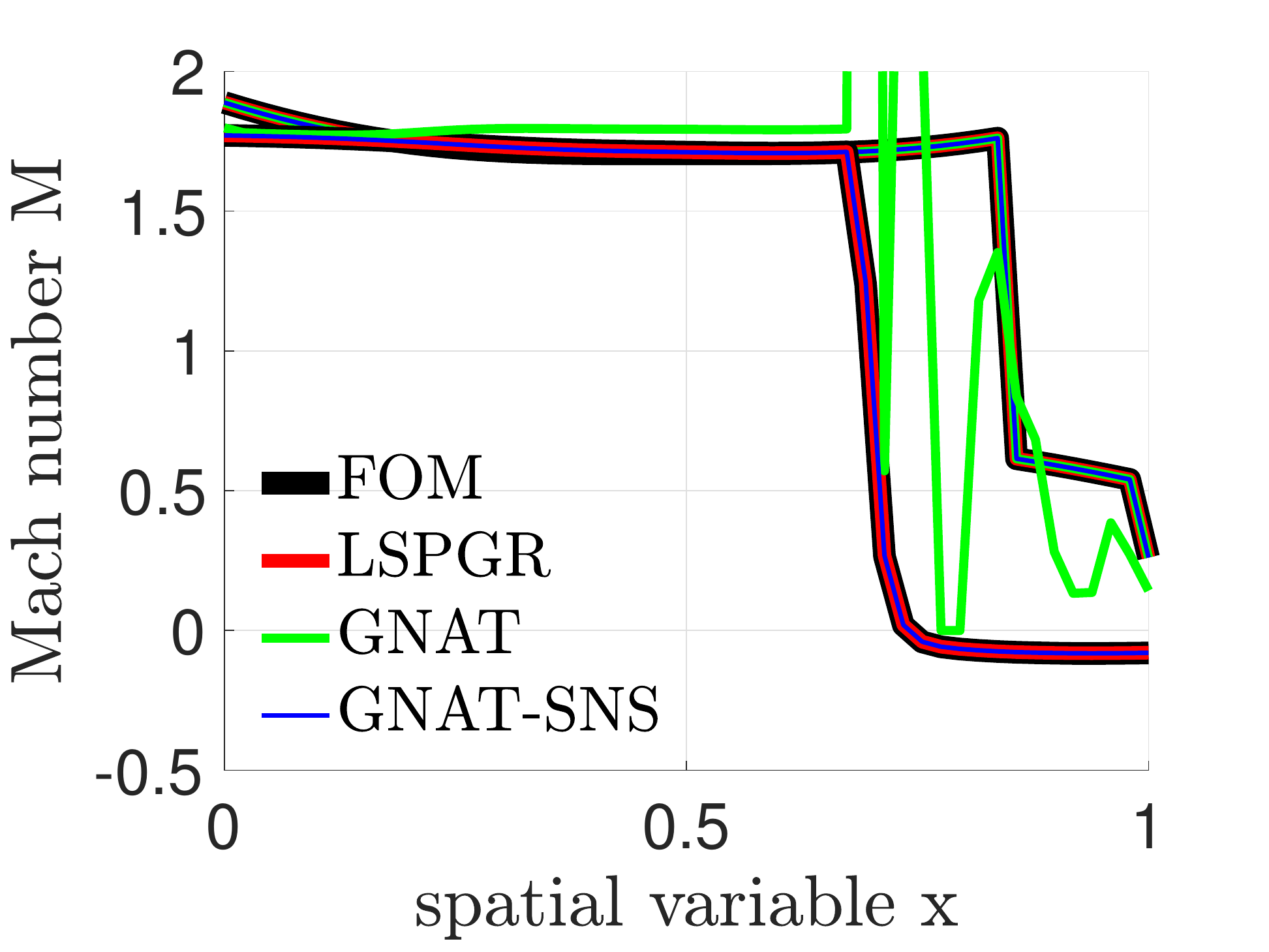}
  }~~~~~~~
  \subfigure[$\nbasisres = 60$]{
    \includegraphics[width=0.3\textwidth]{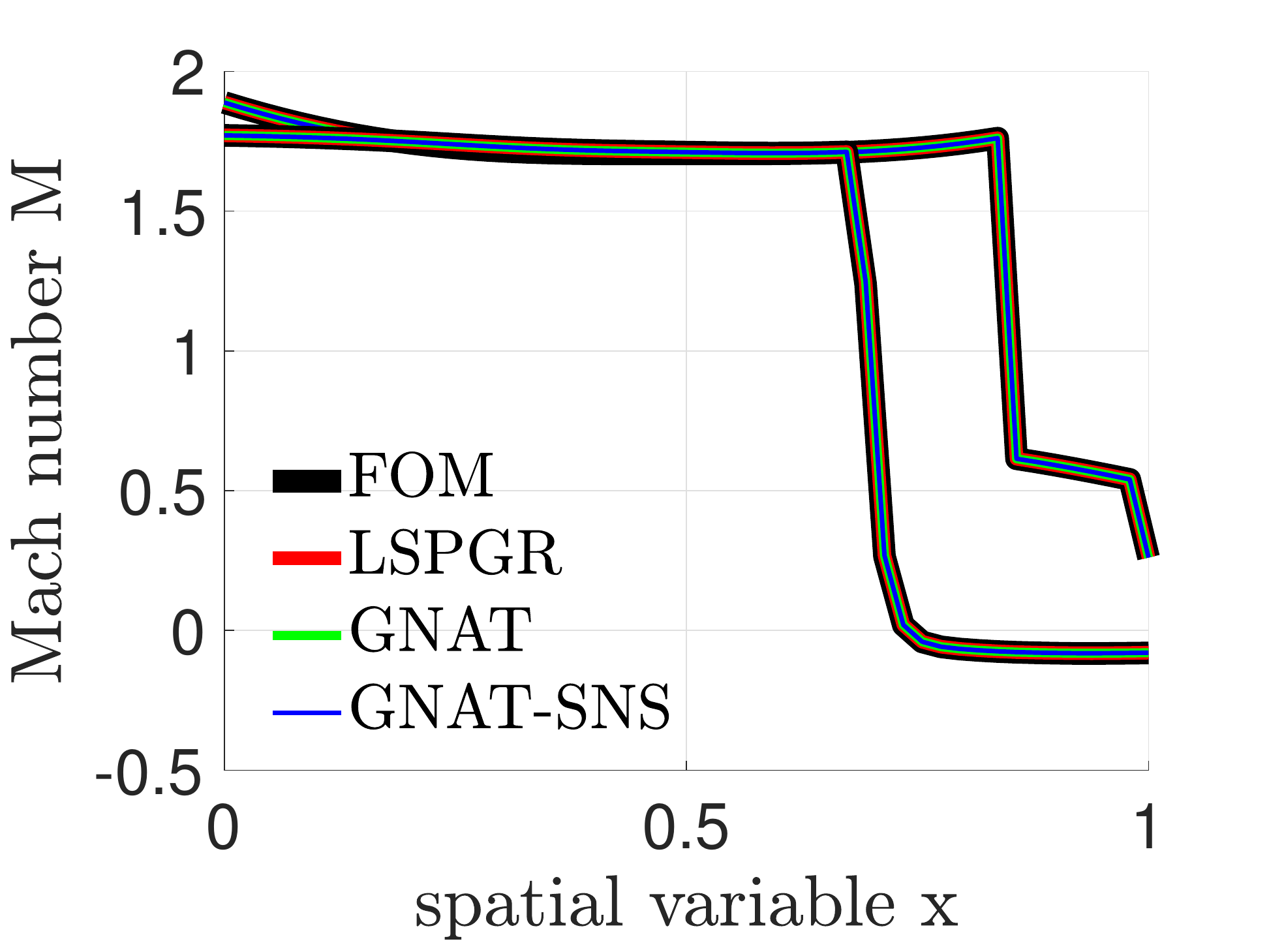}
  }
  \caption{Solution snapshots at $t\in\{0,0.6\}$, $\nressample = 90$ with the
  backward Euler time integrator.}
    \label{fig:GNAT_euler_snapshots}
  \end{figure}
  Figs.~\ref{fig:GNAT_euler_snapshots} compare the solution snapshots of several
  methods: the LSPG, GNAT, and GNAT-SNS methods with the FOM solution snapshots.
  Fig.~\ref{fig:GNAT_euler_snapshots}(a) is generated with $\nbasisres = 30$ and
  $\nressample = 90$, while
  Fig.~\ref{fig:GNAT_euler_snapshots}(b) is generated with $\nbasisres = 60$ and
  $\nressample = 90$.
  All the methods are able to generate almost the same solutions as the FOM
  solutions except for the GNAT method with $\nbasisspace = \nbasisres = 30$.   
  Surprisingly, the GNAT-SNS method does not suffer when $\nbasisspace =
  \nbasisres$ as shown in Figs.~\ref{fig:implicit_GNAT_burgers} of the Burgers'
  example. Again, we do not know why the GNAT-SNS the achieve an accuracy as good as the
  the one of the LSPG method in this particular example.

  \begin{figure}[htbp]
    \centering
    \subfigure[relative error]{
    \includegraphics[width=0.3\textwidth]{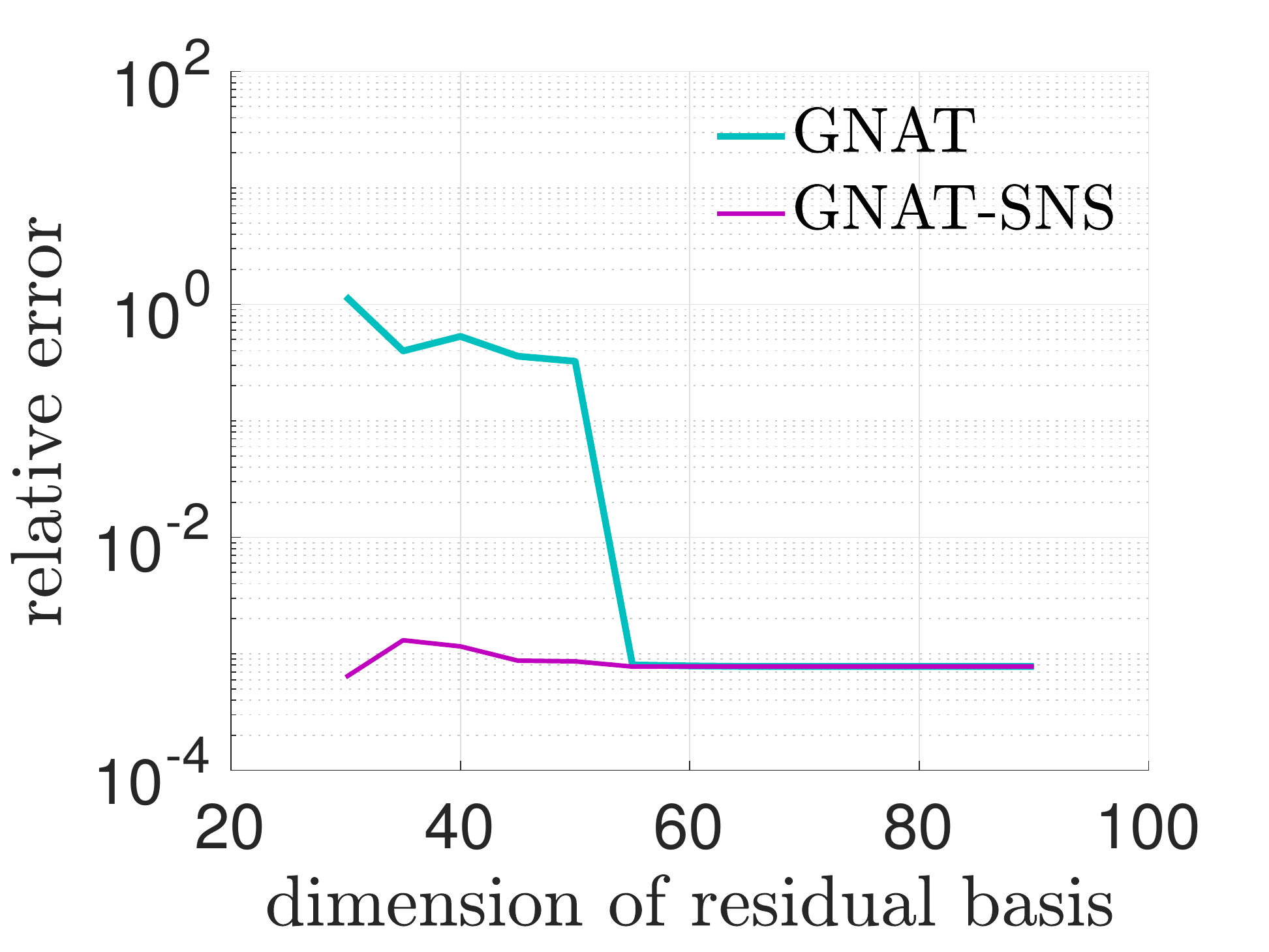}
  }~~~~~~~
  \subfigure[offline time]{
    \includegraphics[width=0.3\textwidth]{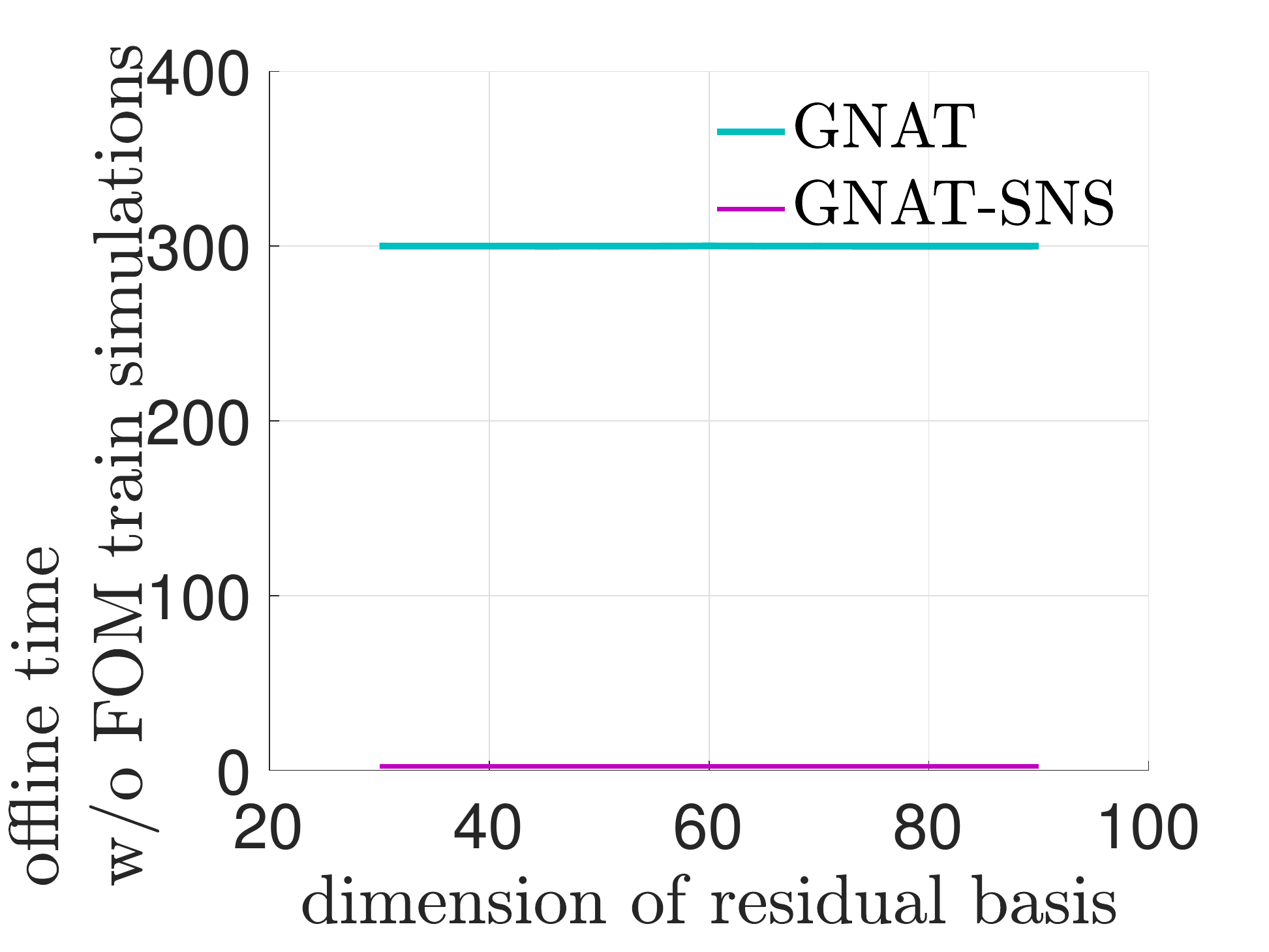}
  }
    \caption{Relative errors and offline time with the backward Euler time integrator}
    \label{fig:GNAT_euler}
  \end{figure}
  In Fig.~\ref{fig:GNAT_euler}(a), the relative errors of both the GNAT and
  GNAT-SNS methods are plotted as the dimension of residual basis increases from
  30 to 90 by 5 with a fixed number of sample indices, 90.  The figures show
  that the GNAT-SNS method is comparable to the GNAT method in terms of
  accuracy; the order of relative errors of the GNAT-SNS method is $10^{-3}$ for
  the whole range of the residual basis dimensions considered here.  On the other
  hand, the relative errors of the GNAT method are bigger than the ones of the GNAT-SNS
  method when the residual basis dimensions are between 30 and 55.
  Fig.~\ref{fig:GNAT_euler}(b) shows the offline time of the GNAT and
  GNAT-SNS methods. The offline time of the GNAT method includes the time of two SVDs
  for the solution and residual bases construction, the time of the training LSPG
  simulations for collecting residual snapshots, and the time of constructing
  sample indices. The offline time of the GNAT-SNS method includes the time of `one'
  SVD for the solution and residual bases construction and the time of constructing sample
  indices. Because the GNAT-SNS method does not need to solve the training LSPG
  simulations for collecting residual snapshots and it requires only one SVD,
  the offline time of the GNAT-SNS method is less than the one of the GNAT
  method. Here, we get a factor of 130 speed-up in the offline computational
  time with the GNAT-SNS method.

\subsubsection{ST-GNAT-SNS versus ST-GNAT}\label{sec:euler-stgnat}
  The relative error and the offline time of both the ST-GNAT and ST-GNAT-SNS
  methods are plotted as
  the dimension of nonlinear residual term basis $\nbasisres$ increases.
  For the ST-GNAT-SNS method, we set $\basismatrest = \stbasismat$ if $\nbasisst =
  \nbasisres$, while $\basismatrest = \stbasismatext$ is used if $\nbasisst <
  \nbasisres$.

  \begin{figure}[htbp]
    \centering
    \subfigure[ST-HOSVD]{
    \includegraphics[width=0.3\textwidth]{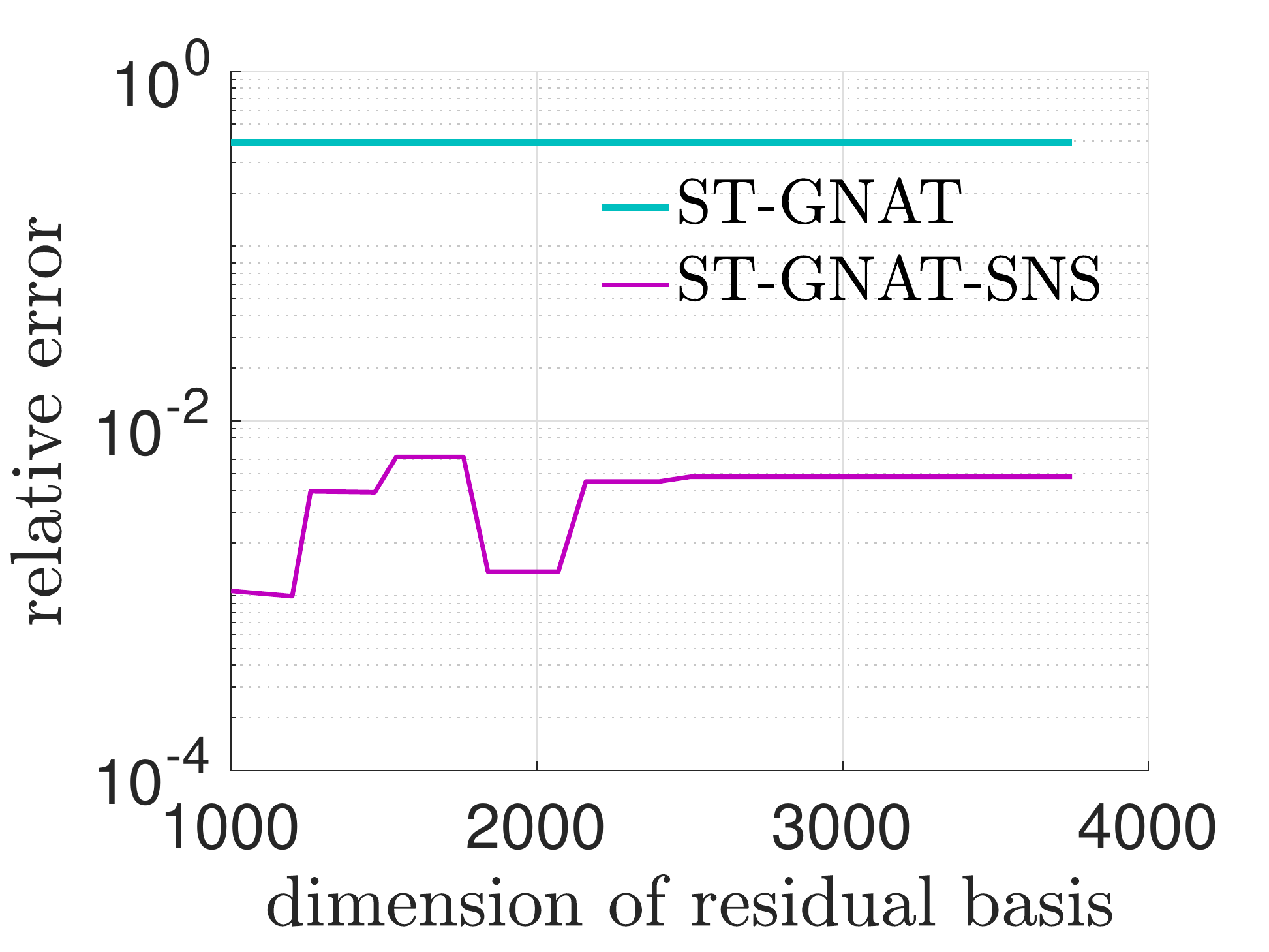}
  }~~~~~~~
    \subfigure[LL1]{
    \includegraphics[width=0.3\textwidth]{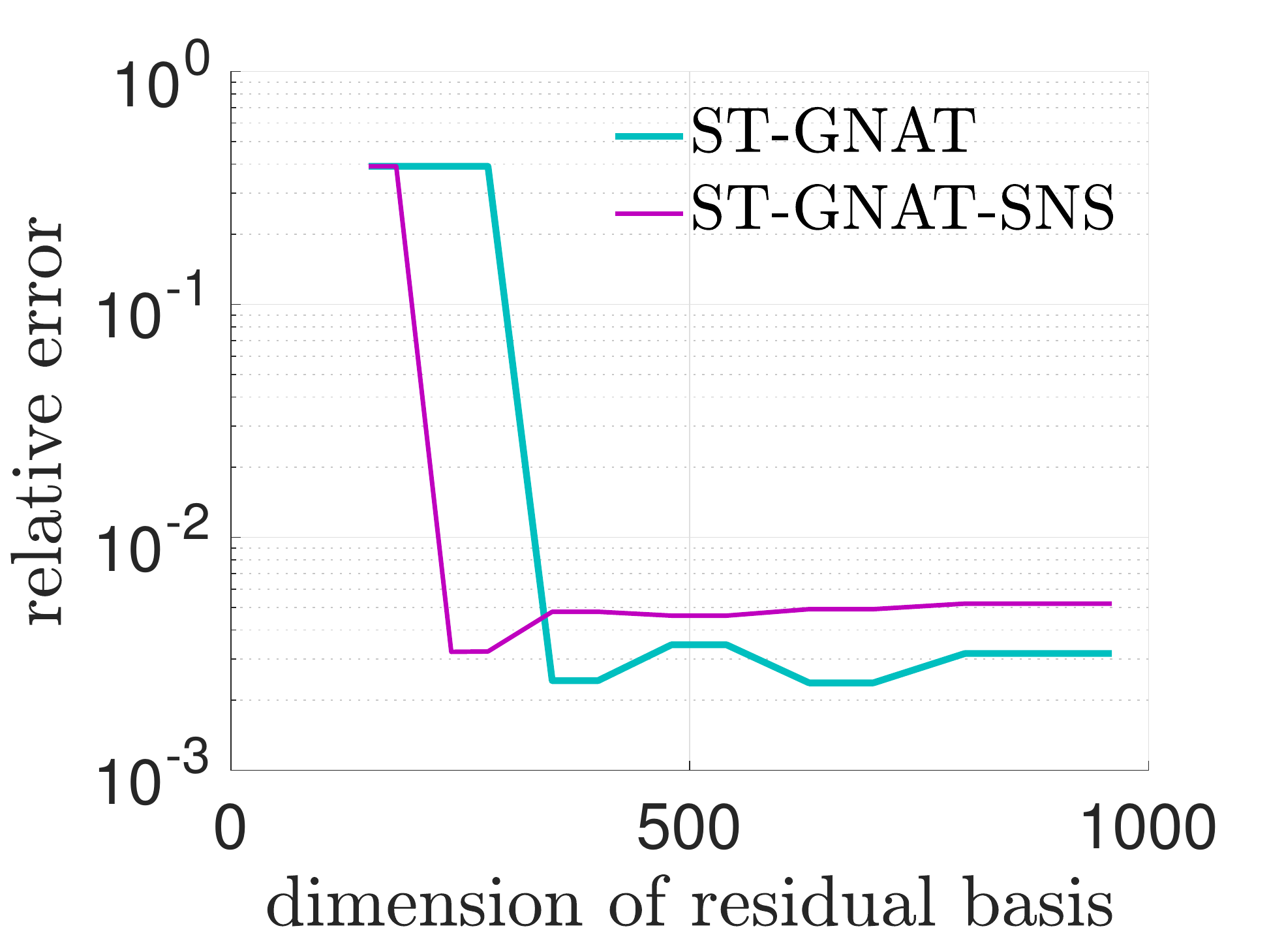}
  }~~~~~~~
    \caption{Relative errors and offline time with the backward Euler time integrator}
    \label{fig:relerror_ST-GNAT_euler}
  \end{figure}
  In Figs.~\ref{fig:relerror_ST-GNAT_euler}, the relative errors of 
  both the ST-GNAT and ST-GNAT-SNS methods are shown
  for two different tensor decompositions: ST-HOSVD and LL1.
  For the ST-HOSVD, we use $\nbasisst = 1\,000$ and $\nbasisres\in\{ 1\,000$, 
    $1\,200$, $1\,260$, $1\,470$, $1\,540$, $1\,760$, $1\,840$, $2\,070$,
  $2\,160$, $2\,400$, $2\,500$, $3\,000$, $3\,500$, $3\,750\}$. 
  For the LL1 decomposition, we use $\nbasisst = 150$ and $\nbasisres\in\{150$, $180$, $240$,
  $280$, $350$, $400$, $480$, $540$, $630$, $700$, $800$, $960\}$.
  For the ST-HOSVD, the ST-GNAT-SNS method achieves two orders of magnitude
  better accuracy than the ST-GNAT method. For the LL1 decomposition, the
  ST-GNAT-SNS method generate results as good as the ST-GNAT method.

  \begin{figure}[htbp]
    \centering
    \subfigure[ST-HOSVD]{
    \includegraphics[width=0.3\textwidth]{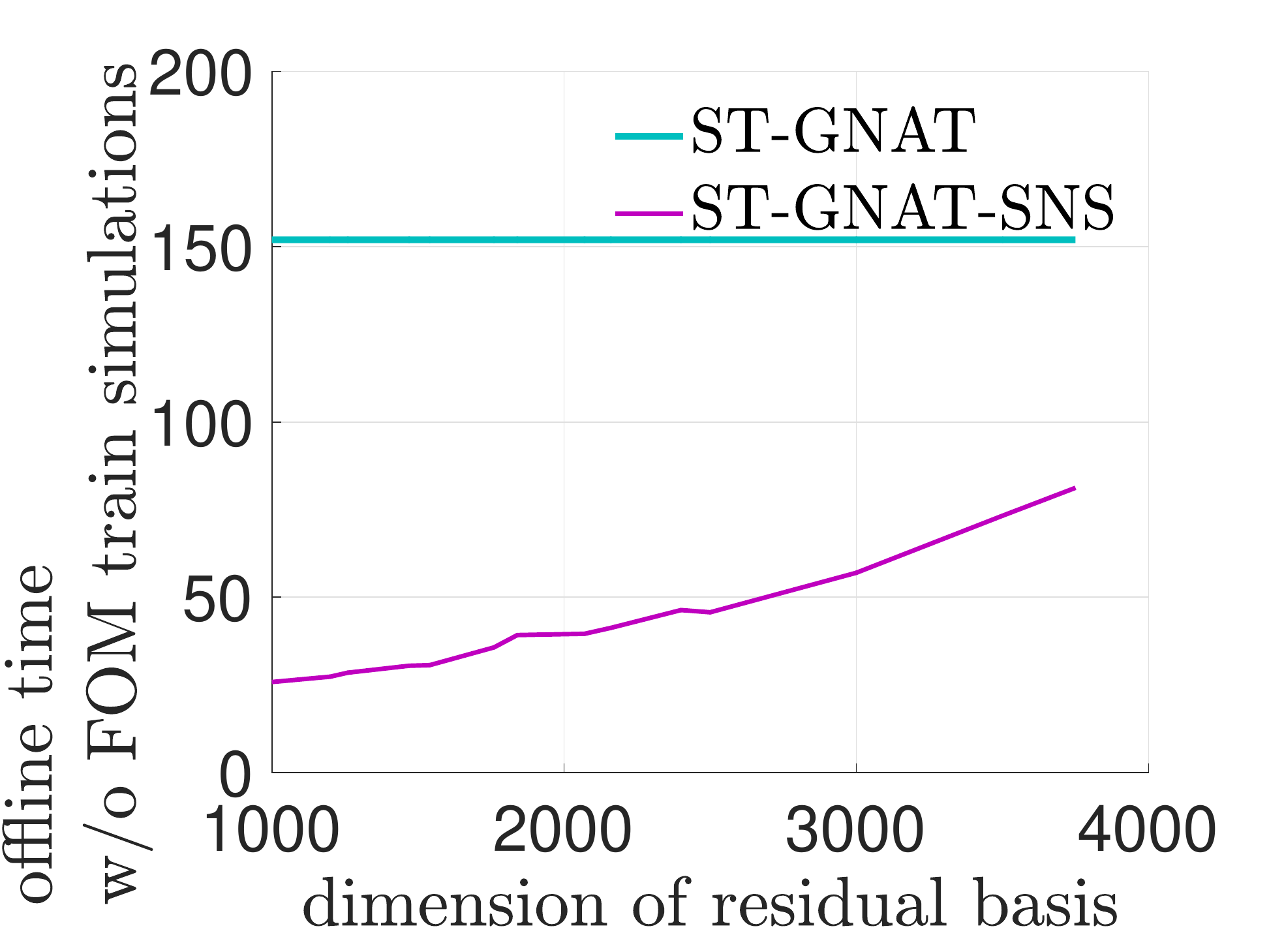}
  }~~~~~~~
  \subfigure[LL1]{
    \includegraphics[width=0.3\textwidth]{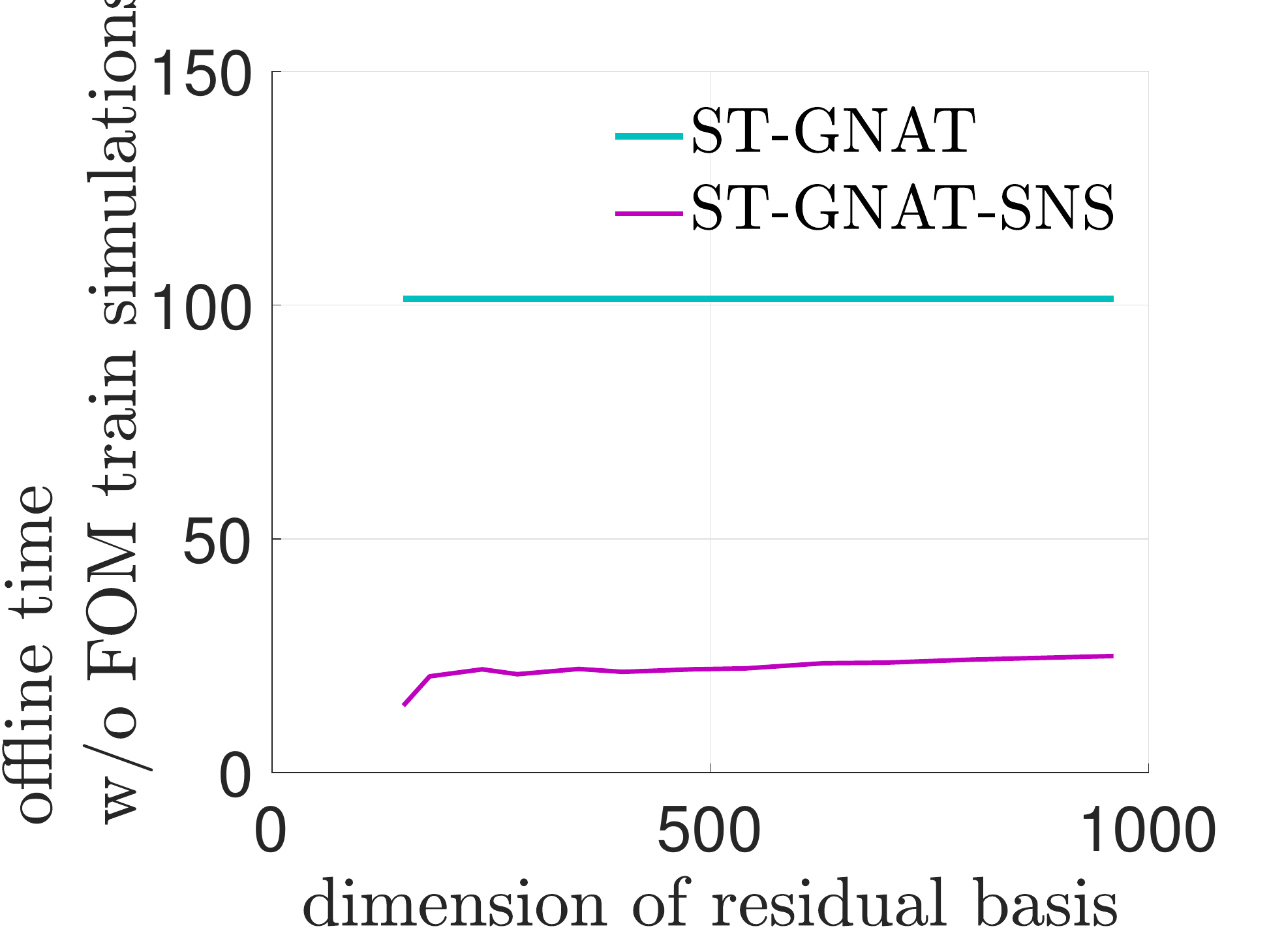}
  }~~~~~~~
      \caption{Offline time of the space--time ROMs for solving Euler equation}
    \label{fig:STGNAT_euler_offlinetime}
  \end{figure}
  Figs.~\ref{fig:STGNAT_euler_offlinetime} show the offline time of the ST-GNAT
  and ST-GNAT-SNS methods.
  The offline time of the ST-GNAT method includes the time of two tensor decompositions 
  (e.g., ST-HOSVD or LL1) for the
  solution and residual bases construction, 
  the time of the training ST-LSPG simulations for collecting residual snapshots, and
  the time of constructing sample indices. The offline time of the ST-GNAT-SNS
  method includes the time of `one' tensor
  decomposition for the solution and residual bases construction 
  and the time of constructing sample indices. Because the
  ST-GNAT-SNS does not need to solve the training ST-LSPG simulations 
  for collecting residual snapshots and it only
  requires one tensor decomposition, 
  the offline time of the ST-GNAT-SNS method is less than the
  one of the ST-GNAT method. Here, we get a factor of six to seven speed-up in
  the offline computational time with the SNS methods. 

\section{Conclusion}\label{sec:conclusion} 
We have introduced a new way of constructing nonlinear term basis using solution
snapshots to construct a projection-based reduced order model for solving a
nonlinear dynamical system of equations, which is characterized by an ordinary
differential equation.  Our proposed method, the SNS method, is theoretically
justified by the conforming subspace condition and the subspace inclusion
relation. 

The two main advantages of the SNS method over the traditional hyper-reduction
methods considered here are 1) the avoidance of collecting nonlinear term
snapshots and 2) the avoidance of the second data compression process. These
advantages result in the offline computational time reduction, which is
demonstrated in numerical experiments. The benefits of the SNS method are more vivid when
it is compared with GNAT and ST-GNAT than DEIM because the GNAT and ST-GNAT methods
require collecting the nonlinear residual snapshots from the corrsponding LSPG
and ST-LSPG simulations. These benefits are also demonstrated in numerical
experiments, where parametric GNAT and ST-GNAT are compared with the SNS
methods. There, we have shown that a considerable speed-up in the offline
computational time is achieved by the SNS methods. 
Also, the error analysis contributes to the theoretical insight about the effects of
the volume matrix $\mass$ on the oblique projection of the SNS method, revealing
that the conditioner number of $\mass$ can affect the error due to the oblique
projection. This issue can be addressed by pre-orthogonalization process. 

In numerical experiments, we observe that the SNS methods
produce a more accurate solution than the DEIM, GNAT, and ST-GNAT methods with a
smaller number of nonlinear term basis vectors, especially for the GNAT-SNS
method. The reason for this attractive feature of the GNAT-SNS method will be
investigated in future work. 

\section*{Acknowledgments}
The authors acknowledge the helpful comments provided by the anonymous reviewers.
This work was performed at Lawrence Livermore National Laboratory and was
supported by the LDRD program (project 43320) and the Engineering Summer Intern
program (project 45415).
Lawrence Livermore National Laboratory is operated by Lawrence
Livermore National Security, LLC, for the U.S. Department of Energy,
National Nuclear Security Administration under Contract DE-AC52-07NA27344
and LLNL-JRNL-767241.

\appendix
\section{Subspace inclusion relation}\label{sec:appendix}
The forward and backward Euler time integrators and their corresponding subspace
inclusion relations are shown in Section~\ref{sec:FOM}.
Here, we show several other time integrators and corresponding subspace
inclusion relation between the subspace spanned by the nonlinear term snapshots
and the subspace spanned by the solution snapshots.

\subsection{The Adams--Bashforth methods}\label{sec:AB}
      The second order Adams--Bashforth (AB) method numerically solves 
      Eq.~\eqref{eq:fom}, by solving the following nonlinear system of
      equations for $\solArg{n}$ at $n$-th time step:
      \begin{equation} \label{eq:AF2}
        \mass \solArg{n} - \mass \solArg{n-1} = \dt\left (
        \frac{3}{2}\fluxArg{n-1}-\frac{1}{2}\fluxArg{n-2}\right ),
      \end{equation} 
      Eq.~\eqref{eq:AF2} implies the following subspace inclusions:
      \begin{equation}\label{eq:AF2_partial_inclusion}
        \Span{\fluxArg{n-2},\fluxArg{n-1}} \subseteq 
        \Span{\mass \solArg{n-1},\mass \solArg{n}}.
      \end{equation}
      By induction, we conclude that
      \begin{equation}\label{eq:AF2_total_inclusion}
        \Span{\fluxArg{0},\dots,\fluxArg{\ntimedof-1}} \subseteq \Span{\mass \solArg{0},\ldots,\mass \solArg{\ntimedof}},
      \end{equation}
      which shows that the span of nonlinear term snapshots is 
      a subspace of the span of $\mass$-scaled solution snapshots. 
      The residual function of the second AB method is defined as
      \begin{align}\label{eq:residual_AF2} 
      \begin{split}
        \resn_{\ADB}(\solArg{n};\solArg{n-1},\param) &\defeq 
        \mass(\solArg{n} - \solArg{n-1}) 
        -\dt\left ( \frac{3}{2}\fluxArg{n-1}-\frac{1}{2}\fluxArg{n-2}\right ).
      \end{split}
      \end{align} 
\subsection{The Adams--Moulton methods}\label{sec:AM}
      The second order Adams--Moulton (AM) method numerically solves 
      Eq.~\eqref{eq:fom}, by solving the following nonlinear system of
      equations for $\solArg{n}$ at $n$-th time step:
      \begin{equation} \label{eq:AM2}
        \mass \solArg{n} - \mass \solArg{n-1} = \frac{1}{2}\dt(
        \fluxArg{n}+\fluxArg{n-1}),
      \end{equation} 
      Eq.~\eqref{eq:AM2} implies the following subspace inclusions:
      \begin{equation}\label{eq:AM2_partial_inclusion}
        \Span{\fluxArg{n-1},\fluxArg{n}} \subseteq \Span{\mass \solArg{n-1},\mass \solArg{n}}.
      \end{equation}
      By induction, we conclude that
      \begin{equation}\label{eq:AM2_total_inclusion}
        \Span{\fluxArg{0},\dots,\fluxArg{\ntimedof}} \subseteq \Span{\mass \solArg{0},\ldots,\mass \solArg{\ntimedof}},
      \end{equation}
      which shows that 
      the span of nonlinear term snapshots is 
      a subspace of the span of $\mass$-scaled solution snapshots. 
      The residual function of the second AM method is defined as
      \begin{align}\label{eq:residual_AM2} 
      \begin{split}
        \resn_{\ADM}(\solArg{n};\solArg{n-1},\param) &\defeq 
        \mass(\solArg{n} - \solArg{n-1}) 
        -\dt\frac{1}{2}(\fluxArg{n}+\fluxArg{n-1}).
      \end{split}
      \end{align} 
\subsection{The backward differentiation formulas}\label{sec:BDF}
      The 2-step BDF numerically solves Eq.~\eqref{eq:fom}, by solving the following nonlinear system of equations for $\solArg{n}$ at n-th time step:
      \begin{equation} \label{eq:BDF2}
        \mass \solArg{n} - \frac{4}{3}\mass \solArg{n-1} + \frac{1}{3}\mass \solArg{n-2}
        = \frac{2}{3}\dt\fluxArg{n},
      \end{equation}
      Eq.~\eqref{eq:BDF2} implies the following subspace inclusions:
      \begin{equation}\label{eq:BDF2_partial_inclusion}
        \Span{\fluxArg{n}} \subseteq
        \Span{\mass \solArg{n-2}, \mass \solArg{n-1}, \mass \solArg{n}}.
      \end{equation}
      By induction, we conclude that
      \begin{equation}\label{eq:BDF2_total_inclusion}
        \Span{\fluxArg{1},\dots,\fluxArg{\ntimedof-1}} \subseteq
        \Span{\mass \solArg{0},\ldots,\mass \solArg{\ntimedof}},
      \end{equation}
      which shows that the span of nonlinear term snapshots is
      a subspace of the span of $\mass$-scaled solution snapshots.
      The residual function of the two-step BDF method is defined as
      \begin{align}\label{eq:residual_BDF2}
      \begin{split}
        \resn_{\BDF}(\solArg{n};\solArg{n-1},\solArg{n-2},\param) &\defeq
        \mass(\solArg{n} - \frac{4}{3}\solArg{n-1} + \frac{1}{3}\solArg{n-2}) 
        -\frac{2}{3}\dt\fluxArg{n}.
      \end{split}
      \end{align}
\subsection{The midpoint Runge--Kutta method}\label{sec:RK}
    The midpoint method, a $2$-stage Runge--Kutta method, 
    takes the following two stages to advance at n-th time step of 
    Eq.~\eqref{eq:fom}:
    \begin{equation}\label{eq:RK2}
      \begin{aligned}
        \mass \solArg{n-\frac{1}{2}} &= \mass \solArg{n-1} +
          \frac{\dt}{2}\fluxArg{n-1} \\
        \mass \solArg{n} &= \mass \solArg{n-1} + \dt\fluxArg{n-\frac{1}{2}}.
      \end{aligned}
    \end{equation}
    These lead to the following two subspace inclusion relations:
    \begin{equation}\label{eq:RK2_partial_inclusion}
      \begin{aligned}
        \Span{\fluxArg{n-1}} &\subseteq
        \Span{\mass \solArg{n-1}, \mass \solArg{n-\frac{1}{2}}} \\
        \Span{\fluxArg{n-\frac{1}{2}}} &\subseteq
        \Span{\mass \solArg{n-1}, \mass \solArg{n}}. \\
      \end{aligned}
    \end{equation}
    These in turn lead to the following subspace inclusion relation:
    \begin{equation}\label{eq:RK2_partial_inclusion2}
        \Span{\fluxArg{n-1},\fluxArg{n-\frac{1}{2}}} \subseteq
      \Span{\mass \solArg{n-1}, \mass \solArg{n-\frac{1}{2}},\mass\solArg{n}}. 
    \end{equation}
    By induction, we conclude that
    \begin{equation}\label{eq:RK2_inclusion}
        \Span{\fluxArg{0},\fluxArg{\frac{1}{2}},\ldots,
        \fluxArg{\ntimedof-1},\fluxArg{\ntimedof-\frac{1}{2}}} \subseteq
      \Span{\mass \solArg{n-1}, \mass \solArg{n-\frac{1}{2}},\mass\solArg{n},\ldots,
        \mass \solArg{\ntimedof-1},\mass\solArg{\ntimedof-\frac{1}{2}},
      \mass\solArg{\ntimedof}}.
    \end{equation}

\bibliographystyle{plain}
\bibliography{references}

\begin{thebibliography}{10}

\bibitem{blast-library}
{BLAST}: High-order finite element hydrodynamics.
\newblock \url{https://computation.llnl.gov/projects/blast/icf-like-implosion}.

\bibitem{mfem-library}
{MFEM}: Modular finite element methods library.
\newblock \url{mfem.org}.

\bibitem{an2008optimizing}
Steven~S An, Theodore Kim, and Doug~L James.
\newblock Optimizing cubature for efficient integration of subspace
  deformations.
\newblock {\em ACM transactions on graphics (TOG)}, 27(5):165, 2008.

\bibitem{anderson2018high}
Robert~W Anderson, Veselin~A Dobrev, Tzanio~V Kolev, Robert~N Rieben, and
  Vladimir~Z Tomov.
\newblock High-order multi-material ale hydrodynamics.
\newblock {\em SIAM Journal on Scientific Computing}, 40(1):B32--B58, 2018.

\bibitem{astrid2008missing}
Patricia Astrid, Siep Weiland, Karen Willcox, and Ton Backx.
\newblock Missing point estimation in models described by proper orthogonal
  decomposition.
\newblock {\em IEEE Transactions on Automatic Control}, 53(10):2237--2251,
  2008.

\bibitem{barrault2004empirical}
Maxime Barrault, Yvon Maday, Ngoc~Cuong Nguyen, and Anthony~T Patera.
\newblock An ‘empirical interpolation’method: application to efficient
  reduced-basis discretization of partial differential equations.
\newblock {\em Comptes Rendus Mathematique}, 339(9):667--672, 2004.

\bibitem{berkooz1993proper}
Gal Berkooz, Philip Holmes, and John~L Lumley.
\newblock The proper orthogonal decomposition in the analysis of turbulent
  flows.
\newblock {\em Annual review of fluid mechanics}, 25(1):539--575, 1993.

\bibitem{buttari2008parallel}
Alfredo Buttari, Julien Langou, Jakub Kurzak, and Jack Dongarra.
\newblock Parallel tiled qr factorization for multicore architectures.
\newblock {\em Concurrency and Computation: Practice and Experience},
  20(13):1573--1590, 2008.

\bibitem{carlberg2017galerkin}
Kevin Carlberg, Matthew Barone, and Harbir Antil.
\newblock Galerkin v. least-squares petrov--galerkin projection in nonlinear
  model reduction.
\newblock {\em Journal of Computational Physics}, 330:693--734, 2017.

\bibitem{carlberg2011efficient}
Kevin Carlberg, Charbel Bou-Mosleh, and Charbel Farhat.
\newblock Efficient non-linear model reduction via a least-squares
  petrov--galerkin projection and compressive tensor approximations.
\newblock {\em International Journal for Numerical Methods in Engineering},
  86(2):155--181, 2011.

\bibitem{carlberg2013gnat}
Kevin Carlberg, Charbel Farhat, Julien Cortial, and David Amsallem.
\newblock The gnat method for nonlinear model reduction: effective
  implementation and application to computational fluid dynamics and turbulent
  flows.
\newblock {\em Journal of Computational Physics}, 242:623--647, 2013.

\bibitem{chaturantabut2010nonlinear}
Saifon Chaturantabut and Danny~C Sorensen.
\newblock Nonlinear model reduction via discrete empirical interpolation.
\newblock {\em SIAM Journal on Scientific Computing}, 32(5):2737--2764, 2010.

\bibitem{choi2019space}
Youngsoo Choi and Kevin Carlberg.
\newblock Space--time least-squares petrov--galerkin projection for nonlinear
  model reduction.
\newblock {\em SIAM Journal on Scientific Computing}, 41(1):A26--A58, 2019.

\bibitem{de2011blind}
Lieven De~Lathauwer.
\newblock Blind separation of exponential polynomials and the decomposition of
  a tensor in rank-(l\_r,l\_r,1) terms.
\newblock {\em SIAM Journal on Matrix Analysis and Applications},
  32(4):1451--1474, 2011.

\bibitem{drmac2016new}
Zlatko Drmac and Serkan Gugercin.
\newblock A new selection operator for the discrete empirical interpolation
  method---improved a priori error bound and extensions.
\newblock {\em SIAM Journal on Scientific Computing}, 38(2):A631--A648, 2016.

\bibitem{drmac2018discrete}
Zlatko Drmac and Arvind~Krishna Saibaba.
\newblock The discrete empirical interpolation method: Canonical structure and
  formulation in weighted inner product spaces.
\newblock {\em SIAM Journal on Matrix Analysis and Applications},
  39(3):1152--1180, 2018.

\bibitem{elmroth2000high}
Erik Elmroth and Fred Gustavson.
\newblock High-performance library software for qr factorization.
\newblock In {\em International Workshop on Applied Parallel Computing}, pages
  53--63. Springer, 2000.

\bibitem{everson1995karhunen}
Richard Everson and Lawrence Sirovich.
\newblock Karhunen--loeve procedure for gappy data.
\newblock {\em JOSA A}, 12(8):1657--1664, 1995.

\bibitem{farhat2014dimensional}
Charbel Farhat, Philip Avery, Todd Chapman, and Julien Cortial.
\newblock Dimensional reduction of nonlinear finite element dynamic models with
  finite rotations and energy-based mesh sampling and weighting for
  computational efficiency.
\newblock {\em International Journal for Numerical Methods in Engineering},
  98(9):625--662, 2014.

\bibitem{farhat2015structure}
Charbel Farhat, Todd Chapman, and Philip Avery.
\newblock Structure-preserving, stability, and accuracy properties of the
  energy-conserving sampling and weighting method for the hyper reduction of
  nonlinear finite element dynamic models.
\newblock {\em International Journal for Numerical Methods in Engineering},
  102(5):1077--1110, 2015.

\bibitem{grepl2007efficient}
Martin~A Grepl, Yvon Maday, Ngoc~C Nguyen, and Anthony~T Patera.
\newblock Efficient reduced-basis treatment of nonaffine and nonlinear partial
  differential equations.
\newblock {\em ESAIM: Mathematical Modelling and Numerical Analysis},
  41(3):575--605, 2007.

\bibitem{gu1996efficient}
Ming Gu and Stanley~C Eisenstat.
\newblock Efficient algorithms for computing a strong rank-revealing qr
  factorization.
\newblock {\em SIAM Journal on Scientific Computing}, 17(4):848--869, 1996.

\bibitem{hernandez2017dimensional}
Joaquin~Alberto Hernandez, Manuel~Alejandro Caicedo, and Alex Ferrer.
\newblock Dimensional hyper-reduction of nonlinear finite element models via
  empirical cubature.
\newblock {\em Computer methods in applied mechanics and engineering},
  313:687--722, 2017.

\bibitem{hinze2005proper}
Michael Hinze and Stefan Volkwein.
\newblock Proper orthogonal decomposition surrogate models for nonlinear
  dynamical systems: Error estimates and suboptimal control.
\newblock In {\em Dimension reduction of large-scale systems}, pages 261--306.
  Springer, 2005.

\bibitem{hotelling1933analysis}
Harold Hotelling.
\newblock Analysis of a complex of statistical variables into principal
  components.
\newblock {\em Journal of educational psychology}, 24(6):417, 1933.

\bibitem{khairallah2014mesoscopic}
Saad~A Khairallah and Andy Anderson.
\newblock Mesoscopic simulation model of selective laser melting of stainless
  steel powder.
\newblock {\em Journal of Materials Processing Technology}, 214(11):2627--2636,
  2014.

\bibitem{kunisch2002galerkin}
Karl Kunisch and Stefan Volkwein.
\newblock Galerkin proper orthogonal decomposition methods for a general
  equation in fluid dynamics.
\newblock {\em SIAM Journal on Numerical analysis}, 40(2):492--515, 2002.

\bibitem{loeve1955}
Michel Loeve.
\newblock {\em Probability Theory}.
\newblock D. Van Nostrand, New York, 1955.

\bibitem{maccormackNote}
Robert MacCormack.
\newblock Numerical computation of compressible viscous flow.
\newblock Technical report, Lecture notes for AA214b and AA214c, Stanford
  University, 2007.

\bibitem{nguyen2008efficient}
Ngoc~C Nguyen and Jaime Peraire.
\newblock An efficient reduced-order modeling approach for non-linear
  parametrized partial differential equations.
\newblock {\em International Journal for Numerical Methods in Engineering},
  76(1):27--55, 2008.

\bibitem{rewienski2003trajectory}
Micha{\l}~Jerzy Rewie{\'n}ski.
\newblock {\em A trajectory piecewise-linear approach to model order reduction
  of nonlinear dynamical systems}.
\newblock PhD thesis, Massachusetts Institute of Technology, Department of
  Electrical Engineering and Computer Science, 2003.

\bibitem{ryckelynck2005priori}
David Ryckelynck.
\newblock A priori hyperreduction method: an adaptive approach.
\newblock {\em Journal of computational physics}, 202(1):346--366, 2005.

\bibitem{sorber2013optimization}
Laurent Sorber, Marc Van~Barel, and Lieven De~Lathauwer.
\newblock Optimization-based algorithms for tensor decompositions: Canonical
  polyadic decomposition, decomposition in rank-(l\_r,l\_r,1) terms, and a new
  generalization.
\newblock {\em SIAM Journal on Optimization}, 23(2):695--720, 2013.

\bibitem{zahr2010comparison}
Matthew~J Zahr, K~Carlberg, D~Amsallem, and C~Farhat.
\newblock Comparison of model reduction techniques on high-fidelity linear and
  nonlinear electrical, mechanical, and biological systems.
\newblock {\em University of California, Berkeley}, 2010.

\end{thebibliography}

\end{document}